\documentclass[11pt,oldfontcommands,longfontname]{article}

\RequirePackage[OT1]{fontenc} \RequirePackage{amsthm,amsmath}

\usepackage{graphics,epsfig,epsf,psfrag}
\usepackage{caption}
\usepackage{subcaption}
\usepackage{url}
\usepackage{color, appendix}
\usepackage{graphicx}
\usepackage{amsmath, amsfonts, amsthm,amssymb}
\usepackage{epstopdf}
\usepackage{hyperref}
\usepackage{enumerate, framed}
\usepackage{graphicx,lscape,rotate}
\usepackage{algorithmic, algorithm}
\usepackage{epsfig}
\usepackage{mathrsfs}
\usepackage{amsfonts}
\usepackage{mathtools}
\usepackage{tasks}
\usepackage{ifthen}
\usepackage{bbm}
\graphicspath{{eps/}}
\usepackage{amsthm}
\usepackage[utf8]{inputenc}
\usepackage[english]{babel}
\usepackage{bm}
\usepackage[nobysame]{amsrefs}

\usepackage{scalerel,stackengine}
\stackMath

\makeatletter
\newcommand{\opnorm}{\@ifstar\@opnorms\@opnorm}
\newcommand{\@opnorms}[1]{%
  \left|\mkern-1.5mu\left|\mkern-1.5mu\left|
   #1
  \right|\mkern-1.5mu\right|\mkern-1.5mu\right|
}
\newcommand{\@opnorm}[2][]{%
  \mathopen{#1|\mkern-1.5mu#1|\mkern-1.5mu#1|}
  #2
  \mathclose{#1|\mkern-1.5mu#1|\mkern-1.5mu#1|}
}
\makeatother

\newcommand\reallywidehat[1]{%
\savestack{\tmpbox}{\stretchto{%
  \scaleto{%
    \scalerel*[\widthof{\ensuremath{#1}}]{\kern-.6pt\bigwedge\kern-.6pt}%
    {\rule[-\textheight/2]{1ex}{\textheight}}
  }{\textheight}%
}{0.5ex}}%
\stackon[1pt]{#1}{\tmpbox}%
}
\numberwithin{equation}{section}
\theoremstyle{plain}
\newtheorem{theorem}{Theorem}[section]
\newtheorem{assumption}{Assumption}[section]
\newtheorem{claim}{Claim}[section]
\newtheorem{remark}{Remark}
\newtheorem{proposition}{Proposition}[section]

\newtheorem{lemma}{Lemma}[section]
\newtheorem{definition}{Definition}
\newtheorem{example}[theorem]{Example}
\theoremstyle{remark}

\newcommand{\R}{ \mathbb R}

\newcommand{\N}{ \mathbb{N}}
\newcommand{\Pp}{ \mathscr{P}}

\newcommand{\Ps}{P_{sc}}
\newcommand{\s}{\mathcal{L}}
\newcommand{\tx}{\tilde{x}}
\newcommand{\A}{ \mathscr{A}}

\newcommand{\eq}{\begin{equation}}
\newcommand{\eeq}{\end{equation}}
\newcommand{\p}{\partial}

\newcommand{\I}{\mathscr{I}}
\newcommand{\Hi}{\mathcal{H}}

\begin{document}


\title{The three-quasi-particle scattering problem: asymptotic completeness for short-range systems  }
\author{Avy Soffer and Xiaoxu Wu}


\date{}  

\maketitle

\begin{abstract}
\noindent 
We develop an approach to scattering theory for generalized $N$-body systems.
In particular we consider a general class of three quasi-particle systems, for which we prove Asymptotic Completeness.
\end{abstract}
{
\hypersetup{linkcolor=black}
\tableofcontents}
\section{Introduction}\label{sec:intro}
In this note we study the scattering theory of a three-quasi-particle system. The system is described by the following equation:
\eq
\begin{cases}
i\partial_t\psi(x,t)=(H_0+V(x))\psi(x,t)\\
\psi(x,0)=\psi_0(x)\in \s^2_x(\mathbb{R}^9)
\end{cases},\quad x=(x_1,x_2,x_3)\in \mathbb{R}^{9}.\label{SE}
\eeq
Here, we define:
\begin{itemize}
    \item  $x_j\in \R^3 $(with $j=1,\cdots,3$) represents the position variable of the $j$th particle. 
    \item $P_j:=-i\nabla_{x_j}$ represents the wave number operator of the $j$th particle, often referred to as the quasi-momentum.
    \item $H_0=\sum\limits_{j=1}^3\omega_j(P_j)$, where $\omega_j: \R^3\to \R, \eta\mapsto \omega_j(\eta)$, denotes the kinetic energy operator, the dispersion relation.
    
    \item The solution $\psi(x,t)$ is a complex-valued function of $(x,t)\in \mathbb{R}^{9+1}$.
    \item The term $V(x)$, reflecting the interaction among three particles, has the form 
$$
V(x)=\sum\limits_{1\leq j<l\leq 3}V_{jl}(x_j-x_l),
$$
for some real-valued functions $V_{jl}, 1\leq j<l\leq 3$. Additionally, we assume that $V_{lj}(x_l-x_j)=V_{jl}(x_j-x_l)$ for $1\leq j<l\leq 3$.

\end{itemize}
We present two classical examples of the equation \eqref{SE}:
\begin{enumerate}
    \item When $\omega_j(P_j)=P_j^2=-\Delta_{x_j}, j=1,2,3$, it is the standard $3$-body system, which describes a system of $3$ non-relativistic particles interacting with each other.
    \item When $\omega_j(P_j)=\sqrt{m_j^2+P_j^2}$, $j=1,2,3,$ \eqref{SE} describes the system of 3 relativistic particles.
\end{enumerate}

There are several significant applications involving quasi-particles. For instance, a particle moving through a medium—like a periodic ionic crystal—exhibits a complex and implicit effective dispersion relation. This complexity is similarly observed in particles within Quantum Field Theory (QFT), where the effects of renormalization make the particle mass a momentum-dependent function in a complicated manner. There are also quasi-particles that are not elementary particles but rather derived constructs. A classic example involves the dynamics of the Heisenberg spin model, which can be described in terms of spin-wave excitations, also known as Magnons. In such cases, a typical dispersion relation can be expressed as
\eq
\omega(k)=4(3-\cos k_1-\cos k_2 -\cos k_3)
\eeq

Now let us introduce the interaction $V(x)$. Define $\langle x\rangle$ as $\sqrt{1+|x|^2}$ for $x\in \R^n$. When dealing with potentials or operators represented by functions, one typically assumes either of the following conditions:
\begin{itemize}
    \item (short-range potentials) For all $1\leq j<l\leq 3$, $|V_{jl}(\eta)|\leq C\langle \eta\rangle^{-1-\epsilon}$ for some $\epsilon>0$ and some constant $C>0$, or
    \item (long-range potentials) For some $1\leq j<l\leq3$, $|V_{jl}(\eta)|\leq C\langle \eta\rangle^{-\epsilon}$ only holds only for some $\epsilon \in (0,1]$.
\end{itemize}
In this note we focus on systems characterized by short-range interactions. System \eqref{SE} represents a multi-particle or three-body system. 

We study the scattering theory for system \eqref{SE}. This is the study of long-time behavior of solutions to system \eqref{SE}. The understanding of the large time behavior for complex multi-particle systems dates back to the early days of Quantum Mechanics. While there is a fairly good understanding of the $N$-body scattering problem, many other directions were left open. A primary objective in scattering theory is to establish asymptotic completeness (AC), which will be detailed later. For AC of the standard $N$-body systems, one can refer references such as \cites{HS1,SS1,SS1993, SS1990,SS1994,TH1993,D1993} and the cited works within them. For general $2$-body systems, please refer to \cite{S1990}. A consistent tool employed in the aforementioned studies is the Mourre estimate, a dispersive estimate. In these works,  this estimate is an integral part of the proof of AC. 

 Extending the Mourre estimate to cases where the particle dynamics are not described by the standard Schrödinger equation presents a challenge. Partial results have been obtained in this direction  over the years \cites{gerard1991mourre, Der1990, BreeS2019, zielinski1997asymptotic }. In \cite{Der1990} the Mourre estimate is proven for cases where the usual Dilation generator applies as a conjugate operator. This is established by introducing a structure assumption on the form of the commutator of the Dilation operator with all sub-Hamiltonians.  For a Mourre estimate of the relativistic case with uniform mass, where $\omega_j(P_j)=\sqrt{P_j^2+m^2}$ for $j=1,2,3,$ see \cite{gerard1991mourre}. For the AC of scattering states with negative energy, specifically when the Hamiltonian is defined as $H=P^2_x+|P_y|+V_{12}(x)+V_{13}(y)+V_{23}(x-y),$ where $x,y\in \R^3$, please refer to \cite{BreeS2019}. In \cite{zielinski1997asymptotic}
 the charge tranfer Hamiltonian for a quasi-particle dynamics is developed.

 In this work, we employ new ideas for deriving scattering, which have recently been developed in the context of general nonlinear dispersive and hyperbolic equations \cites{Liu-S,LS2021,SW20221,SW20222}. Importantly, our proof does not rely on either the Mourre estimate or dispersive estimates with respect to the full Hamiltonian. In \cite{Sof-W5}, the authors also developed a new method that utilizes asymptotic completeness (AC) and a compactness argument to establish local decay estimates.

%
\subsection{Problems}
 We investigate the scattering theory of system \eqref{SE}. Within scattering theory, the primary objective is to establish asymptotic completeness (AC). To be precise, let us define $\tilde{x}:=(x_1-x_2,x_2-x_3,x_3-x_1)\in \R^3$ and $H:=H_0+V(x)$. A state $\psi\in \s^2_x(\R^9)$ is termed a scattering state if, for all $M>0$,
\eq
t^{-1}\int_0^t ds\| \chi(|\tilde{x}|\leq M)\psi(x,s)\|^2_{\s^2_x(\mathbb{R}^9)}  \to0\quad \text{ in }\R\label{MRuelle}
\eeq
as $t\to \infty$. The collection of all such scattering states forms a linear space. We can characterize this space by proving Ruelle's Theorem. For an in depth discussion, refer to subsections \ref{2.2} and \ref{2.3}. Let $P_{sc}=P_{sc}(H)$ denote the projection on the space of all scattering states. AC implies that all scattering states will evolve asymptotically as a superposition of various simple solutions. Mathematically, this is represented as: if $\psi(x,0)=P_{sc}\psi(x,0)$, then
\eq
\| \psi(x,t)-\sum\limits_{a}\psi_{a,\pm}(x,t)\|_{\s^2_x(\R^9)}\to 0
\eeq
as $t\to \pm\infty,$ where $\psi_{a,\pm}(x,t)$ denotes some simple solution. For example, $\psi_{a,\pm}$ could be a free wave. This pertains to the study of the long-time behavior of solutions to system \eqref{SE}. 
\begin{definition}[Short-Range Potentials]
Assume that for all $1 \leq j < l \leq 3$, there exists a constant $C > 0$ and a positive value $\epsilon$ such that $|V_{jl}(\eta)| \leq C\langle \eta\rangle^{-1-\epsilon}$.
\end{definition}
In this note, we study the long-time behavior of solutions to \eqref{SE} in the presence of short-range potentials for a general class of $H$. Without relying on dispersive estimates, we will establish AC for all scattering states under certain assumptions, which we will detail subsequently.

\subsection{Preliminaries, assumptions and results}
Now let us introduce our assumptions and results. Consider the set 
\eq
L:=\{ (12)(3),(13)(2),(23)(1), (1)(2)(3), (123), 0\}.
\eeq
We assume $V(x)$ is of the form:
\eq
V(x)=\sum\limits_{a=(jk)(l)\in L}V_{jk}(x_j-x_k)\label{V3}
\eeq
subject to the condition:
\eq
\opnorm{V}:=\sum\limits_{a=(jk)(l)\in L}\| \langle \eta\rangle^{5} V_{jk}(\eta)\|_{\s^\infty_\eta(\R^3)}<\infty.
\eeq

\begin{figure}
     \centering 
     \begin{subfigure}[h]{0.3\textwidth}
         \centering
         \includegraphics[width=\textwidth]{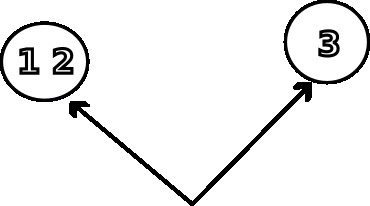}
         \caption{An example of a non-free scattering state}
         \label{p1}
     \end{subfigure}

     \begin{subfigure}[h]{0.3\textwidth}
         \centering
         \includegraphics[width=\textwidth]{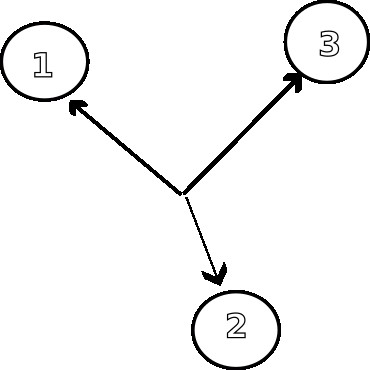}
         \caption{A free scattering state}
         \label{p2}
     \end{subfigure}
     \caption{Examples of three-body scattering state profiles}
\end{figure}

\begin{assumption}\label{asp: short}
Let \( V(x) \) be as defined in \eqref{V3} and ensure that \( \opnorm{V} < \infty \).
\end{assumption}
In this note, we consider $\omega_j(\eta)\in W^{4,\infty}_{loc}(\R^3)$ for $j=1,2,3$. The group velocity of the $j$th particle is defined as 
\eq
v_j(\eta):=\nabla_{\eta}[\omega_j(\eta)], \quad \text{ where }\eta\equiv P_j.
\eeq
The operator $P_j = -i\nabla_{x_j}$ (where $j=1,2,3$) represents the wave number operator, or quasi-momentum of the $j$th particle. A specific $\omega_j(\eta_0)$ is referred to as the threshold point for the $j$th particle when the group velocity of particle $j$ is equal to zero, i.e., $v_j(\eta)\vert_{\eta=\eta_0}=0$. When the energy of a free particle $j$ is close to its threshold point(s), its behavior may deviate from the usual behavior. In this note, we assume that for each $j=1,2,3$, $\omega_j(\eta)$ is smooth, and each particle has only a finite number of threshold points:

\begin{assumption}\label{asp: critical}
For every $j=1,2,3$, $\omega_j(\eta)$ adheres to the following:
\begin{itemize}
    \item[(a)]$\omega_j(\eta)\in W^{4,\infty}_{loc}(\R^3)$ .
    \item[(b)] The set $\{\eta: v_j(\eta)=0\}$ is finite.  
    
\end{itemize}
\end{assumption}
\noindent Examples that satisfy Assumption \ref{asp: critical} include:
\eq
\omega_j=\sqrt{P_j^2+m_j^2}, \quad m_j\geq 0, j=1,2,3,
\eeq
and 
\eq
\omega_j=|P_j|^4, j=1,2,3.
\eeq
Under Assumption \ref{asp: critical}, the set of all critical points of $\omega_j(\eta)$ (where $j=1,2,3)$ can be expressed as:
\begin{equation}
\mathscr{C}_j:=\{\eta: \nabla\eta[\omega_j(\eta)]=0\}=\{\eta_{j1}, \cdots, \eta_{jN_j}\}.\label{set: threshold}
\end{equation}
For a deeper understanding of group velocity, please consult Section \ref{section: curve}. We define 
\begin{align*}
V_a(x)&:=V_{jk}(x_j-x_k), \quad V(x)\equiv \sum_a V_a\\
H_{a}&:=H_0+V_a(x)
\end{align*}
for each $a=(jk)(l)\in L$. An in-depth exploration of $H_a$ can be found in Section \ref{subsection: Channels}. Additionally, for $a=(1)(2)(3)$, we define:
\eq
H_a:=H_0.
\eeq
Let's delve into the profiles of the asymptotic part, $\psi_{a,\pm}(t)$. Due to \eqref{MRuelle}, which defines a scattering state, it's clear that given sufficient time, the scattering solution will exit any region where $|\tilde{x}|$ is bounded. Under these circumstances, two scenarios emerge:
\begin{itemize}
    \item The distances between pairs of particles become significant. Refer to Figure \ref{p1}. To put it simply, in this scenario $H  \approx H_0 $.
    \item For some $1\leq j<k\leq 3$, particles $j$ and $k$ are cohesively bound, while the remaining particle drifts away from this particle cluster. Refer to Figure \ref{p2} for the case when $j=1$ and $k=2$. Here, in a broad sense, $H  \approx H_a $.
\end{itemize}
In this note, we focus on the case where $t>0$. Building upon our discussion, we can rephrase the statement of AC as follows: for all scattering states of system \eqref{SE}, the corresponding solutions adhere to:
\eq
\| \psi(t)-\sum\limits_{a\in L_{sc}}e^{-itH_a}\psi_{a}(x)\|_{\s^2_x(\R^n)}\to 0\quad \text{ as }t\to \infty,\label{AC}
\eeq
where $L_{sc}:=L-\{(123), 0\}$ and for each $a\in L_{sc},$ $\psi_a$ denotes the corresponding asymptotic state. 

We will now introduce the concept of channel wave operators, which are essential for capturing all asymptotic states. In the three-body system described in equation \eqref{SE}, a scattering state can either evolve as a free wave or follow a non-free scattering pattern, as depicted in Figures \ref{p1} and \ref{p2}. Thus, this is a multi-channel problem, and we require multi-channel scattering theory. The key elements in this theory are the multi-channel wave operators, which play a crucial role in capturing all asymptotic states.

\begin{definition}[Channel Wave Operators] The channel wave operator for the $H_a$ channel is defined as:
\begin{equation}
\Omega_a^*\psi(0) = s\text{-}\lim_{t\to \infty} e^{itH_a}J_ae^{-itH}P_{sc}\psi(0),
\end{equation}
where $\psi(0)\in\mathcal{H}=\s^2_x(\mathbb{R}^9)$, and $P_{sc}$ represents the projection onto the space of all scattering states. The set ${ H_a}$ refers to the collection of all sub-Hamiltonians of $H$, and ${J_a}$ is a set of smooth cut-off functions of $x$, $P$, and $t$ satisfying:
\begin{equation}
\sum_{a} J_a = 1.
\end{equation}
\end{definition}

Our proof of the existence of all channel wave operators relies on the introduction of $J_a$, Assumption \ref{asp: short}, and the dispersive estimates of the free flow of each particle. We will now detail the assumption associated with the dispersive estimate of the free flow for each particle. To be precise, we introduce $F(\lambda)$, $F_c(\lambda)$, and $F_1(\lambda)$ as smooth cutoff functions over the interval $[1,+\infty)$. We then define:
\begin{equation}
F(\lambda\leq a):=1-F(\lambda/a), \quad F_c(\lambda\leq a):=1-F_c(\lambda/a), \quad F_1(\lambda>a):=F_1(\lambda/a), 
\end{equation}
and
\begin{equation}
F_c(\lambda\leq a):=F_c(\lambda/a\leq 1):=1-F_c(\lambda/a), \quad F_1(\lambda\leq a):=F_1(\lambda/a\leq 1):=1-F_1(\lambda/a).
\end{equation}
Furthermore, it's posited that:
\eq
\begin{cases}
F_c(k)=F_1(k)=F(k)=1& \text{ when } k\geq 1\\
F_c(k)=F_1(k)=F(k)=0& \text{ when } k\leq 1/2
\end{cases}.\label{supF}
\eeq 
For operators $A_j$, where $j=1,\cdots,n$, we present the product $\Pi_{j=1}^{j=n}A_j$ as: 
\eq
\Pi_{j=1}^{j=n}A_j=A_1\times A_2\times \cdots \times A_n
\eeq
for operators $A_j$, $j=1,\cdots,n$. Throughout this paper, we maintain the assumption that for $j=1,2,3,$ $\omega_j(P_j)$ adheres to the following:
\begin{assumption}\label{asp: omega}Given $\eta_{jk}(k=1,\cdots, N_j)$ as in \eqref{set: threshold}, it is postulated that for every $j=1,2,3,$ an $\alpha>0$ exists such that for all $l\in \{1,2,3\}-\{j\}$,
\eq
E_{jl}:=\int_1^\infty dt\|F_{c,j,\alpha}(x_j,t,P_j)e^{it\omega_j(P_j)}\langle x_j-x_l\rangle^{-5}\|_{\s^2_x(\R^9)\to \s^2_x(\R^9)} <\infty,
\eeq
with
\eq
F_{c,j,\alpha}(x_j,t,P_j):=F_c(|x_j|\leq t^\alpha)\left(\Pi_{k=1}^{N_j}F_1(|P_j-\eta_{jk}|>\frac{1}{t^{\alpha/2}})\right)F_1(|P_j|\leq t^{\alpha/2}).\label{Fcjalpha}
\eeq
 We define $E$ as:
\eq
E:=\max\limits_{j,l\in \{1,2,3\},j\neq l} E_{jl},
\eeq
and  setting 
\eq
J_a=e^{-it\omega_j(P_j)}F_{c,j,\alpha}(x_j,t,P_j)e^{it\omega_j(P_j)}.
\eeq
\end{assumption}
\begin{remark}
    Indeed, Assumption \ref{asp: omega} can be derived using the non-stationary phase method, coupled with the fact that there exists only a finite number of vectors $\eta$ for which $\nabla_\eta[\omega_j(\eta)]=0$ $($where $j=1,2,3).$ 
\end{remark}
We prove the existence of all channel wave operators under Assumptions \ref{asp: short}, \ref{asp: critical}, and \ref{asp: omega}. Specifically, for every $a=(jk)(l)\in L$, and for all $\alpha$ that meet Assumption \ref{asp: omega}, the following holds:
\eq
\Omega_{a,\alpha}^{*}:=s\text{-}\lim\limits_{t\to \infty} F_{c,l,\alpha}(x_l,t,P_l)e^{itH_a}e^{-itH},\quad \text{ on }\s^2_x(\R^9)
\eeq
exists and the new free channel wave operator
\eq
\Omega_{free,\alpha}^{*}:=s\text{-}\lim\limits_{t\to \infty} \left(\Pi_{l=1}^{l=3}F_{c,l,\alpha}(x_l,t,P_l)\right)e^{itH_0}e^{-itH},\quad \text{ on }\s^2_x(\R^9)
\eeq
exists as well, where $F_{c,l,\alpha}(x_l,t,P_l)$ is defined in \eqref{Fcjalpha}. We state this result formally in the following theorem:
\begin{theorem}\label{thm: channel}If Assumptions \ref{asp: short}, \ref{asp: critical} and \ref{asp: omega} are satisfied, then  $\Omega_{a,\alpha}^{*}$ exists for all $a=(jk)(l)\in L$, and $\Omega_{ free,\alpha}^*$ exists as well.  
\end{theorem}
The proof of Theorem \ref{thm: channel} is deferred to Section \ref{Channel}. According to Theorem \ref{thm: channel}, for each $a=(jk)(l)\in L_{sc}$, $\psi_a$ is given by 
\eq
\psi_a=\Omega_{a,\alpha}^{*}\psi,\label{psia}
\eeq
and when $a=(1)(2)(3)$,
\eq
\psi_{a}=\Omega_{free,\alpha}^*\psi. \label{psifree}
\eeq
A key element in proving Theorem \ref{thm: channel} is the application of the propagation estimates from \cite{SW1}. Here, we point out that if both $\alpha$ and $\alpha'$ satisfy Assumption \ref{asp: omega}, then for all $a=(jk)(l)\in L$, we have
\begin{align*}
\Omega_{a,\alpha}^{*}&=\Omega_{a,\alpha'}^{*},\\
\Omega_{free,\alpha}^{*}&=\Omega_{free,\alpha'}^{*}.
\end{align*}
In other word, both $\Omega_{a,\alpha}^{*}$ and $\Omega_{free,\alpha}^*$ are independent on the specific choice of $\alpha$. Instead of using $\Omega_{a,\alpha}^{*}$ and $\Omega_{free,\alpha}^*$, we employ $\Omega_{a}^{*}$ and $\Omega_{free}^*$, respectively. A detailed discussion can be found in Lemma \ref{Lem: weak}. Based on the concept of channel wave operators, we introduce the channel projections as follows. For each $a=(jk)(l)\in L$,
\eq
P_{a}=s\text{-}\lim\limits_{t\to \infty} e^{itH}e^{-itH_a}F_{c,l,\alpha}(x_l,t,P_l)e^{itH_a}e^{-itH},\quad \text{ on }\s^2_x(\R^9),\label{channel Pa}
\eeq
and
\eq
P_{free}=s\text{-}\lim\limits_{t\to \infty} e^{itH}e^{-itH_0}\left(\Pi_{l=1}^{l=3}F_{c,l,\alpha}(x_l,t,P_l)\right)e^{itH_0}e^{-itH},\quad \text{ on }\s^2_x(\R^9).\label{channel Pfree}
\eeq
Here, we omit the subscript $\alpha$ in both \eqref{channel Pa} and \eqref{channel Pfree} since the channel projections are independent of $\alpha$ as well. See Lemma \ref{Lem: Pab: weak}. The existence of $P_a$ and $P_{free}$ follow by using Cook's method and propagation estimates introduced in \cite{SW20221}. Please refer to the proof of Proposition \ref{cor 1} in Section \ref{sec: channel projection} for a detailed discussion. We would also like to remind the reader that, based on the definition of $F_{c,j,\alpha}(x_j,t,P_j)$ in \eqref{Fcjalpha}, the following identity holds for all $a=(jk)(l)\in L$: 
\eq
e^{-itH_a}(1-F_{c,l,\alpha}(x_l,t,P_l))e^{itH_a}=e^{-it\omega_l(P_l)}(1-F_{c,l,\alpha}(x_l,t,P_l))e^{it\omega_l(P_l)}.\label{intro: id1}
\eeq 
By employing \eqref{intro: id1}, we can conclude that the existence of all channel projections, mentioned in \eqref{channel Pa} and \eqref{channel Pfree}, implies the existence of $P_\mu$:
\eq
P_\mu:= s\text{-}\lim\limits_{t\to \infty} e^{itH}\left(\Pi_{l=1}^{l=3} e^{-it\omega_l(P_l)}(1-F_{c,l,\alpha}(x_l,t,P_l))e^{it\omega_l(P_l)}\right)e^{-itH}\quad \text{ on }\s^2_x(\mathbb{R}^9),
\eeq
where $\alpha$ satisfies Assumption \ref{asp: omega}, see Propositions \ref{cor 1} and \ref{Prop: Pmu: exist}.

 According to the statement in \eqref{AC} and the formulations for $\psi_a$ given by \eqref{psia} and \eqref{psifree}, we can infer that if $P_\mu P_{sc}=0$, then AC holds true.

\begin{theorem}
If $P_\mu P_{sc}=0$ and if Assumptions \ref{asp: short}, \ref{asp: critical}, and \ref{asp: omega} are satisfied, then AC holds for system \eqref{SE}.
\end{theorem}
Therefore, the challenge of proving AC is reduced to demonstrating that $P_\mu P_{sc}=0$. Before delving into the validity of $P_\mu P_{sc}=0$, let's take a closer look at $P_{sc}.$ Once we establish Ruelle's Theorem for system \eqref{SE}, we can provide a precise description of the space of all scattering states. This space is a subspace of $\s^2_x(\mathbb{R}^9)$ and corresponds to the continuous spectrum of $H$ on each fiber. Its complement in $\s^2_x(\mathbb{R}^9)$ is referred to as the space of all bound states or eigenfunctions. More specifically, the projection onto the space of all scattering states is given by the fiber integral $P_{sc}: \s^2_x(\mathbb{R}^9)\to \s^2_x(\mathbb{R}^9)$:

\begin{equation}
\Ps f = P_c(P_1+P_2+P_3)f,
\end{equation}
where $P_c(Q)$ denotes the projection onto the space of the continuous spectrum of $H$ with a total wave number operator $\sum\limits_{j=1}^3 P_j=Q$ (Here, we use the fact that $[P_1+P_2+P_3,H]=0$, which implies that the total wave number operator $\sum\limits_{j=1}^3P_j$ is invariant under the evolution of system \eqref{SE}). For detailed information, please see Section \ref{2.2}. We refer to $\psi_0(x)\in \s^2_x(\mathbb{R}^9)$ as a bound state of $H$ if it is orthogonal to all scattering states. Here, the space of all bound states, the complement of the space of all scattering states in $\s^2_x(\mathbb{R}^9)$, is equal to the discrete spectrum of $H$ on each fiber. In quantum mechanics, an eigenfunction is commonly known as a bound state. For more details, please see Proposition \ref{Rwhole} in Section \ref{2.3}.

Now, let's discuss the validity of $P_\mu P_{sc}=0$. Interestingly, under further assumptions on the two-body sub-Hamiltonians, $P_\mu$ acts as the projection onto the space of all bound states. This claim is formally presented in Proposition \ref{Prop: Pmu}, with its proof deferred to Sections \ref{sec 3} through \ref{section: proof free}. Before presenting Proposition \ref{Prop: Pmu}, let's first discuss the assumptions outlined within it. To begin with, the discussion about the assumptions requires the Feynman-Hellmann Theorem:

\begin{theorem}[Feynman-Hellmann Theorem]\label{Thm: FH0}For $ I=(a_1,b_1)\times \cdots (a_m,b_m)\subseteq \R^m$, let $\{H(p)\}_{p\in I}$ be a class of self-adjoint operators on $\s^2_x(\R^n)$. {Assume that $\lambda(p)$ is an eigenvalue of $H(p)$. If $H(p)$ is differentiable with respect to $p$,} then $\lambda(p)$ is differentiable with respect to $p$ at $p=p_0$:
\eq
\frac{\partial\lambda(p)}{\partial p_k}\vert_{p=p_0}=(\psi_\lambda, \frac{\partial H(p)}{\partial p_k}\psi_{\lambda})_{\s^2_x(\R^n)}\vert_{p=p_0}, \quad k=1,\cdots, m\label{con:FH0}
\eeq
where $\psi_\lambda$ stands for a normalized eigenfunction $(\|\psi_\lambda\|_{\s^2_x(\R^n)}=1)$ with an eigenvalue $\lambda(p).$

\end{theorem}
It's clear that the Feynman-Hellmann Theorem cannot be universally applied. Consider, for example, a family of parameterized operators denoted as $\{H(p)\}_{p\in I}$, where $H(p) = -\Delta_x + pW(x)$ and $p\in I=(-2,2)$. In this case, given $W(x)\leq 0$ to ensure that $-\Delta_x + W(x)$ has an eigenvalue $\lambda$ in $\s^2_x(\mathbb{R}^n)$, it follows that $H(1)$ has an eigenvalue $\lambda(1):=\lambda$. However, $H(0)=-\Delta_x$ does not have an eigenvalue in $\s^2_x(\mathbb{R}^n)$. This implies that the eigenvalue $\lambda(p)$ disappears at some $p\in (0,1)$. Thus, in such scenarios, the Feynman-Hellmann Theorem isn't applicable across the entire index set $(-2,2)$. If the eigenvalue does not disappear at $p=p_0$, then we can apply the Feynman-Hellmann Theorem, and we say that the Feynman-Hellmann Theorem is applicable to this eigenvalue at $p=p_0$. In this note, we state that the Feynman-Hellmann Theorem is \textbf{applicable} to $\{ H(p)\}_{p\in I}$ at $p=p_0$ if all eigenvalues of $H(p)$ are defined in a neighborhood of $p_0$ in $\mathbb{R}^m.$

Let's now discuss the conditions under which an eigenvalue might disappear. To begin, let's consider non-embedded eigenvalues. We refer to $\lambda$ as a non-embedded eigenvalue of a one-body self-adjoint Hamiltonian $H_1 = \omega(P) + W(x)$ if $\lambda$ is an eigenvalue of $H_1$ and if
\eq
\lambda < \min\limits_{q\in \R^3}\omega(q). 
\eeq
Here, we refer to $\min\limits_{q\in \R^3} \omega(q)$ as the threshold energy of $\omega(P)$. Remarkably, we've found that the non-embedded eigenvalues of our two-body Hamiltonians, which can be simplified to one-body forms, vanish only as they near the threshold energy of $\omega(P)$. This is provided that both the position variable $x_j$ of the $j$th particle and the group velocity remains controlled when its energy is finite.
\begin{assumption}
    \label{asp: critical b}We assume that:
   
    \begin{itemize}
     \item[(a)] $\omega_j(\eta)\to \infty$ as $|\eta|\to \infty$. 
        \item[(b)] For any $E\in \R$, we assume that for all $t\in \R$ and $d=1,2,3,4,$
    \eq
    \| \langle x_k\rangle^d e^{itH^a(\eta)}F(H^a(\eta)< 4E)\langle x_k\rangle^{-d}\|_{\s^2_{x_k}(\R^3)\to \s^2_{x_k}(\R^3)}\lesssim_E \langle t\rangle^d 
    \eeq
    and 
    \eq
\sup\limits_{s\in \R}\| \nabla_\eta[\omega_j(\eta-P_k)]e^{isH^a(\eta)} F(H^a(\eta)<E)\|_{\s^2_{x_k}(\R^3)\to \s^2_{x_k}(\R^3)}\lesssim_{E}1
\eeq
    hold true.
    \end{itemize}
\end{assumption}
\begin{remark}We define $\psi\in \s^2$ to have a finite energy if $\psi=F(H\leq M)\psi$ for some $M\geq1$. In Assumption \ref{asp: critical b}, we assume that for $j=1,2,3,$ $\omega_j(\eta)\to \infty$ as $|\eta|\to\infty$. This is because we aim to establish that $P_\mu P_{sc}=0$ on $\s^2_x(\R^9)$ by demonstrating $P_\mu P_{sc}\psi=0 $ for every $\psi\in \s^2_x(\R^9)$ with a finite energy. 
    
\end{remark}
\begin{lemma}\label{Lem: dis}
If Assumptions \ref{asp: short} and \ref{asp: critical b} are valid, then all non-embedded eigenvalues of $H^a(\eta)=\omega_j(\eta-P_k)+\omega_k(P_k)+V_{jk}(-x_k)$, for $a=(jk)(l)\in L$ and $\eta\in \R^3$, disappear at threshold points of $H^a(\eta)$.
\end{lemma}

We defer the proof of Lemma \ref{Lem: dis} to subsection \ref{sec: FH}. We refer to $E\in \R$ as a threshold point of the one-body Hamiltonian $H_1$ if $E=\omega(\eta_0)$ with $\nabla_\eta[\omega(\eta)]\vert_{\eta=\eta_0}=0.$ Then the threshold energy of $\omega(P)$ is a threshold point of $H_1$. We know that a two-body Hamiltonian can be simplified to a one-body Hamiltonian through a coordinate transformation. Assuming the eigenvalues of all two-body sub-Hamiltonians in our three-body Hamiltonian $H$ do not reach threshold points, the embedded eigenvalues of these two-body sub-Hamiltonians will not disappear. This assumption, that all threshold points of a two-body Hamiltonian are neither an eigenvalue nor a resonance, is also termed the regularity assumption for the two-body Hamiltonian. Here, a resonance of a one-body Hamiltonian $H_1$ is a distributional solution of $H_1\psi=0$ such that $\psi\notin \s^2_x$ but $\langle x\rangle^{-\sigma}\psi \in \s^2_x$ for some $\sigma >0$. Refer to \cite{JK1979} for one-body Schr\"odinger operators.

\begin{assumption}\label{asp: Re} All threshold points are neither an eigenvalue nor a resonance of $H^a(\eta):=\omega_j(\eta-P_k)+\omega_k(P_k)+V_{jk}(-x_k)$ for all $\eta\in \R^3$ and $a=(jk)(l)\in L$. 
    
\end{assumption}
Under Assumption \ref{asp: Re}, due to Lemma \ref{Lem: dis}, we conclude that all non-embedded eigenvalues of  $H^a(\eta)$, for $\eta\in \R^3, a=(jk)(l)\in L$, will not disappear and Feynman-Hellmann Theorem is then \textbf{applicable}. Consequently, we can express the projection onto the space of all eigenfunctions with non-embedded eigenvalues of $H^a(\eta)$ as 
\eq
P_b(H^a(\eta))=\sum\limits_{d=1}^{\tilde{N}^a} \tilde{P}_{b,d}(H^a(\eta))
\eeq
where $\tilde{N}^a=\tilde{N}^a(a)\in \N$ is independent on $\eta\in \R^3$. Here, for $d=1,\cdots,\tilde{N}^a$, $\tilde{P}_{b,d}(H^a(\eta))$ denotes rank-one projection with an eigenvalue $\lambda_{a,d}(\eta)$. We say that for $1\leq d_1<d_2\leq \tilde{N}^a$ $\lambda_{a,d1}\neq \lambda_{a,d2}$ if $\lambda_{a,d_1}(\eta)\neq \lambda_{a,d_2}(\eta)$ for some $\eta\in \R^3$. Otherwise, we say that $\lambda_{a,d_1}$ and $\lambda_{a,d_2}$ are the same. By putting the rank-one projections $\tilde{P}_{b,d}(H^a(\eta))$ with the same $\lambda_{a,d} $(same for all $\eta\in \R^3$) into one projection, we get
\eq
P_b(H^a(\eta))=\sum\limits_{d=1}^{N^a} P_{b,d}(H^a(\eta))\label{def: Na}
\eeq
where $N^a=N^a(a)\in \N$ denotes the number of distinct $\lambda_{a,d}$. 
\begin{remark}Here, we have both $\tilde{N}^a$ and $N^a$ are finite for all $a=(jk)(l)\in L$. This arises due to \eqref{asp:Pbd} as presented in Assumption \ref{asp: known}, which will be detailed later. Put simply, when all eigenfunctions are both smooth and spatially localized as suggested by \eqref{asp: known}, the dimension of the space containing all eigenfunctions of $H^a(\eta)$ remains finite for all $\eta\in \R^3, a=(jk)(l)\in L$. This is because $\langle x\rangle^{-4}\langle P\rangle^{-4}$ is a compact operator on $\s^2_x(\R^n)$ for all $n=1,2,\cdots.$ It's also noteworthy that for a one-body Schr\"odinger operator, $H_1=-\Delta_x+W(x)$, given some constant $b>0$, if the potential $W(x)\geq -\frac{b}{4}|x|^{-2}$ for $|x|\geq R$ with some sufficiently large $R\geq 1$, then $H_1$ has at most finitely many eigenvalues. See Theorem XIII.6 on page 87 in \cite{RS1978}.
    
\end{remark}
For simplicity, we assume that for all $\eta\in \R^3$, $a=(jk)(l)\in L$, $d_1, d_2\in \{1,2,3\}$ and $ l_1,l_2\in \{0,1\}$, we have
\eq
\mathscr{C}_{a,b}:=\sum\limits_{l_1=0}^1\sum\limits_{l_2=0}^1\| \p_{\eta_{d_1}}^{l_1}\p_{\eta_{d_2}}^{l_2}[ P_{b,d}(H^a(\eta))]\|_{\s^2_{x_k}(\R^3)\to \s^2_{x_k}(\R^3)}\lesssim 1.\label{Pb: eq}
\eeq
\begin{assumption}\label{asp: Pb}For all $\eta\in \R^3$, $a=(jk)(l)\in L$, $d_1, d_2\in \{1,2,3\}$ and $ l_1,l_2\in \{0,1\}$, \eqref{Pb: eq} is valid. 
    
\end{assumption}
Now let us discuss the embedded eigenvalues of $H^a$ for $a=(jk)(l)\in L$. In the context of a one-body Hamiltonian, the Fermi's golden rule (as discussed in \cite{SW1998}, for instance) indicates that these embedded eigenvalues are unstable. This means that adding a small perturbation, the embedded eigenvalues go away. Let
\eq
I_{a,em}:=\{ \eta\in \R^3: H^a(\eta)\text{ has an embedded eigenvalue}\}.
\eeq
For simplicity, we assume that $I_{a,em}$ is contained in the union of fintely many hypersurfaces of $\R^3$:
\begin{assumption}\label{asp: emb}For all $a=(jk)(l)\in L$, $I_{a,em}$ is contained in the union of fintely many hypersurfaces of $\R^3$.
    
\end{assumption}
Given Assumption \ref{asp: emb}, we infer that for all $a=(jk)(l)\in L$, $P_b(H^a)$ can be viewed as the projection onto the space of all bound states of $H^a$ in $\s^2_x(\R^9)$.

We now turn our attention to the assumptions made on all two-body sub-Hamiltonians. We begin with the introduction of cut-off functions of operators. Let us consider a finite set $A=\{a_1,\cdots,a_n\}\subseteq \R$. For any given $y\in \R$, we define 
\eq
[y-A]_d:=(y-a_1,\cdots, y-a_n).
\eeq
Given any $c\in (0,1)$, we introduce the following notation
\eq
\bar{F}([y-A]_d,c):=\left( \Pi_{j=1}^{j=n} F(|y-a_j|> c)\right)\times F(|y|\leq\frac{1}{c}).
\eeq
When $A= \emptyset$, we adopt
\eq
\bar{F}([y-A]_d,c):=F(|y|\leq\frac{1}{c} ).
\eeq
Here, $\bar{F}$ means "away". In other word, in the support of $\bar{F}([y-A]_d,c)$, the value of $y$ is bounded and it is distinctly separated from the set $A$ by a distance greater than $c/2$. Let's also consider a family of finite sets in $\R$, represented as $\{ S(p)\}_{p\in \R^3}$. When $n\geq2$, we define $H_n$ as a quasi-$n$-particle Hamiltonian:
$H_n=\sum\limits_{j=1}^n\omega_j(P_j)+V(x)$ where for $j=1,\cdots,n$, $x_j\in \R^3$, $P_j:=-i\nabla_{x_j}$ and 
\eq
V(x)=\sum\limits_{1\leq j<l\leq n}V_{jl}(x_j-x_l).
\eeq
For a given $c>0$ and any quasi-$n$-particle Hamiltonian $H_n$(where $n=2,3$), we define 
\eq
\bar{F}_{S,n}(H_n,c):=\bar{F}([H_n, S(\sum\limits_{j=1}^n P_j)]_d, c)F(|\sum\limits_{j=1}^nP_j|\leq \frac{1}{c}).
\eeq
Moreover, we define $F_{S,n}:=1-\bar{F}_{S,n}$ for $n=2,3$. To express the asymptotic relative velocity between particle $j$ and particle $k$, we introduce the notion $v_a$, given $a=(jk)(l)\in L$. This is defined by
\eq
v_{a}:=\Omega^a(v_j(P_j)-v_k(P_k))\Omega^{a,*},
\eeq 
where $\Omega^a:=s\text{-}\lim\limits_{t\to \infty} e^{itH_a}e^{-itH_0}$ and $\Omega^{a,*} $ stands for the conjugate of $\Omega^a$. For $a=(jk)(l)\in L$ and $\eta\in \R^3$, let $P_c(H^a(\eta))$ denote the projection on the space of the continuous spectral subspcae of $H^a(\eta)$. We define $P_c(H^a)$ as $P_c(H^a(P_j+P_k))$. Roughly speaking, our Assumption \ref{asp: subH} is that when we fix $\sum\limits_{j=1}^3P_j$, except for a finite number of points, the magnitude of $|v_a|$ is approximately $ | P_j-P_k|$ when $P_j$ is close to $P_k$. Furthermore, we assume that both $e^{-itH_a}P_c(H^a)$ and $e^{-itH_0}$ adhere to local decay estimates in a manner similar to how $e^{-it(-\Delta)}$ does.
\begin{assumption}\label{asp: subH}For all $a=(jk)(l)\in L$, we assume that there exists a collection of finite sets in $\R$ denoted by $\{ I_a(p)\}_{p\in \R^3}$, with $\sup\limits_{p\in \R^3} |I_a(p)|<\infty,$ such that for all $\beta>0$ and some $C=C(c)>0$
\eq
\|  \langle x_j-x_k\rangle^{-2} F(|v_{a}|\leq \langle t\rangle^{-\beta})\bar{F}_{I_a,3}(H_a,c)\|_{\s^2_x(\R^9)\to \s^2_x(\R^9)}\leq \frac{C}{\langle t\rangle^{3\beta/2}},\label{asp: subH: eq1}
\eeq
    \eq
    \| \langle x_k-x_j\rangle^{-2} e^{-itH_a}\bar{F}_{I_a,3}(H_a,c)P_c(H^a)\langle x_k-x_j\rangle^{-2 }\|_{\s^2_{x}(\R^9)\to \s^2_{x}(\R^9)}\leq \frac{C}{\langle t\rangle^{3/2} }\label{local decay 2body}
    \eeq
and
\eq
    \| \langle x_j-x_k\rangle^{-2} e^{-itH_0}\bar{F}_{I_a,3}(H_0, c) \|_{\s^2_x(\R^9)\to \s^2_{x,t}(\R^{9+1})}\leq C.\label{free: smoothing}
    \eeq

\end{assumption}
\begin{remark}Assumption \ref{asp: subH} indicates that when we fix the $\sum\limits_{j=1}^3P_j$, aside from a finite number of energy values, all two-body sub-Hamiltonians and the free flow exhibit normal decay estimates similar to those seen with Schr\"odinger operators. For Schr\"odinger operators, an estimate like \eqref{free: smoothing} follows from the local resolvent estimate: for some $C=C(c)>0$
\eq
\| \langle x_j-x_k\rangle^{-2} \frac{1}{H_0-y\pm i0} \bar{F}_{I_a,3}(H_0,c)\|_{\s^2_x(\R^9)\to \s^2_{x,y}(\R^{9+1})}\leq C.
\eeq
See \cite{CS1988}. 
    
\end{remark}
Using a pointwise in time local decay estimate, as given in \eqref{local decay 2body}, we can obtain further results, including resolvent estimates and minimum velocity bounds for Schr\"odinger operators. Refer to \cite{HSS1999} for details on the minimum velocity bounds for Schr\"odinger operators. For a one-body Schr\"odinger operator, it is also well-known that if $0$ is neither an eigenvalue nor a resonance and if the potential $W(x)$ is well-localized (e.g. $|W(x)|\leq \frac{C}{\langle x\rangle^{2+\epsilon}}$ for some $\epsilon>0$ and constant $C>0$), the dimension of the space of all eigenfunctions is finite and all eigenfunctions are localized in space variable $x$.  These outcomes can be analogously derived for the general one-body Hamiltonian. For the sake of clarity in this note, we choose to assume these results, which can be validated using known techniques, in Assumption \ref{asp: known}. We defer their proofs to a subsequent paper. We assume that the set of all threshold points of $H^a(\eta)$(where $\eta\in\R^3$ and $a=(jk)(l)\in L$) are finite uniformly in $\eta\in \R^3$ as well. 
\begin{assumption}\label{asp: Hthreshold}For $a=(jk)(l)\in L$, let 
\eq
I_{aE}(\eta):=\{E: E=\omega_j(\eta-q_k)+\omega_k(q_k), \nabla_{q_k}[\omega_j(\eta-q_k)+\omega_k(q_k)]=0 \},\quad \text{for }\eta\in \R^3.
\eeq
We postulate the existence of a positive integer $N$ such that $|I_{aE}|\leq N<\infty$.
    
\end{assumption}
\begin{assumption}\label{asp: known}For $a=(jk)(l)\in L$, let $I_a, I_{aE}$ be as in Assumptions \ref{asp: subH} and \ref{asp: Hthreshold}, respectively. 
We further assume the existence of constants $C_1=C_1(c,\epsilon)>0$ and $C_2=C_2(c)>0$ such that
    \eq
    \| \langle x_k-x_j\rangle^{-2} e^{-itH_a}\bar{F}_{S,3}(H_a,c)P_c(H^a)F(|v_a|> \frac{1}{\langle t\rangle^{1/2-\epsilon}})\langle x_k-x_j\rangle^{-2 }\|_{\s^2_{x}(\R^9)\to \s^2_{x}(\R^9)}\leq \frac{C_1}{\langle t\rangle^{3/2} },
    \eeq
     \eq
    \| \langle x_j-x_k\rangle^{-2} e^{-itH_a}\bar{F}_{S,3}(H_a, c) \|_{\s^2_x(\R^9)\to \s^2_{x,t}(\R^{9+1})}\leq C_2
    \eeq
     and
    \eq
    \| \langle P_k\rangle^4\langle x_k\rangle^4 P_{b,d}(H^a(\eta)) \|_{\s^2_{x_k}(\R^3)\to \s^2_{x_k}(\R^3)}\lesssim 1,\label{asp:Pbd}
    \eeq
    where $d=1,\cdots,N^a$. Additionally, we assume that for all values $c>0$, $ \sigma \in (0,4]$, $\epsilon\in (0,1/2)$ and $M\geq 1$, 
    \eq
    \| \langle x_j-x_k\rangle^{-\sigma} \bar{F}_{I_{aE},2}(H^a,c)e^{-itH_a}P_c(H^a)\langle x_j-x_k\rangle^{-\sigma}\|_{\s^2_{x}(\R^9)\to \s^2_x(\R^9)}\leq \frac{C}{\langle t\rangle^\sigma},
    \eeq
   and 
    \begin{align}
   & \| F(|x_l-x_j|\geq M)\langle x_j-x_k\rangle^{-\sigma} F(|v_a|>\frac{1}{\langle t\rangle^{1/2-\epsilon}})F_{S,2}(H^a,c)e^{-itH_a}\nonumber\\
   &\times \langle x_j-x_l\rangle^{-\sigma}\chi(|x_j-x_k|\leq \langle t\rangle^{1/2})\|_{\s^2_{x}(\R^9)\to \s^2_x(\R^9)}\nonumber\\
   &\leq \frac{C}{\langle t+M\rangle^{\sigma/2}}
    \end{align}
    hold true for some constants $C_3=C_3(c)>0$ and $C_4=C_4(\epsilon)>0$.

\end{assumption}
Let's return to the statement $P_\mu P_{sc}=0$. Our approach to prove this is as follows: when $\sum\limits_{j=1}^3 P_j=v_T$ is fixed, we aim to identify a finite set of points such that if the energy of the $\s^2$ data $\psi$ is away from these finitely many values, then $P_\mu P_{sc}\psi=0$. We refer to these finitely many points as threshold points of $H$. We will now discuss the assumptions required to ensure that, for a given total momentum $v_T\in \R^3$, the set of all threshold points for $H$(defined in subsection \ref{sec: threshold}) is finite. Given a finite set $A\subseteq \R$, we let $|A|$ denote the cardinality of elements in $A$. 
\begin{assumption}\label{asp: threshold}We assume that for all $v_T\in \R^3$, the set 
    \eq
    \tau_{free}(v_T):=\{ (q_1,q_2,q_3)\in \R^9: v_1(q_1)=v_2(q_2)=v_3(q_3), \sum\limits_{j=1}^3 q_j=v_T \}
    \eeq
    is finite with 
    \eq
    \sup\limits_{v_T\in \R^3} |\tau_{free}(v_T)|<\infty. 
    \eeq

\end{assumption}
According to Assumptions \ref{asp: critical}, \ref{asp: subH}, \ref{asp: Hthreshold} and \ref{asp: threshold}, we prove that when $\sum\limits_{j=1}^3P_j=v_T\in \R^3$ is fixed, the set of the threshold points for $H$(defined in subsection \ref{sec: threshold}) is finite. Let $\tau(v_T)$ denote the set of all threshold points for $H$ when $\sum\limits_{j=1}^3P_j=v_T\in \R^3$.
\begin{lemma}\label{Lem: tauf}If Assumptions \ref{asp: critical}, \ref{asp: subH}, \ref{asp: Hthreshold} and \ref{asp: threshold} hold, then there exists a positive integer $N$ such that $|\tau(v_T)|\leq N$ for all $v_T\in \R^3$ and all $a=(jk)(l)\in L$. 
    
\end{lemma}
We defer the proof of Lemma \ref{Lem: tauf} to subsection \ref{sec: threshold}.

We're now ready to present the final assumption. Given the explanations that follow, we consider this assumption to be reasonable. Let $C_1, C_2\in \R$ two distinct constants. If $H_a$ is very close to $C_1$ and $H_0$ is very close to $C_2$, then roughly speaking, $V_{jk}(x_j-x_k)=H_a-H_0$ should be strictly away from $0$. This implies that it occurs in a region where $|x_j-x_k|$ is bounded. Similarly, if $H_a$ is very close to $C_1$ and $H$ is very close to $C_2$, then roughly speaking, $V_{jl}(x_j-x_l)+V_{kl}(x_k-x_l)=H-H_a$ should be strictly away from $0$. This indicates that this situation happens in the region where either $|x_j-x_l|$ or $|x_k-x_l|$ is bounded. Now, let's delve into the details of this assumption. For each $a=(jk)(l)\in L$, we define $\s^2_a(\R^9)$ as a Banach space:
\eq
\{f(x) : f=f_1+f_2, \text{ for some }f_1, f_2 \text{ with } \| \langle x_j-x_l\rangle^3 f_1\|_{\s^2_x(\R^9)}, \| \langle x_k-x_l\rangle^3 f_2\|_{\s^2_x(\R^9)} \},
\eeq
endowed with the norm:
\eq
\| f\|_{\s^2_a(\R^9)}:=\inf\limits_{(f_1,f_2)\in Ind_f} \| \langle x_j-x_l\rangle^3 f_1\|_{\s^2_x(\R^9)}+ \| \langle x_k-x_l\rangle^3 f_2\|_{\s^2_x(\R^9)},
\eeq
where $Ind_f$ denotes the set of all $(f_1,f_2)$ such that $f=f_1+f_2$ with $\langle x_j-x_l\rangle^3 f_1, \langle x_k-x_l\rangle^3 f_2\in \s^2_x(\R^9).$ Let $\tau^a(\eta)$ be the set of all threshold points of $H^a(\eta)$, as defined in Definition \ref{TPHa} in Section \ref{sec: threshold}. Given Assumption \ref{asp: Hthreshold}, we also find that $\sup\limits_{v_T\in \R^3}|\tau^a(v_T)|<\infty.$ Recall that in Lemma \ref{Lem: tauf}, we also have $\sup\limits_{v_T\in \R^3} |\tau(v_T)|<\infty$. Given that these threshold points are finite within each fiber, we introduce the subsequent assumption.
\begin{assumption}\label{asp: subHH} We assume that for all $a=(jk)(l)\in L$ and $c>0$, we have:
\begin{itemize}
\item[a)]
\eq
\| \langle x_j-x_k\rangle^{3}\bar{F}_{\tau}(H_a,c)F_{\tau}(H_0,\frac{c}{10})\|_{\s^2_x(\R^9)\to \s^2_x(\R^9)}\lesssim_{c}1,\label{Lem2.9: a}
\eeq
and
\eq
\| \langle x_j-x_k\rangle^3 F_{\tau^a}(H^a, c)\bar{F}_{\tau^a}(H_0^a, 10 c)\|_{\s^2_x(\R^6)\to \s^2_x(\R^6)}\lesssim_{ c}1,
\eeq
where $H^a_0:=H_0-\omega_l(P_l)$.
\item[b)] 
\eq
\| \bar{F}_{\tau}(H_a,c)F_{\tau}(H,\frac{c}{10})\|_{\s^2_x(\R^9)\to \s^2_a(\R^9)}\lesssim_{c}1,
\eeq
and 
\eq
\| F_{\tau}(H_a,c)\bar{F}_{\tau}(H,10c)\|_{\s^2_x(\R^9)\to \s^2_a(\R^9)}\lesssim_{c}1.
\eeq
\end{itemize}

\end{assumption}
\begin{proposition}\label{Prop: Pmu}If Assumptions \ref{asp: short} to \ref{asp: subHH} hold, then $P_\mu \Ps=0$ and $(1-P_\mu) \Ps=\Ps$.
\end{proposition}
Furthermore, we use Proposition \ref{Prop: Pmu} as a basis for introducing a new method of defining the projection onto the space of all bound states:
\begin{definition}[A new characterization of the projection on the space of all $3$-body bound states]\label{def:bd}Let $\alpha$ be as in Assumption \ref{asp: omega}. The projection on the space of all $3$-body bound states is defined by
\eq
P_\mu:= s\text{-}\lim\limits_{t\to \infty} e^{itH}[\Pi_{j=1}^{j=3} e^{-it\omega_j(P_j)}(1-F_{c,j,\alpha}(x,t,P))e^{it\omega_j(P_j)}]e^{-itH}\quad \text{ on }\s^2_x(\mathbb{R}^9).
\eeq
\end{definition}
To continue our analysis, we require the concept of following projections. Given $a=(jk)(l), b=(j'k')(l')\in L$, we define the following:
\eq
\Pp^{l'}_l:=s\text{-}\lim\limits_{t\to \infty}e^{itH_b} e^{-itH_a}F_{c,l,\alpha}(x_l,t,P_l)e^{itH_a}e^{-itH_b},\quad \text{ on }\s^2_x(\R^9),
\eeq 
\eq
P_{sc}(H_a):=s\text{-}\lim\limits_{t\to \infty} e^{itH_a}e^{-itH_0}F_{c,j,\alpha}(x_j,t,P_j)e^{itH_0} e^{-itH_a},\quad \text{ on }\s^2_x(\R^9)
\eeq
and 
\eq
P_{bs}(H_a):=1-P_{sc}(H_a).
\eeq
\begin{proposition}\label{Prop: Pab}If Assumptions \ref{asp: short}, \ref{asp: critical} and \ref{asp: omega} hold, then $\Pp_l^{l'}, P_{sc}(H_a), P_{bs}(H_a)$ and $\Omega_{\alpha}^{a,*}$ exist for all $a=(jk)(l),b=(j'k')(l')\in L-\{(1)(2)(3)\}$.
\end{proposition}
The proof of Proposition \ref{Prop: Pab} follows, similarly, by using Cook's method and propagation estimates, as detailed in Section \ref{Channel}.

Based on Theorem \ref{thm: channel}, Proposition \ref{Prop: Pmu} and Proposition \ref{Prop: Pab}, we arrive at our main result in this note:
\begin{theorem}\label{Thm}If Assumptions \ref{asp: short} to \ref{asp: subHH} hold, then we have AC for system \eqref{SE}: for all $\psi(0)=\Ps \psi(0)\in \s^2_x(\R^9)$, we have
\eq
\lim\limits_{t\to \infty}\|\psi(t)-\sum\limits_{a\in L}e^{-itH_a}\psi_{a,+}(x)\|_{\s^2_x(\R^9)}=0,\label{thm: AC}
\eeq
where for all $a\in L$, $\psi_{a,+}(x)$ are given by:
\begin{subequations}
For $a_0=(1)(2)(3)$, $a_l=(jk)(l)\in L$, 
\eq
\psi_{a_0,+}=\Omega^{a_1,*}\Omega_{a_1}^{*}\psi_0+\Omega^{a_2,*}(1-\Pp_{1}^{2}) \Omega_{a_2}^{*}\psi_0+\Omega^{a_3,*}(1-\Pp_{1}^{3})(1-\Pp_{2}^{3}) \Omega_{a_3}^{*}\psi_0,
\eeq
\eq
\psi_{a_1,+}:=P_{bs}(H_{a_1})\Omega_{a_1}^{*}\psi_0,
\eeq
\eq
\psi_{a_2,+}:=P_{bs}(H_{a_2})( 1-\Pp_{1}^{2})\Omega_{a_2}^{*}\psi_0
\eeq
and 
\eq
\psi_{a_3,+}:=P_{bs}(H_{a_3})(1-\Pp_1^{3})(1-\Pp_{2}^{3})\Omega_{a_3}^{*}\psi_0.
\eeq
\end{subequations}
\end{theorem}
\begin{proof}[Proof of Theorem \ref{Thm}]We define
\eq
a_1=(23)(1), \quad a_2=(13)(2)\text{ and }a_3=(12)(3).
\eeq
For $a=(jk)(l)\in L$, we set
\eq
J_{\alpha, a}(t):=e^{-itH_a}F_{c,l,\alpha}(x_l,t,P_l) e^{itH_a}.
\eeq
We choose $\psi(0)=\psi_0=\Ps \psi_0\in \s^2_x(\R^9)$. We can write the solution $\psi(t) = e^{-itH}\psi_0$ as
 \begin{align}
 e^{-itH}\psi_0=&J_{\alpha,a_1}(t)e^{-itH}\psi_0+(1-J_{\alpha,a_1}(t))J_{\alpha,a_2}(t)e^{-itH}\psi_0\nonumber\\
 &+(1-J_{\alpha,a_1}(t))(1-J_{\alpha,a_2}(t))J_{\alpha,a_3}(t)e^{-itH}\psi_0\nonumber\\
 &+(1-J_{\alpha,a_1}(t))(1-J_{\alpha,a_2}(t))(1-J_{\alpha,a_3}(t))e^{-itH}\psi_0.
 \end{align}
Using Theorem \ref{thm: channel} and Proposition \ref{Prop: Pmu}, we have that as $t\to \infty$:
\eq
\|(1-J_{\alpha,a_1}(t))(1-J_{\alpha,a_2}(t))(1-J_{\alpha,a_3}(t))e^{-itH}\psi_0\|_{\s^2_x(\R^9)}\to0
\eeq
and
\eq
\|e^{-itH}\psi_0-e^{-itH_{a_1}}\Omega_{a_1}^{*}\psi_0-\psi_{\alpha,2}(t)\|_{\s^2_x(\mathbb{R}^9)}\to 0,
\eeq
where
\eq
\psi_{\alpha,2}(t):=(1-J_{\alpha,a_1}(t))e^{-itH_{a_2}}\Omega_{a_2}^{*}\psi_0+(1-J_{\alpha,a_1}(t))(1-J_{\alpha,a_2}(t))e^{-itH_{a_3}}\Omega_{a_3}^{*}\psi_0.
\eeq
Using Proposition \ref{Prop: Pab}, we have that as $t\to \infty$:
\eq
\| (1-J_{\alpha,a_1}(t))e^{-itH_{a_2}} \Omega_{a_2}^{*}\psi_0 -e^{-itH_{a_2}}( 1-\Pp_{1}^{2})\Omega_{a_2}^{*}\psi_0 \|_{\s^2_x(\mathbb{R}^9)}\to 0,
\eeq
and
\eq
\|(1-J_{\alpha,a_1}(t))(1-J_{\alpha,a_2}(t))e^{-itH_{a_3}} \Omega_{a_3}^{*}\psi_0 -e^{-itH_{a_3}}(1-\Pp_1^{3})(1-\Pp_{2}^{3})\Omega_{a_3}^{*}\psi_0 \|_{\s^2_x(\mathbb{R}^9)}\to 0.
\eeq
Proposition \ref{Prop: Pab} implies that \eqref{thm: AC} holds by choosing
\eq
\psi_{a_0,+}=\Omega^{a_1,*}\Omega_{a_1}^{*}\psi_0+\Omega^{a_2,*}(1-\Pp_{1}^{2}) \Omega_{a_2}^{*}\psi_0+\Omega^{a_3,*}(1-\Pp_{1}^{3})(1-\Pp_{2}^{3}) \Omega_{a_3}^{*}\psi_0,
\eeq
\eq
\psi_{a_1,+}:=P_{bs}(H_{a_1})\Omega_{a_1}^{*}\psi_0,
\eeq
\eq
\psi_{a_2,+}:=P_{bs}(H_{a_2})( 1-\Pp_{1}^{2})\Omega_{a_2}^{*}\psi_0
\eeq
and 
\eq
\psi_{a_3,+}:=P_{bs}(H_{a_3})(1-\Pp_1^{3})(1-\Pp_{2}^{3})\Omega_{a_3}^{*}\psi_0.
\eeq
We finish the proof.
\end{proof}

\section{Basic dynamics}\label{sec:setting}
Let $\s^2(X_{v_T})$ be the space of all $\s^2$ functions with $\sum\limits_{j=1}^3P_j=v_T$, where $v_T \in \R^3$ is fixed. $\s^2(X_{v_T})$ is a fiber Hilbert space, which we will delve into later in this section. This section focuses on configuration space and three key properties of Schr\"odinger operators $H$ on the fiber Hilbert space $\s^2(X_{v_T})$: self-adjointness, Ruelle's theorem, and the set of all threshold points for $H$. By using Ruelle's theorem in $\s^2(X_{v_T})$, we get a Ruelle's theorem in $\s^2_x(\R^9)$. We can not get the Ruelle's theorem in $\s^2_x(\R^9)$ directly because, in the region where $|x_1-x_2|$, $|x_1-x_3|$, and $|x_2-x_3|$ are bounded, there exist points $x\in \R^9$ for which $|x|$ goes to infinity. Throughout this section, we assume the validity of Assumptions \ref{asp: short}- \ref{asp: threshold}. The structure of this section follows \cite{HS1}.
\subsection{Configuration space}\label{sec: space}
Given that $[P_1+P_2+P3, H]=0$, the total quasi-momentum of the three particles is conserved under the flow $e^{-itH}$ of the system \eqref{SE}. By decomposing $\s^2_x(\R^9)$ based on the total quasi-momentum, we obtain the fiber space, $\s^2(X_{v_T})$, which is the $\s^2$ space associated with a total quasi-momentum $v_T$. 

\subsubsection{The fiber Hilbert space}\label{sec: fiber space}
The configuration space $X$ of a quasi-three-particle system is a Euclidean space with scalar product denoted by $x\cdot y$. Specifically, we have
\eq
X\equiv \mathbb{R}^9\quad \text{ and }\quad x\cdot y\equiv \sum\limits_{j=1}^3 (x_j\cdot y_j)_{\mathbb{R}^3}.
\eeq
The quasi-momentum conjugate to $ X,$ is denoted by $p$. In quantum mechanics, this quasi-momentum has the conventional form in Cartesian coordinates (not particle coordinates) of $X$, which is given by $p=-i\nabla_x$. For a given $v_T\in \mathbb{R}^3$, we introduce $P_{v_T}$ as a framework representing the total quasi-momentum $v_T$. This can be symbolically expressed as:
\eq
P_{v_T}\equiv \{ p=(p_1,p_2,p_3) : p_j\in \mathbb{R}^3\text{ and } \sum\limits_{j=1}^3 p_j=v_T \}.
\eeq
The fiber Hilbert space $\s^2(P_{v_T})$ is defined as the space of all $\s^2$ functions on $P_{v_T}$:
\eq
\s^2(P_{v_T})\equiv \{ \hat{f}(p): \int_{\{p_1+p_2+p_3=v_T\}}d^9p |\hat{f}(p)|^2<\infty  \},
\eeq
which can also be expressed in the form:
\eq
\s^2(P_{v_T})=\{ \hat{f}(p): \hat{f}(p)= g_{v_T}(p_1,p_2) \quad \text{for some }g_{v_T}(p_1,p_2)\in \s^2_{p_1,p_2}(\mathbb{R}^6)\}.
\eeq
Subsequently, we characterize the fiber Hilbert space, $\s^2(X_{v_T})$, as the dual of $\s^2(P_{v_T})$:
\eq
\s^2(X_{v_T}):=\{f_{v_T}(x): f_{v_T}(x)=\tilde{f}_{v_T}(x_1-x_2,x_2-x_3)e^{iv_T\cdot x_3}, \text{ for some }\tilde{f}_{v_T}(y_1,y_2)\in \s^2_{y_1,y_2}(\mathbb{R}^6)\}.
\eeq
\begin{remark}
We can view $\s^2(X_{v_T})$ as $\s^2_{x_1,x_2}(\mathbb{R}^6)$ due to the homeomorphism $\mathcal{L}_{v_T}: \s^2(X_{v_T})\to \s^2_{x_1,x_2}(\mathbb{R}^6), f_{v_T}\mapsto \mathcal{L}_{v_T}(f_{v_T}):=\tilde{f}_{v_T}(x_1,x_2)$. Specifically, this relationship is given by:
\eq
\mathcal{L}_{v_T}(f_{v_T})=e^{-ix_3\cdot v_T} e^{ix_3\cdot P_2}e^{ix_2\cdot P_1}f_{v_T}.
\eeq 
\end{remark}
For all $1\leq p\leq \infty$, we define $\s^p(X_{v_T})$ as:
\eq
\s^p(X_{v_T}):=\{f_{v_T}(x): f_{v_T}(x)=\tilde{f}_{v_T}(x_1-x_2,x_2-x_3)e^{iv_T\cdot x_3}, \text{ for some }\tilde{f}_{v_T}(y_1,y_2)\in \s^p_{y_1,y_2}(\mathbb{R}^6)\}.
\eeq
In this paper, we employ the notations $W^{k,p}(X_{v_T})$, $\Hi^k(X_{v_T})$, $W^{k,p}(\mathbb{R}^9)$, and $\Hi^k(\mathbb{R}^9)$ to denote the respective Sobolev spaces. Throughout this section, we adopt the convention that for all $f_{v_T}\in \s^p(X_{v_T})$, $\Hi^k(X_{v_T})$, or $W^{k,p}(X_{v_T})$,

\eq
\tilde{f}_{v_T}(x_1-x_2,x_2-x_3):=f_{v_T}(x)e^{-iv_T\cdot x_3}.\label{convention}
\eeq
Recall that $P_j=-i\nabla_{x_j}$, $j=1,2,3$.
\begin{lemma}\label{def:lem1}For all  $v_T\in \mathbb{R}^3$ and $f_{v_T}\in \s^2(X_{v_T})$, 
\eq
(P_1+P_2+P_3)f_{v_T}(x)=v_T f_{v_T}(x).
\eeq
\end{lemma}
\begin{proof}This conclusion arises from the following computation:
\begin{align}
(P_1+P_2+P_3)f_{v_T}=&(-i)e^{ix_3\cdot v_T}(\nabla_1-\nabla_1+\nabla_2-\nabla_2)[\tilde{f}_{v_T}(x_1-x_2,x_2-x_3)]+v_T f_{v_T}\\
=&v_T f_{v_T}.
\end{align}
Here, $\nabla_1\tilde{f}_{v_T}(y_1,y_2):=\nabla_{y_1}\tilde{f}_{v_T}$ and $\nabla_2\tilde{f}_{v_T}(y_1,y_2):=\nabla_{y_2}\tilde{f}_{v_T}$. 
\end{proof}
The Hilbert space $\s^2(X_{v_T})$ is endowed with an inner product, which is defined as 
\eq
(f_{v_T}(x), g_{v_T}(x))_{\s^2(X_{v_T})}:=\int d^3y_1d^3y_2(\tilde{f}_{v_T}(y_1,y_2))^*\tilde{g}_{v_T}(y_1,y_2)\label{l2Xvt}
\eeq
for $f_{v_T}(x)=\tilde{f}_{v_T}(x_1-x_2,x_2-x_3)e^{iv_T\cdot x_3},$ $g_{v_T}(x)=\tilde{f}_{v_T}(x_1-x_2,x_2-x_3)e^{iv_T\cdot x_3}\in \s^2(X_{v_T})$. Therefore, we have
\eq
\begin{split}
\| f_{v_T}(x)\|_{\s^2(X_{v_T})}=&\| \tilde{f}_{v_T}(x_1,x_2) \|_{\s^2_{x_1,x_2}(\mathbb{R}^6)}\\
=&\|\tilde{f}_{v_T}(x_1-x_2,x_2)   \|_{\s^2_{x_1,x_2}(\mathbb{R}^6)}\\
=&\|   e^{ix_3\cdot P_2}e^{ix_3\cdot P_1}f_{v_T}(x)\|_{\s^2_{x_1,x_2}(\mathbb{R}^6)}
\end{split}\label{Feb.22.eq.1}
\eeq
Due to the following lemma and based on \eqref{l2Xvt}, $\s^2(X_{v_T})$ is a proper Hilbert space for system \eqref{SE} since Lemma \ref{Lem: H close} implies that $e^{-itH}$ is closed on $\s^2(X_{v_T})$:
\begin{lemma}\label{Lem: H close}$[V(x),P_1+P_2+P_3]=0$.
\end{lemma}
\begin{proof}
It follows from that for all $1\leq j<l\leq 3$, we have
\eq
[(x_j-x_l), P_1+P_2+P_3]=0.
\eeq
\end{proof}
\subsubsection{Channels and Hamiltonians}\label{subsection: Channels}
 The definition of channels is similar to what was done in \cite{HS1}. In the configuration space of a quantum $N$-body system, $X$, there is a distinguished, finite lattice $L$ of subspaces $a,b,\cdots$(channels). $L$ is closed under intersections and contains at least $a=\{0\}$ and $a=X$. In the case of $H=H_0+\sum\limits_{1\leq j<l\leq N}V_{jl}(x_j-x_l)$ the channels correspond to all partitions of $(1,\cdots,N)$ into clusters. For example, if $N=4$:
\eq
\begin{split}
\text{partition} &~~~~\text{ channel}\\
(12)(34) &  \leftrightarrow a= \{ x\vert x_1=x_2; x_3=x_4\},
\end{split}
\eeq
and if $N=3$:
\eq
\begin{split}
\text{partition} &~~~~\text{ channel}\\
(12)(3) &  \leftrightarrow a= \{ x\vert x_1=x_2\}.
\end{split}
\eeq
For simplicity, we write 
\eq
a=\text{the corresponding partition of }(1,\cdots,N).
\eeq
For example, when $N=3$, we write $a=(12)(3)$ if $a= \{ x\vert x_1=x_2\}$, $a=(1)(2)(3)$ if $a=X$ and $a=0$ if $a=\emptyset$. In general the partial ordering of $L$ is defined by 
\eq
a<b \leftrightarrow a \subset b; \quad a\neq b.
\eeq
For $a\in L,$ let $\I_a:=\{(j,l) : 1\leq j<l\leq N, \{ x: x_j=x_l\}\supset a\}$,
\eq
\tilde{\I}_a:=\{j: (j,l)\in \I_a \text{ for some }l \text{ or } j=\max\limits_{(l,k)\in \I_a} k \},
\eeq
\eq
 V_a(x):=\sum\limits_{(j,l)\in \I_a} V_{jl}(x_j-x_l),\label{3.8.1}
\eeq
\eq
H_a:=H_0+V_a(x),\label{3.8.2}
\eeq
\eq
H_0^a:=\sum\limits_{j\in \tilde{\I}_a} \omega_j(P_j)\label{3.8.3}
\eeq
and 
\eq
H^a:=H_0^a+V_a(x).\label{3.8.4}
\eeq
In the case of three quasi-particles, when $a = (jk)(l) \in L$, the definitions of $H^a$ and $H_a$ align with what we discussed in the introduction: 
\eq
\I_a=\{ (j,k)\},
\eeq
\eq
\tilde{\I}_a=\{j,k\},
\eeq
\eq
V_a(x)=V_{jk}(x_j-x_k),
\eeq
\eq
H_a=H_0+V_{jk}(x_j-x_k),
\eeq
\eq
H^a_0=\omega_j(P_j)+\omega_k(P_k)
\eeq
and 
\eq
H^a=H_0^a+V_{jk}(x_j-x_k).
\eeq
Here is an example of the case when $N=3$ and $a=(12)(3)$:

\begin{example}If $N=3$ and $a=(12)(3)$, then $V_a(x)=V_{12}(x_1-x_2)$, $H_a=H_0+V_{12}(x_1-x_2)$ and $H^a=\omega_1(P_1)+\omega_2(P_2)+V_{12}(x_1-x_2)$.
\end{example}
\begin{proof}The expressions provided can be deduced from equations \eqref{3.8.1} to \eqref{3.8.4} when choosing $a=(12)(3)$. In such scenarios, as $|x_1-x_3|\to \infty$ and $|x_2-x_3|\to \infty$, the influences of $V_{13}$ and $V_{23}$ is weaken, allowing us to consider $H$ as $H_a$. Geometrically, the evolution of $e^{-itH}\psi$ towards $e^{-itH_a}\psi_a$ for some $\psi_a\in \s^2_x(\R^9)$ as $t\to \infty$ can be visualized as particle 1, bound with particle 2, distancing itself from particle 3. Refer to Figure \ref{p3} for an illustration.
\vspace{10mm}
\begin{figure}[h]
\centering
\includegraphics[width=8cm]{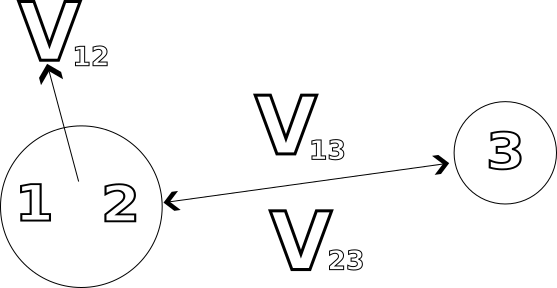}

\caption{Visualization of the $e^{-itH_a}$ flow for $a=(12)(3)$}
\label{p3}
\end{figure}
    
\end{proof}
\subsection{Ruelle's Theorem in the Context of Fiber Space}\label{2.2}
In this section, we present a proof of Ruelle's theorem in fiber space $\s^2(X_{v_T})$ is presented, utilizing self-adjointness and local compactness, which are the two essential components of the proof of Ruelle's theorem (see e.g., Theorem 2.4 in \cite{HS1}). Ruelle's theorem provides a characterization of all scattering states in $\s^2(X_{v_T})$. \par
We initiate our discussion with the self-adjointness of $H$. Define 
\eq
\mathcal{D}(H):=\{ f_{v_T}(x)\in \s^2(X_{v_T}): Hf_{v_T} \in \s^2(X_{v_T})\}.
\eeq
For any given pair $f_{v_T}, g_{v_T}\in \mathcal{D}(H)$, the following identity holds true:
\eq
(f_{v_T}(x), Hg_{v_T}(x) )_{\s^2(X_{v_T})}=(\tilde{f}_{v_T}(x_1-x_2,x_2), \tilde{H}(v_T,3)\tilde{g}_{v_T}(x_1-x_2,x_2))_{\s^2_{x_1,x_2}(\mathbb{R}^6)}.\label{Jan.1}
\eeq
Here,
\eq
\tilde{H}_0(v_T,3):=\omega_1(P_1)+\omega_2(P_2)+\omega_{3}(v_T-P_1-P_2),\label{tH0}
\eeq
and
\eq
\tilde{V}_3(x_1,x_2):=V_{12}(x_1-x_2)+V_{13}(x_1)+V_{23}(x_2),~\text{ and }\tilde{H}(v_T,3):=\tilde{H}_0(v_T)+\tilde{V}_3(x_1,x_2).\label{tH}
\eeq
\begin{lemma}\label{Lem:jan1}\eqref{Jan.1} is true.
\end{lemma}
\begin{proof}First of all, we have the following identity:
\begin{align}
H=&e^{-ix_3\cdot P_1}e^{-ix_3\cdot P_2}(\omega_1(P_1)+\omega_2(P_2)+\omega_3(P_3-P_2-P_1) +\tilde{V}(x_1,x_2))e^{ix_3\cdot P_1}e^{ix_3\cdot P_2}\nonumber\\
=&e^{-ix_3\cdot P_1}e^{-ix_3\cdot P_2} \tilde{H}(P_3,3)e^{ix_3\cdot P_1}e^{ix_3\cdot P_2}.\label{Lem2.3: eq1}
\end{align}
Employing the convention from \eqref{convention}, it's evident that for all $f_{v_T}\in \s^2(X_{v_T})$,
\eq
\tilde{f}_{v_T}(x_1-x_2,x_2)=e^{-iv_T\cdot x_3}e^{ix_3\cdot P_1}e^{ix_3\cdot P_2}f_{v_T}(x).\label{Lem2.3: eq2}
\eeq
Furthermore, based on \eqref{l2Xvt}, we can deduce:
\begin{align}
(f_{v_T}(x), g_{v_T}(x))_{\s^2(X_{v_T})}=& \int d^3y_1d^3y_2(\tilde{f}_{v_T}(y_1,y_2))^*\tilde{g}_{v_T}(y_1,y_2)  \nonumber\\
(\text{take }y_1=x_1-x_2,y_2=x_2)=& \int d^3x_1d^3x_2(\tilde{f}_{v_T}(x_1-x_2,x_2))^*\tilde{g}_{v_T}(x_1-x_2,x_2).\label{Lem2.3: eq3}
\end{align}
According to \eqref{Lem2.3: eq1}, \eqref{Lem2.3: eq2} and \eqref{Lem2.3: eq3}, their inner product can be expressed as follows:
\begin{align}
&(f_{v_T}(x), Hg_{v_T}(x) )_{\s^2(X_{v_T})}\\
=&(\tilde{f}_{v_T}(x_1-x_2,x_2), \widetilde{H g}_{v_T}(x_1-x_2,x_2))_{\s^2_{x_1,x_2}(\mathbb{R}^6)}\\
=&(\tilde{f}_{v_T}(x_1-x_2,x_2)e^{iv_T\cdot x_3},e^{-iv_T\cdot x_3}e^{ix_3\cdot P_1}e^{ix_3\cdot P_2}Hg_{v_T}(x))_{\s^2_{x_1,x_2}(\mathbb{R}^6)}\\
=&(\tilde{f}_{v_T}(x_1-x_2,x_2)e^{iv_T\cdot x_3},\tilde{H}(v_T,3)e^{-iv_T\cdot x_3}e^{ix_3\cdot P_1}e^{ix_3\cdot P_2}g_{v_T}(x))_{\s^2_{x_1,x_2}(\mathbb{R}^6)}\\
=&(\tilde{f}_{v_T}(x_1-x_2,x_2), \tilde{H}(v_T,3) \tilde{g}_{v_T}(x_1-x_2,x_2))_{\s^2_{x_1,x_2}(\mathbb{R}^6)}.
\end{align}
\end{proof}
Lemma \ref{Lem:jan1} implies that the self-adjointness of $H$ on $\mathcal{D}(H)$ corresponds to the self-adjointness of $\tilde{H}(v_T, 3)$ on $\mathcal{D}(\tilde{H}(v_T, 3)):=\{f(x_1,x_2)\in \s^2(\mathbb{R}^6): \tilde{H}(v_T, 3)f \in \s^2(\mathbb{R}^6) \}$, which has been well established in literature (see e.g. page 20 in \cite{CL2002} or \cite{RS4}). Consequently, $H$ is self-adjoint on $\mathcal{D}(H)$. Recall that $\tx=(x_1-x_2,x_2-x_3,x_3-x_1)\in \mathbb{R}^9$.
{\begin{lemma}[Local compactness property] Let $f(x):=\tilde{f}(x_1-x_2,x_2-x_3)\in \s^\infty_x(\mathbb{R}^9)$ with 
\eq
f(x)=\tilde{f}(x_1-x_2,x_2-x_3)\to 0 \label{asp: f}
\eeq
as $|\tilde{x}|\to \infty$. Under Assumption \ref{asp: critical b} part (a), the operator 
\eq
f(x)(z-H)^{-1}\text{ is compact on }\s^2(X_{v_T}) 
\eeq
for any $z$ in the resolvent set $\rho(H)$. We refer to this by saying that $H$ has the local compactness property. 
\end{lemma}
\begin{proof}Here, $|\tilde{x}|$ is equivalent to $|x_1-x_2|+|x_1-x_3|$ because 
\eq
|x_1-x_2|+|x_2-x_3|\lesssim \sqrt{|x_1-x_2|^2+|x_1-x_3|^2+|x_2-x_3|^2}
\eeq
and 
\begin{align}
\sqrt{|x_1-x_2|^2+|x_1-x_3|^2+|x_2-x_3|^2}\lesssim& |x_1-x_2|+|x_1-x_3|+|x_2-x_3|\nonumber\\
\lesssim &|x_1-x_2|+|x_2-x_3|
\end{align}
by using triangle inequality. Let $g_{v_T}\in \s^2(X_{v_T})$. Based on the definition of $\s^2(X_{v_T})$, we have $f(x)(H-z)^{-1}g_{v_T}\in \s^2(X_{v_T})$. We rewrite $f(x)(H-z)^{-1}g_{v_T}$ as
\eq
f(x)(H-z)^{-1}g_{v_T}=e^{-ix_3\cdot P_1}e^{-ix_3\cdot P_2} \tilde{f}(x_1,x_2) (\tilde{H}(v_T,3)-z)^{-1}\tilde{g}_{v_T}(x_1,x_2)e^{iv_T\cdot x_3}.\label{proof: eq1}
\eeq
Using the second resolvent identity, we obtain 
\begin{align}
\tilde{f}(x_1,x_2)(\tilde{H}(v_T,3)-z)^{-1}=&\tilde{f}(x_1,x_2)(\tilde{H}_0(v_T,3)-z)^{-1}\nonumber\\
&-\tilde{f}(x_1,x_2)(\tilde{H}_0(v_T,3)-z)^{-1}\tilde{V}_3(x_1,x_2)(\tilde{H}(v_T,3)-z)^{-1}.\label{Feb.21.1}
\end{align}
The compactness of $\tilde{f}(x_1,x_2)(\tilde{H}_0(v_T,3)-z)^{-1}$ on $\s^2_{x_1,x_2}(\mathbb{R}^6)$ is due to Assumption \ref{asp: critical b} part (a) and \eqref{asp: f}. By utilizing \eqref{Feb.21.1} and the compactness of $\tilde{f}(x_1,x_2)(\tilde{H}_0(v_T,3)-z)^{-1}$ on $\s^2_{x_1,x_2}(\mathbb{R}^6)$, we can conclude that $f(x)(\tilde{H}(v_T,3)-z)^{-1}$ is also compact on $\s^2(X_{v_T})$. Using \eqref{proof: eq1} and the compactness of $f(x)(\tilde{H}(v_T,3)-z)^{-1}$ on $\s^2(X_{v_T})$, we can further conclude that $f(x)(H-z)^{-1}$ is also compact on $\s^2(X_{v_T})$.
\end{proof}}
\begin{lemma}[Ruelle's theorem in fiber space]\label{Rfiber}Suppose that $H=H^*$ on $\mathcal{D}(H)$ has the {\emph{local compactness property}}. Let $\mathcal{H}_B=\mathcal{H}_B(v_T)$ be the subspace spanned by all eigenvectors of $H$, and $\mathcal{H}_C=\mathcal{H}_C(v_T)=\Hi_B^\perp$. If $\chi_R(\tx)$ is the characteristic function of some ball $|\tx|<R$, then 
\eq
\varphi\in \Hi_B\iff \lim\limits_{R\to \infty}\| (1-\chi_R)e^{\pm itH}\varphi\|_{\s^2(X_{v_T})}=0\text{ uniformly in }0\leq t<\infty;
\eeq
\eq
\psi \in \Hi_C \iff \lim\limits_{t\to \infty} t^{-1}\int_0^tds \|\chi_Re^{\pm isH}\psi\|_{\s^2(X_{v_T})}^2=0\quad \text{ for any }R<\infty.
\eeq
In particular, there exists a sequence of time $\{t_n\}_{n=1}^{n=\infty}$, with $t_n=t_n(\psi, R)\uparrow  \infty$ as $n\to \infty$, such that 
\eq
\|\chi_R e^{-it_nH} P_c\psi\|_{\s^2(X_{v_T})}<1/n,\quad \text{ for all }n \in \mathbb{N},
\eeq
where $P_c=P_c(v_T)$ denotes the projection on $\mathcal{H}_C$.
\end{lemma}
\begin{proof}See the proof of Theorem 2.4 in \cite{HS1}. 
\end{proof}
By applying Ruelle's theorem in $\s^2(X_{v_T})$, we obtain a precise definition of $P_{c}(v_T)$, which is the projection onto the space of all scattering states in $\s^2(X_{v_T})$.
\subsection{Ruelle's theorem in the full space and the scattering projection}\label{2.3}
In this section, we demonstrate how to prove Ruelle's theorem in $\s^2_x(\mathbb{R}^9)$ by utilizing Ruelle's theorem in fiber space $\s^2(X_{v_T})$ as described in Lemma \ref{Rfiber}. Let $\psi\in C^\infty_0(\mathbb{R}^9)$. We can write $\psi$ as
\eq
\psi(x_1,x_2,x_3)=e^{-ix_3\cdot P_1}e^{-ix_3\cdot P_2}\psi(x_1+x_3,x_2+x_3,x_3).
\eeq
Define $\tilde{\psi}(x_1,x_2,x_3):= \psi(x_1+x_3,x_2+x_3,x_3)$. Then $\tilde{\psi}\in \s^2_x(\mathbb{R}^9)$ and for any $t\in \mathbb{R}$, we have 
\eq
e^{-itH}\psi=\frac{1}{(2\pi)^{3/2}}e^{-ix_3\cdot P_1}e^{-ix_3\cdot P_2} \int d^3q_3 e^{-it\tilde{H}(q_3,3)} \tilde{\psi}(x_1,x_2,\hat{q}_3),
\eeq
where $\tilde{\psi}(x_1,x_2,\hat{q}_3)$ denotes the Fourier transform in $x_3$ variable 
\eq
\tilde{\psi}(x_1,x_2,\hat{q}_3):=\mathscr{F}_{x_3}[\tilde{\psi}](q_3)=\frac{1}{(2\pi)^{3/2}}\int d^3x_3 e^{-iq_3\cdot x_3} \tilde{\psi}(x_1,x_2,x_3)
\eeq
and by using Fourier inversion theorem, we have
\eq
\tilde{\psi}(x)=\mathscr{F}_{q_3}^{-1}[\tilde{\psi}](x_3)=\frac{1}{(2\pi)^{3/2}}\int d^3q_3 e^{iq_3\cdot x_3} \tilde{\psi}(x_1,x_2,\hat{q}_3).
\eeq
Here, for $\tilde{H}(q_3,3)$, see \eqref{tH0} and \eqref{tH}. $\psi$ has following orthogonal decomposition 
\eq
\psi=\psi_{sc}+\psi_{bs}
\eeq
where 
\begin{align}
\psi_{sc}(x):=&P_c(P_1+P_2+P_3)e^{-ix_3\cdot P_1}e^{-ix_3\cdot P_2}\tilde{\psi}(x)\nonumber\\
=&e^{-ix_3\cdot P_1}e^{-ix_3\cdot P_2} P_c(P_3)\tilde{\psi}(x)
\end{align}
and 
\begin{align}
\psi_{bs}(x):=&(1-P_c(P_1+P_2+P_3))e^{-ix_3\cdot P_1}e^{-ix_3\cdot P_2}\tilde{\psi}(x)\nonumber\\
=&e^{-ix_3\cdot P_1}e^{-ix_3\cdot P_2} (1-P_c(P_3))\tilde{\psi}(x).
\end{align}
We use $P_c(P_1+P_2+P_3)$ here because $P_1+P_2+P_3$ is really the total quasi-momentum of three particles. Here, we also use
\eq
P_c(P_1+P_2+P_3)e^{-ix_3\cdot P_1}e^{-ix_3\cdot P_2}=e^{-ix_3\cdot P_1}e^{-ix_3\cdot P_2}P_c(P_3),
\eeq
where $P_c(P_3)$ on $\s^2_x(\mathbb{R}^9)$ is understood as 
\eq
P_c(P_3)f=\frac{1}{(2\pi)^{3/2}}\int d^3q_3 e^{ix_3\cdot q_3} P_c(q_3) f(x_1,x_2,\hat{q}_3),\quad f\in \s^2_x(\mathbb{R}^9),
\eeq
and $P_c(q_3)$ is defined in Lemma \ref{Rfiber} for all $q_3\in \mathbb{R}^3$ by taking $H=\tilde{H}(q_3,3)$. Given the equivalence between $\tilde{H}(P_3, 3)$ and our three-body Hamiltonian $H$ on $\s^2_x(\mathbb{R}^9)$, we define $P_{sc}$ as the projection onto the space of all scattering states of $H$:
\eq
P_{sc}:=P_c(P_1+P_2+P_3).
\eeq
We also define $P_{bs}$ as the projection onto the spac eof all bound states of $H$:
\eq
P_{bs}:=1-P_{sc}.
\eeq
Consequently, we can represent the scattering component of $\psi$ as 
\eq
\psi_{sc}=P_{sc}\psi, \label{psic}
\eeq
and the non-scattering component of $\psi$ as 
\eq
\psi_{bs}=P_{bs}\psi.\label{psib}
\eeq
Note that $\psi_{sc}(x)$ is orthogonal to $\psi_{bs}(x)$ in $\s^2_x(\mathbb{R}^9)$ due to Plancherel's theorem: for some constant $c>0$, we have
\begin{align}
(\psi_{sc}(x), \psi_{bs}(x))_{\s^2_x(\mathbb{R}^9)}=&(\tilde{\psi}_{sc}(x), \tilde{\psi}_{bs}(x))_{\s^2_x(\mathbb{R}^9)}\nonumber\\
=&c \int d^3q_3 (P_c(q_3)\tilde{\psi}_{sc}(x_1,x_2,\hat{q}_3), (1-P_c(q_3))\tilde{\psi}_{bs}(x_1,x_2,\hat{q}_3) )_{\s^2_{x_1,x_2}(\mathbb{R}^6)}\nonumber\\
=&0.
\end{align}
We can apply the same orthogonal decomposition to all $\psi\in \s^2_x(\mathbb{R}^9)$ by using the density of $C_0^\infty(\mathbb{R}^9)$ in $\s^2_x(\mathbb{R}^9)$:
\begin{proposition}For all $\psi\in \s^2_x(\mathbb{R}^9)$, $\psi$ has an orthogonal decomposition
\eq
\psi(x)=\psi_{sc}(x)+\psi_{bs}(x)
\eeq
where $\psi_{sc}$ and $\psi_{bs}$ are defined in \eqref{psic} and \eqref{psib}, respectively.
\end{proposition}
In the following, we provide a proof of Ruelle's theorem in $\s^2_x(\mathbb{R}^9)$ space:
\begin{proposition}[Ruelle's theorem in $\s^2_x(\mathbb{R}^9)$]\label{Rwhole}If Assumptions \eqref{asp: short} and \ref{asp: critical b} part (a) hold, and if $\chi_R(\tx)$ is the characteristic function of some ball $|\tx|<R$, then for all $\psi\in \s^2_x(\mathbb{R}^9)$,
\eq
 \lim\limits_{R\to \infty}\| (1-\chi_R)e^{\pm itH}\psi_{bs}\|_{\s^2_x(\mathbb{R}^9)}=0\text{ uniformly in }0\leq t<\infty;\label{Feb.22.goal1}
\eeq
\eq
 \lim\limits_{t\to \infty} t^{-1}\int_0^tds \|\chi_Re^{\pm isH}\psi_{sc}\|_{\s^2_x(\mathbb{R}^9)}^2=0\quad \text{ for any }R<\infty.\label{Feb.22.goal2}
\eeq
In particular, there exists a sequence of time $\{t_n\}_{n=1}^{n=\infty}$, with $t_n=t_n(\psi, R)\uparrow  \infty$ as $n\to \infty$, such that 
\eq
\|\chi_R e^{-it_nH} \psi_{sc}\|_{\s^2(\mathbb{R}^9)}<1/n,\quad \text{ for all }n \in \mathbb{N}.
\eeq
\end{proposition}
\begin{proof}Choose $\psi\in \s^2_x(\mathbb{R}^9)$. For any $\epsilon>0$, there exists $\psi_{\epsilon}\in C^\infty_0(\mathbb{R}^9)$ such that $\| \psi-\psi_\epsilon\|_{\s^2_x(\mathbb{R}^9)}<\epsilon$. Our goal is to show that 
\eq
 \limsup\limits_{R\to \infty}\| (1-\chi_R)e^{\pm itH}[\psi_\epsilon]_{bs}\|_{\s^2_x(\mathbb{R}^9)}=0\quad\text{ uniformly in }0\leq t<\infty\label{Feb.22.1}
\eeq
and 
\eq
 \limsup\limits_{t\to \infty} t^{-1}\int_0^tds \|\chi_Re^{\pm isH}[\psi_\epsilon]_{sc}\|_{\s^2_x(\mathbb{R}^9)}^2=0\quad \text{ for any }R<\infty,\label{Feb.22.2}
\eeq
where $\psi_\epsilon=[\psi_\epsilon]_{sc}+[\psi_\epsilon]_{bs}$ with $ [\psi_\epsilon]_{sc}$ and $ [\psi_\epsilon]_{bs} $ defined in \eqref{psic} and \eqref{psib}, respectively, by substituting $\psi$ with $\psi_\epsilon$. Once we establish \eqref{Feb.22.1} and \eqref{Feb.22.2}, we can obtain the following inequalities for any $\epsilon >0$:
\eq
 \limsup\limits_{R\to \infty}\| (1-\chi_R)e^{\pm itH}\psi_{bs}\|_{\s^2_x(\mathbb{R}^9)}<\epsilon\text{ uniformly in }0\leq t<\infty\label{Feb.22.3}
\eeq
and 
\eq
 \limsup\limits_{t\to \infty} t^{-1}\int_0^tds \|\chi_Re^{\pm isH}\psi_{sc}\|_{\s^2_x(\mathbb{R}^9)}^2<\epsilon^2\quad \text{ for any }R<\infty,\label{Feb.22.4}
\eeq
since due to unitarity of $e^{-itH}$ on $\s^2_x(\mathbb{R}^9)$, we have
\eq
\begin{split}
 t^{-1}\int_0^tds \|\chi_Re^{\pm isH}(\psi_{sc}-[\psi_\epsilon]_{sc})\|_{\s^2_x(\mathbb{R}^9)}^2\leq &\|\psi_{sc}-[\psi_\epsilon]_{sc}\|_{\s^2_x(\mathbb{R}^9)}^2 \\
= &\| P_c(P_1+P_2+P_3)(\psi-\psi_\epsilon)\|_{\s^2_x(\mathbb{R}^9)}^2\\
\leq& \| \psi-\psi_\epsilon\|_{\s^2_x(\mathbb{R}^9)}^2
\end{split}
\eeq
and
\eq
\begin{split}
\| (1-\chi_R)e^{\pm itH}(\psi_{bs}-[\psi_\epsilon]_{bs})\|_{\s^2_x(\mathbb{R}^9)}\leq &\|\psi_{bs}-[\psi_\epsilon]_{bs}\|_{\s^2_x(\mathbb{R}^9)}\\
= &\| (1-P_c(P_1+P_2+P_3))(\psi-\psi_\epsilon)\|_{\s^2_x(\mathbb{R}^9)}\\
\leq& \| \psi-\psi_\epsilon\|_{\s^2_x(\mathbb{R}^9)}.
\end{split}
\eeq
Therefore, we conclude \eqref{Feb.22.goal1} and \eqref{Feb.22.goal2}. Now let us prove \eqref{Feb.22.1} and \eqref{Feb.22.2}. Let us define
\begin{align}
\tilde{\chi}_R(x_1,x_2):=&e^{ix_3\cdot P_2}e^{ix_3\cdot P_1} \chi_R(\tilde{x})e^{-ix_3\cdot P_1}e^{-ix_3\cdot P_2}\nonumber\\
=&e^{ix_3\cdot P_2}e^{ix_3\cdot P_1} \chi_R(x_1-x_2,x_2-x_3,x_3-x_1)e^{-ix_3\cdot P_1}e^{-ix_3\cdot P_2}\nonumber\\
=&\chi_R(x_1-x_2, x_2, -x_1),
\end{align}
\eq
\tilde{\psi}_\epsilon=e^{ix_3\cdot P_2}e^{ix_3\cdot P_1}\psi_\epsilon,\label{psisctilde}
\eeq
\eq
[ \tilde{\psi}_\epsilon]_{sc}:=P_c(P_3)\tilde{\psi}_\epsilon
\eeq
and 
\eq
[ \tilde{\psi}_\epsilon]_{bs}:=(1-P_c(P_3))\tilde{\psi}_\epsilon.
\eeq
Since $\psi_\epsilon\in C_0^\infty(\mathbb{R}^9)$, it follows that $\tilde{\psi}_\epsilon(x_1,x_2,\hat{q}_3) \in \s^2_{x_1,x_2}(\mathbb{R}^6)$ for all $q_3\in \mathbb{R}^3$. Hence, we have
\eq
e^{-ix_3\cdot P_1}e^{-ix_3\cdot P_2} \tilde{\psi}_\epsilon(x_1,x_2,\hat{q}_3)e^{ix_3\cdot q_3}=\tilde{\psi}_\epsilon(x_1-x_3, x_2-x_3)e^{ix_3\cdot q_3}\in \s^2(X_{q_3}),
 \eeq
since, by taking 
\eq
\bar{\psi}_\epsilon(y_1,y_2)=\tilde{\psi}_\epsilon(y_1+y_2,y_2),
\eeq
we have $\tilde{\psi}_\epsilon(x_1-x_3, x_2-x_3)=\bar{\psi}_\epsilon(x_1-x_2,x_2-x_3)\in \s^2_{y_1,y_2}(\mathbb{R}^6)$ satisfying $$\|\tilde{\psi}_\epsilon(y_1, y_2)\|_{\s^2_{y_1,y_2}(\mathbb{R}^6)}=\|\bar{\psi}_\epsilon(y_1, y_2)\|_{\s^2_{y_1,y_2}(\mathbb{R}^6)}.$$
Using Assumptions \ref{asp: short} and \ref{asp: critical} part (b), we conclude that for all $q_k\in \R^3$, $\tilde{H}(q_k,3)$ is self-adjoint and has local compactness. Given the self-adjointness and local compactness of $\tilde{H}(q_k,3)$, the conclusion of Lemma \ref{Rfiber} is valid for $\tilde{H}(q_k,3)$. Therefore, we have that for all $q_3\in \mathbb{R}^3$ by applying \eqref{psisctilde}, \eqref{convention} and \eqref{Feb.22.eq.1},
\eq
\begin{split}
&\lim\limits_{R\to \infty}\| (1-\chi_R)e^{\pm itH}[\psi_\epsilon]_{bs}\|_{\s^2(X_{q_3})}\\
=&\lim\limits_{R\to \infty}\| (1-\chi_R)e^{\pm itH}e^{-ix_3\cdot P_1}e^{-ix_3\cdot P_2}[\tilde{\psi}_{\epsilon}]_{bs}(x_1,x_2,\hat{q}_3)e^{ix_3\cdot q_3}\|_{\s^2(X_{q_3})}\\
=&\lim\limits_{R\to \infty}\| (1-\chi_R)e^{-ix_3\cdot P_1}e^{-ix_3\cdot P_2}e^{\pm it\tilde{H}(P_3,3)}[\tilde{\psi}_{\epsilon}]_{bs}(x_1,x_2,\hat{q}_3)e^{ix_3\cdot q_3}\|_{\s^2(X_{q_3})}\\
=&\lim\limits_{R\to \infty}\| (1-\tilde{\chi}_R(x_1,x_2))e^{\pm it\tilde{H}(q_3,3)}[\tilde{\psi}_{\epsilon}]_{bs}(x_1,x_2,\hat{q}_3)\|_{\s^2_{x_1,x_2}(\mathbb{R}^6)}\\
=&0
\end{split}
\eeq
uniformly in $0\leq t<\infty$, and
\eq
\begin{split}
& \lim\limits_{t\to \infty} t^{-1}\int_0^tds \| \chi_R e^{\pm itH} [\psi_\epsilon]_{sc}\|_{\s^2(X_{q_3})}\\
=&\lim\limits_{t\to \infty} t^{-1}\int_0^tds \|\chi_Re^{\pm itH}e^{-ix_3\cdot P_1}e^{-ix_3\cdot P_2}[\tilde{\psi}_{\epsilon}]_{sc}(x_1,x_2,\hat{q}_3)e^{ix_3\cdot q_3}\|_{\s^2(X_{q_3})}^2\\
=&\lim\limits_{t\to \infty} t^{-1}\int_0^tds \|\chi_Re^{-ix_3\cdot P_1}e^{-ix_3\cdot P_2}e^{\pm it\tilde{H}(P_3,3)}[\tilde{\psi}_{\epsilon}]_{sc}(x_1,x_2,\hat{q}_3)e^{ix_3\cdot q_3}\|_{\s^2(X_{q_3})}^2\\
=&\lim\limits_{t\to \infty} t^{-1}\int_0^tds \|\tilde{\chi}_Re^{\pm is\tilde{H}(q_3,3)}[\tilde{\psi}_{\epsilon}]_{sc}(x_1,x_2,\hat{q}_3)e^{ix_3\cdot q_3}\|_{\s^2_{x_1,x_2}(\mathbb{R}^6)}^2\\
=&0
\end{split}
\eeq
for any $R<\infty$. Set 
\eq
a_R(q_3,t):=\| (1-\chi_R)e^{\pm itH}[\psi_\epsilon]_{bs}\|_{\s^2(X_{q_3})}
\eeq
and 
\eq
b_t(q_3,R):=t^{-1}\int_0^tds \|\chi_R e^{\pm itH} [\psi_\epsilon]_{sc}\|_{\s^2(X_{q_3})}^2.
\eeq
Then $\lim\limits_{R\to \infty} a_R(q_3,t)=0$ for all $q_3\in \mathbb{R}^3, t\geq 0$, and $\lim\limits_{t\to \infty} b_t(q_3,R)=0$ for all $q_3\in \mathbb{R}^3, R\geq 0$. Hence, for all $q_3\in \mathbb{R}^3, t\geq 0$, 
\eq
|a_R(q_3,t)|^2\to 0 \quad \text{ as }R\to \infty
\eeq
and 
\eq
|a_R(q_3,t)|^2\leq \| \tilde{\psi}_\epsilon (x_1,x_2,\hat{q}_3)\|_{\s^2_{x_1,x_2}(\mathbb{R}^6)}^2 \in \s^1_{q_3}(\mathbb{R}^3).
\eeq
Therefore, due to Plancherel Theorem and Dominated Convergence Theorem, we obtain that for all $t\geq 0$, 
\eq
\begin{split}
&\| (1-\chi_R)e^{\pm itH} (1-P_c(P_1+P_2+P_3))\psi_\epsilon\|_{\s^2_x(\mathbb{R}^9)}\\
=&\| (1-\tilde{\chi}_R)e^{\pm it\tilde{H}(P_3,3)} (1-P_c(P_3))\tilde{\psi}_\epsilon\|_{\s^2_x(\mathbb{R}^9)}\\
=&c\| \|(1-\tilde{\chi}_R)e^{\pm it\tilde{H}(q_3,3)} (1-P_c(q_3))\tilde{\psi}_\epsilon(x_1,x_2,\hat{q}_3) \|_{\s^2_{x_1,x_2}(\mathbb{R}^6)}\|_{\s^2_{q_3}(\mathbb{R}^3)}\\
=&c\| a_R(q_3,t)\|_{\s^2_{q_3}(\mathbb{R}^3)}\\
\to & 0
\end{split}
\eeq
as $R\to \infty$.
Therefore, we obtain \eqref{Feb.22.1}. Similarly, since for all $q_3\in \mathbb{R}^3, R\geq 0$, 
\eq
b_t(q_3,R)\to 0 \quad \text{ as }t\to \infty
\eeq
and 
\eq
|b_t(q_3,R)| \leq \| \tilde{\psi}_\epsilon (x_1,x_2,\hat{q}_3)\|_{\s^2_{x_1,x_2}(\mathbb{R}^6)}^2 \in \s^1_{q_3}(\mathbb{R}^3),
\eeq
by utilizing the Plancherel Theorem, Fubini's Theorem, and Dominated Convergence Theorem, we have that for all $R\geq 0$,
\eq
 \begin{split}
 & t^{-1}\int_0^t ds \| \chi_R e^{\pm i sH} P_c(P_1+P_2+P_3)\psi\|_{\s^2_x(\mathbb{R}^9)}^2\\
 =&t^{-1}\int_0^tds \|\tilde{\chi}_Re^{\pm is\tilde{H}(P_3,3)} P_c(P_3)\psi_{\epsilon}\|_{\s^2_x(\mathbb{R}^9)}^2\\
 =&ct^{-1}\int_0^t ds \| \|\tilde{\chi}_Re^{\pm is\tilde{H}(q_3,3)} P_c(q_3)\psi_{\epsilon}(x_1,x_2,\hat{q}_3)\|_{\s^2_{x_1,x_2}(\mathbb{R}^6)}\|_{\s^2_{q_3}(\mathbb{R}^3)}^2\\
 =&c\int d^3q_3 b_t(q_3,R)\\
 \to & 0
 \end{split}
\eeq
as $t\to \infty$, where $c$ is a positive constant. Therefore, we obtain \eqref{Feb.22.2}. We finish the proof.
\end{proof}
 \subsection{Threshold and Feynman-Hellmann Theorem}\label{sec: threshold}
In this section, we will introduce the concept of a threshold point in a two/three-body Hamiltonian and discuss the utilization of the Feynman-Hellmann Theorem. Generally speaking, when the spectrum of a Hamiltonian H is near the threshold points, the dispersive behavior of the operator $e^{-itH}$ is abnormal. The Feynman-Hellmann Theorem provides us with a continuity property for these threshold points.
 \subsubsection{The Threshold set of one/two-body Hamiltonians}\label{section: curve} 
 For a threshold point of a Hamiltonian, a similar notion exists in the case of a one-body Schr\"odinger operator, $H=-\Delta_x+V(x)$, where, for example, $\langle x\rangle^\sigma V(x)\in \s^\infty_x(\R^n)$ holds for some $\sigma>3$. $0$ is considered as a threshold point. This notion for Schr\"odinger operators comes from the observation that when the spectrum of $H$ is near $0$, the dispersive estimate of $e^{itH}$ exhibits abnormal behavior due to the possible presence of the zero eigenvalue, the zero resonance for $H$ or both, which can be considered as an exceptional event. See, for example, \cites{JK1979}, which provides a detailed explanation of the zero threshold point for Schr\"odinger operators. Here, we refer to $f$ as a resonance of a one-body Schr\"odinger operator $H=-\Delta_x+V(x)$ if $f$ is the distributional solution to $H\psi(x)=0$ with $ f\notin \s^2_x(\R^3)$ while $\langle x\rangle^{-\sigma }f\in \s^2_x(\R^3)$ for some $\sigma>1/2$. See, for example, \cites{JK1979} for the concept of resonance. In the case of free Schr\"odinger evolution, the following relationship holds:
 \eq
 \| e^{-it(-\Delta_x)}\|_{\s^{1}_x(\R^3)\to \s^{\infty}_x(\R^3)}\lesssim\frac{1}{|t|^{3/2}}.
 \eeq
However, when $e^{itH}$ behaves "badly" around $0$ due to $0$ being either a resonance or an eigenvalue of $H$, we observe only a point-wise decay of $1/\sqrt{|t|}$: for some initial data $u_0(x)$, which decays sufficiently rapidly at infinity, the solution $u(t,x)=e^{-itH}u_0(x)$ admits as $t\to \infty$ an asymptotic expansion
 \eq
 u(t,x)=\sum\limits_{\text{finite}} a_je^{i\lambda_jt}\phi_j(x)+\frac{1}{\sqrt{|t|}}B_0u_0(x)+\frac{1}{|t|^{3/2}}B_1u_0(x)+\cdots,
 \eeq
 which is valid locally in space. Here $\phi_j(x)$ are eigenfunctions of $H$ with eigenvalues $\lambda_j$, and $B_0\neq 0$ if $0$ is a resonance or an eigenvalue of $H$. The $B_j, j=0,1,\cdots$ are finite-rank operators. For more detailed information and additional references, please refer to \cite{GJY2004} (page 233) and the relevant literature cited therein. $0$ holds a special significance for the Schr\"odinger operator due to the fact that the group velocity of $-\Delta_x$ is given by $2P:=2(-i)\nabla_x$, and $0$ corresponds to the energy where $2P=0$. In the context of general one-/two-body Hamiltonians, we define a \textbf{threshold point} in a similar way. In the case of a three-body Hamiltonian $H$, we define a point $\lambda\in \R$ as a threshold point if the behavior of $e^{-itH}$ becomes "abnormal" when the spectrum of $H$ is in the vicinity of $\lambda$.

Let us start with the one/two-body case. The threshold points are defined based on the concepts of group velocity and relative group velocity of two particles. In quantum mechanics, velocity is characterized by the concept of group velocity.
\begin{definition}[Group velocity of a particle]
Let $\omega_j\in C^1(\mathbb{R}^3)$ for $j=1,2,3$. The group velocity of particle $j$ is defined as
\eq
v_j(P_j):= \nabla_{P_j}[\omega_j(P_j)],\quad P_j:=-i\nabla_{x_j}.
\eeq
\end{definition}
\begin{definition}[Relative group velocity of two particles]
The relative group velocity of particles $j$ and $l$ with $1\leq j<l\leq 3$ is defined as
\eq
v_{jl}(P_j,P_l):= v_j(P_j)-v_l(P_l),\quad P_j:=-i\nabla_{x_j}, P_l:=-i\nabla_{x_l},
\eeq
where $v_j$ and $v_l$ are the group velocities of particles $j$ and $l$, respectively.
\end{definition}

Now let us look at the relationship between two-body Hamiltonians and one-body Hamiltonians. In this section, we assume that $(jkl)$ is a permutation of $(123)$ with $j<k$ fixed. Let $a=(jk)(l)\in L$. Recall that 
\eq
H^a:=\omega_j(P_j)+\omega_k(P_k)+V_{jk}(x_j-x_k)
\eeq
represents a two-body Hamiltonian. We introduce the following definition:
\eq
H^a(\eta):=\omega_{j}(\eta-P_k)+\omega_k(P_k)+V_{jk}(-x_k)\quad \text{ for }\eta\in \mathbb{R}^3.\label{Hab}
\eeq
We can observe the relationship between $H^a$ and $H^a(\eta)$ through the following identity:
\eq
H^a=e^{-ix_j\cdot P_k}H^a(P_j)e^{ix_j\cdot P_k}.\label{Hab2}
\eeq
In \eqref{Hab2}, the appearance of $P_j$ in $H^a(P_j)$ on the right-hand side corresponds to the total quasi-momentum of particle $j$ and particle $k$, which is $P_j+P_k$ in $H^a$, since we have
\eq
P_j+P_k=e^{-ix_j\cdot P_k}P_je^{ix_j\cdot P_k},
\eeq
and
\eq
F(|P_j+P_k-v|\leq \epsilon)=e^{-ix_j\cdot P_k}F(|P_j-v|\leq \epsilon)e^{ix_j\cdot P_k},\quad \text{for all }v\in \R^3, \epsilon>0.
\eeq
Here $F$ denotes a smooth cut-off function. When studying the two-body problem, it is common to simplify it by reducing it to a one-body problem using equation \eqref{Hab2}. The notion of threshold points for a two-body Hamiltonian $H^a$ is established within the framework of fixing the total quasi-momentum. In other words, the total quasi-momentum needs to be fixed prior to discussing the concept of a threshold point for the two-body Hamiltonian $H^a$. Here is the definition of a threshold point for a one-body Hamiltonian $H^a(\eta)$, where $a=(jk)(l)\in L$ and $\eta\in \R^3$:
\begin{definition}[Threshold points for $H^a(\eta)$]\label{defHa1} Given $\eta\in \R^3$ and $a=(jk)(l)\in L$, a threshold point $\lambda^a(\eta)$ of $H^a(\eta)$ occurs if 
\eq
\lambda^a(\eta)=\omega_j(\eta-q_k)+\omega_k(q_k)
\eeq
with $q_k$ satisfying 
\eq
\nabla_{p_k}[ \omega_j(\eta-p_k)+\omega_k(p_k)]\vert_{p_k=q_k}=0.
\eeq
    
\end{definition}
\begin{definition}[Threshold points for $H^a$]\label{defHa2} Given $v\in \R^3$ and $a=(jk)(l)\in L$, a threshold point $\lambda^a(v)$ of $H^a$, with a total quasi-momentum $v$, occurs if $\lambda^a(v)$ is a threshold point of $H^a(v)$. 
    
\end{definition}
According to Definitions \ref{defHa1} and \ref{defHa2}, the collection of threshold points of $H^a(v)$ (where $v\in \R^3$) is equal to the set of threshold points of $H^a$ with a fixed total quasi-momentum $v$. Let $\tau^a(v)$ denote the set of all threshold points for $H^a(v)$, where $v\in \R^3$. 

The definition of a threshold point for a three-body Hamiltonian relies on the concept of a generic one-body Hamiltonian and the notion of eigenfunctions of a one-body Hamiltonian. Now, let us proceed by introducing the notion of a generic one-body Hamiltonian. Given $a=(jk)(l)\in L$ and $\eta\in \R^3$, we let
\eq
H^a_0(\eta):=\omega_j(\eta-P_k)+\omega_k(P_k).
\eeq
\begin{definition}[Eigenfunctions/eigenvalues of a one-body Hamiltonian]\label{def: b-1-H} Given $a=(jk)(l)\in L$ with $j<k$ and $\eta\in \R^3$, $\lambda$ is an eigenvalue of $H^a(\eta)$ if 
\eq
H^a(\eta)f(x_k)=\lambda f(x_k)\label{egn}
\eeq
for some $f(x_k)\in \s^2_{x_k}(\R^3)$. We refer to $f(x_k)$ satisfying \eqref{egn} as an eigenfunction of $H^a(\eta)$ with an eigenvalue $\lambda$. 

\end{definition}
When it comes to the concepts of eigenvalues and eigenfunctions of a two-body Hamiltonian $H^a(a=(jk)(l)\in L)$, we have to fix the total quasi-momentum of the two particles, and then delve into these concepts:
\begin{definition}[Eigenvalues of a two-body Hamiltonian]\label{def: b-2-H} Given $a=(jk)(l)\in L$ with $j<k$ and $v\in \R^3$, $\lambda$ is an eigenvalue of $H^a$ with a total quasi-momentum $v$ if $\lambda$ is an eigenvalue of $H^a(v)$. For an eigenvalue $\lambda$, $f(x_k-x_j)$ is an corresponding eigenfunction of $H^a$ with a total quasi-momentum $v$ if 
\eq
H^a(v)f(x_k)=\lambda f(x_k).
\eeq
\end{definition}
{

{It is known, for example, by using resolvent, that if all threshold points are neither an eigenvalue nor a resonance, for a Schr\"odinger-tpye one-body Hamiltonian, we have some nice results:
\begin{proposition}\label{prop: regular}If all threshold points have neither an eigenvalue nor a resonance for a Schr\"odinger-type one-body Hamiltonian, then such Schr\"odinger-type one-body Hamiltonian enjoys: for example, when $H_s=-\Delta_{x_k}+W(x_k)$ for some $\langle x_k\rangle^\sigma W(x_k)\in \s^\infty_{x_k}(\R^3)$ for some $\sigma \geq 3$,
\begin{itemize}
\item[(a)] There are at most finitely many eigenfunctions and eigenvalues. 
\item[(b)] All of eigenfunctions are localized in space: for all eigenfunctions of $H_s$, $f_e(x_k)$, $f_e(x_k)$ satisfy, at least,
\eq
\|\langle x_k\rangle^{\sigma}f_e(x_k)\|_{\s^2_{x_k}(\R^3)}\lesssim_{\sigma, H_s} \|f_e(x_k)\|_{\s^2_{x_k}(\R^3)}.
\eeq
\item[(c)] $e^{itH_s}$ satisfies $\s^p$ decay estimates:
\eq
\|e^{itH_s}P_c(H_s)f(x_k)\|_{\s^p_{x_k}(\R^3)}\lesssim_{p,H_s} \frac{1}{|t|^{3(\frac{1}{2}-\frac{1}{p})}}\|f(x_k)\|_{\s^{p'}_{x_k}(\R^3)},
\eeq
for all $f(x_k)\in \s^{p'}_{x_k}(\R^3)\cap \s^2_{x_k}(\R^3)$ and all $p\in [2,\infty]$, where $P_c(H_s)$ denotes the projection on the continuous spectrum of $H_s$ and 
$$\frac{1}{p}+\frac{1}{p'}=1.$$
\item[(d)]$e^{itH_s}$ satisfies the local decay property: for all $f(x_k)\in \s^2_{x_k}(\R^3)$ and $\epsilon>0,$, for $\lambda non-threshold$
\eq
\| \langle x_k\rangle^{-1-\epsilon}\frac{1}{H_s-\lambda \pm i0}f(x_k)\|_{\s^2_{x_k,\lambda}(\R^{3+1})}\lesssim_\epsilon \| f(x_k)\|_{\s^2_{x_k}(\R^3)}.
\eeq
\end{itemize}
\end{proposition}
Assumption \ref{asp: subH} ensures that all sub-Hamiltonians $H^a$ with $a=(jk)(l)\in L$ and $H_0$ share 
 similar properties outlined in Proposition \ref{prop: regular}, provided that their energy is away from finitely many points in each fiber space. These points will be used to define threshold points for our three-body Hamiltonian $H$, which will be elaborated upon in section \ref{sec: FH}. 


}

\subsubsection{Feynman-Hellmann Theorem}\label{sec: FH}
The discussion of this part is based on the following assumption on Assumptions \ref{asp: short}- \ref{asp: subH}. Based on Assumption \ref{asp: short} and Theorem XIII.6 on page 87 in \cite{RS1978}, it can be inferred that for every $a=(jk)(l)\in L, \eta\in \R^3$, the operator $H^a(\eta)$ has a finite number of eigenvalues at most, denoted as $\lambda_1\leq \cdots\leq \lambda_m$ (where $\lambda_r, r=1,\cdots,m,$ may be identical to one another), with $m\leq N$ for some constant $N$ which is independent on $\eta\in \R^3$. 

The Feynman-Hellmann Theorem elucidates how these eigenvalues change as $\eta\in \R^3$ varies. In simple terms, the Feynman-Hellmann Theorem supplies a formula for the derivative of the eigenvalues with respect to the parameter. 
\begin{theorem}[Feynman-Hellmann Theorem]\label{Thm: FH}For $ I=(a_1,b_1)\times \cdots (a_m,b_m)\subseteq \R^m$, let $\{H(p)\}_{p\in I}$ be a class of self-adjoint operators on $\s^2_x(\R^n)$. {Assume that $\lambda(p)$ is an eigenvalue of $H(p)$. If $H(p)$ is differentiable with respect to $p$,} then $\lambda(p)$ is differentiable with respect to $p$ at $p=p_0$:
\eq
\frac{\partial\lambda(p)}{\partial p_k}\vert_{p=p_0}=(\psi_\lambda, \frac{\partial H(p)}{\partial p_k}\psi_{\lambda})_{\s^2_x(\R^n)}\vert_{p=p_0}, \quad k=1,\cdots, m\label{con:FH}
\eeq
where $\psi_\lambda$ standards for a normalized eigenfunction $(\|\psi_\lambda\|_{\s^2_x(\R^n)}=1)$ with an eigenvalue $\lambda(p).$

\end{theorem}
When we talk about the derivative of an eigenvalue $\lambda(p)$ with respect to the parameter $p\in \R^m$ at $p=p_0$ for some $p_0\in \R^m$, it means that $\lambda(p)$ exists in a neighborhood of $p_0$ in $\R^m$. Therefore, the conclusion of Feynman-Hellmann Theorem holds on the premise that this eigenvalue $\lambda(p)$ exists in a neighborhood of $p_0$. However, in some situation, $\lambda(p)$ may \textbf{disappear} at $p=p'$ for some $p'\in \R^m$. We can see such phenomenon in the following example:
\begin{example}Assume that: 
\begin{enumerate}
\item $V_{jk}(-x_k)\in C_0^\infty(\R^3)$ and $V_{jk}(-x_k)\leq 0$ for all $x_k\in \R^3$;
\item $-\Delta_{x_k}+V_{jk}(-x_k)$ has an eigenvalue $\lambda=-1$ on $\s^2_{x_k}(\R^3)$. 
\end{enumerate}
If $H^a(\eta)=(2-|\eta|^2)(-\Delta_{x_k})+V_{jk}(-x_k)$, then when $|\eta|=1$, $H^a(\eta)$ has one eigenvalue $\lambda=-1$. When $|\eta|=2$, $-H^a(\eta)=2(-\Delta_{x_k})-V_{jk}(-x_k) $ with $-V_{jk}(-x_k)\geq 0$. Since a Schr\"odinger operator with a well-localized potential does not have embedded eigenvalues (for example, see Theorem XIII.58 in \cite{RS4}), then $H^a(\eta)$ does not have eigenvalues when $|\eta|=2$. Therefore, the eigenvalue of $H^a(\eta)$ \textbf{disappears} at $\eta=\eta_0$ for some $\eta_0$ with $|\eta_0|\in (1,2)$. 
    
\end{example}
\begin{proof}[Proof of Lemma \ref{Lem: dis}] It is sufficient to show that if $lambda$ is a non-embedded eigenvalue, then it will not disappear. Assume that $\lambda$ is a non-embeded eigenvalue of $H^a(\eta)$. Let
\eq
\mathscr{E}^a:=\min\limits_{q\in \R^3}\omega_j(\eta-q)+\omega_k(q),
\eeq
be minimum energy of $H^a_0(\eta):=\omega_j(\eta-P_k)+\omega_k(P_k)$. Then
\eq
c= \mathscr{E}^a- \lambda>0. 
\eeq
Then if we can show that operator \eq
F(H^a(\eta), c):=F(H^a(\eta)<\mathscr{E}^a-\frac{c}{10} )F(H^a(\eta)> \mathscr{E}^a-10c)
\eeq
is continuous in $\eta$ on $\s^2_{x_k}(\R^3)$, then we are done because when $M\neq N$ and $N,M\in \mathbb{N}^+$, a rank $N$ projection can never become a rank $M$ projection continuously. The continuity of $F(H^a(\eta),c)$ on $\s^2_{x_k}(\R^3)$ follows from the following steps:
\begin{itemize}
\item[a)] Let 
\eq
F_1(H^a(\eta),c):=F(H^a(\eta)<4\mathscr{E}^a-\frac{2c}{5} )F(H^a(\eta)> \frac{\mathscr{E}^a}{4}-\frac{5c}{2})
\eeq
and write 
\eq
F(H^a(\eta),c)=F(H^a(\eta),c)F_1(H^a(\eta),c);
\eeq
    \item[b)] Use Fourier transform to represent $F(H^a(\eta),c)$ and $F_1(H^a(\eta),c)$:
\eq
F(H^a(\eta),c)= \int_{\R}dw e^{iw H}\hat{F}(w),\label{c: F1id}
\eeq
and 
\eq
F_1(H^a(\eta),c)= \int_{\R}dw e^{iw H}\hat{F}_1(w),\label{c: F2id}
\eeq
where $\hat{F}(w)$ and $\hat{F}_1(w)$ denotes the Fourier transform of $F(y,c)$ in $y$ variable. Therefore, we have $|w|\hat{F}(w),|w|\hat{F}_1(w)\in \s^1_w(\R)$;
\item[c)] By using the definition of derivative of parameterized operators, we have that for all $\omega\in \R, d=1,2,3,$
\begin{align}
&\p_{\eta_d}[e^{iwH^a(\eta)}]F(H^a(\eta)<4\mathscr{E}^a-\frac{2c}{5} )\nonumber\\
=&i \int_0^w ds e^{i(w-s)H^a(\eta)} \p_{\eta_d}[\omega_j(\eta-P_k)]e^{isH^a(\eta)}F(H^a(\eta)<4\mathscr{E}^a-\frac{2c}{5} )\label{collapse: eq1}
\end{align}
and 
\begin{align}
&F(H^a(\eta)<\mathscr{E}^a-\frac{c}{10} )\p_{\eta_d}[e^{iwH^a(\eta)}]\nonumber\\
=&i \int_0^w ds F(H^a(\eta)<\mathscr{E}^a-\frac{c}{10} )e^{i(w-s)H^a(\eta)} \p_{\eta_d}[\omega_j(\eta-P_k)]e^{isH^a(\eta)}\label{collapse: eq2}
\end{align}
on $\s^2_{x_k}(\R^3)$. By using Assumption \ref{asp: subH}: for all $\omega\in \R$,
\eq
\| \nabla_\eta[\omega_j(\eta-P_k)]e^{isH^a(\eta)} F(H^a(\eta)<4\mathscr{E}^a-\frac{2c}{5} )\|_{\s^2_{x_k}(\R^3)\to \s^2_{x_k}(\R^3)}\lesssim_{\mathscr{E}^a}1
\eeq
and 
\eq
\| F(H^a(\eta)<\mathscr{E}^a-\frac{c}{10} )\nabla_\eta[\omega_j(\eta-P_k)]e^{isH^a(\eta)} \|_{\s^2_{x_k}(\R^3)\to \s^2_{x_k}(\R^3)}\lesssim_{\mathscr{E}^a}1,
\eeq
due to \eqref{c: F1id}, \eqref{c: F2id}, \eqref{collapse: eq1} and \eqref{collapse: eq2}, we conclude that
\eq
\|\p_{\eta_d}[F(H^a(\eta),c)]F(H^a(\eta)<4\mathscr{E}^a-\frac{2c}{5} ) \|_{\s^2_{x_k}(\R^3)\to \s^2_{x_k}(\R^3)}\lesssim_{\mathscr{E}^a} \| w\hat{F}(w)\|_{\s^1_w(\R)}
\eeq
and 
\eq
\|F(H^a(\eta)<\mathscr{E}^a-\frac{c}{10} )\p_{\eta_d}[F_1(H^a(\eta),c)] \|_{\s^2_{x_k}(\R^3)\to \s^2_{x_k}(\R^3)}\lesssim_{\mathscr{E}^a} \| w\hat{F}_1(w)\|_{\s^1_w(\R)}.
\eeq
\item[d)]
Now we can calculate $\p_{\eta_d}[ F(H^a(\eta),c)], d=1,2,3$: by using \eqref{collapse: eq1} and \eqref{collapse: eq2}, we have
\begin{align}
    \p_{\eta_d}[F(H^a(\eta),c)]=&\p_{\eta_d}[F(H^a(\eta),c)F_1(H^a(\eta),c)]  \nonumber\\
    =& \p_{\eta_d}[F(H^a(\eta),c)]F_1(H^a(\eta),c)+F(H^a(\eta),c)\p_{\eta_d}[F_1(H^a(\eta),c)]
\end{align}
and $ \p_{\eta_d}[F(H^a(\eta),c)]$ is bounded on $\s^2_{x_k}(\R^3)$ since 
\begin{align}
    &\| \p_{\eta_d}[F(H^a(\eta),c)]\|_{\s^2_{x_k}(\R^3)\to \s^2_{x_k}(\R^3)}\nonumber\\
    \leq& \|\p_{\eta_d}[F(H^a(\eta),c)]F(H^a(\eta)<4\mathscr{E}^a-\frac{2c}{5} ) \|_{\s^2_{x_k}(\R^3)\to \s^2_{x_k}(\R^3)}\nonumber\\
    &+\|F(H^a(\eta)<\mathscr{E}^a-\frac{c}{10} )\p_{\eta_d}[F_1(H^a(\eta),c)] \|_{\s^2_{x_k}(\R^3)\to \s^2_{x_k}(\R^3)}\nonumber\\
    \lesssim_{\mathscr{E}^a} &1.
\end{align}

\end{itemize}
We finish the proof.  
\end{proof}

Let $\{H(p)\}_{p\in I}$ be as in Theorem \ref{Thm: FH}. In this note, we define the Feynman-Hellmann Theorem as being \textbf{applicable} to $\{ H(p)\}_{p\in I}$ at $p=p_0$ if all eigenvalues of $H(p)$ are defined in a neighborhood of $p_0$ in $\R^m$. Otherwise, we say that the eigenvalues of $\{H(p)\}_{p\in I}$ \textbf{collapse} at $p=p_0$. See the following definition for our $H^a(\eta), a=(jk)(l)\in L, \eta\in \R^3$:
\begin{definition}[Collapse of the eigenvalues of $H^a(\eta)$]\label{def: coll-eigen}For a given $a=(jk)(l)\in L$ and $\eta_0\in \R^3$, we say that the eigenvalues of $H^a(\eta)$ do not \textbf{collapse} at $\eta=\eta_0$ if Feynman-Hellmann Theorem is applicable to $\{ H^a(\eta)\}_{\eta\in \R^3}$ at $\eta=\eta_0$. Otherwise, we say that the eigenvalues of $H^a(\eta)$ \textbf{collapse} at $\eta=\eta_0$.

    
\end{definition}
\begin{remark}Establishing the validity of \eqref{con:FH} is crucial. When using Duhamel's formula to expand $e^{-itH}$ with respect to $e^{-itH_a}$ for some $a=(jk)(l)\in L$, it becomes necessary to have dispersive estimates of $e^{-itH^a}$ to comprehend the Duhamel term. For the component involving $P_b(H^a)e^{-itH_a}$, we cannot derive any dispersion from $H^a$ or from the $x_j-x_k$ direction. If the quasi-momentum $P_j+P_k=v_{jk}\in \R^3$ is fixed, with $P_b(H^a)$, $H^a$ becomes a number dependent on $v_{jk}$. When the total quasi-momentum of three particles is fixed, $v_{jk}$ can be regarded as a unique variable in the space of quasi-momentum. To achieve decay in $t$, we perform integration by parts with respect to the $v_{jk}$ variable. This necessitates understanding the group velocity of $\omega_k(v_T-(P_j+P_k))+\lambda(P_j+P_k)$ with respect to $v_{jk}$, where $\sum\limits_{d=1}^3 P_d=v_T\in \R^3$ is fixed and $\lambda(P_j+P_k)$ represents an eigenvalue of $H^a$.  
    
\end{remark}
\begin{remark}Let $P_{b,d}(H^a(\eta))$ be defined as in \eqref{def: Na} for $d=1,\cdots,N^a$. Given \eqref{con:FH}, if $P_b(H^a(\eta))$ is differentiable at $\eta=\eta_0$, it allows us to compute higher-order derivatives of the eigenvalues of $H^a$ with respect to the parameter. We will see it later in Section \ref{sec 3}. 
    
\end{remark}
Fortunately, with Lemma \ref{Lem: dis}, all eigenvalues of $H^a(\eta)$ do not collapse.

Now let us talk about the {overlapping eigenvalues}. When Feynman-Hellmann Theorem \ref{Thm: FH} is applicable, it tells us that the eigenvalues of $H^a(\eta) (a=(jk)(l)\in L, \eta\in \R^3)$ varies continuously as $\eta\in \R^3$ varies. Since the dimension of the space of all eigenfunctions, with an eigenvalue $\lambda(p)$, may be greater than $1$, and since different $\psi_\lambda$ may have a different
\eq
(\psi_\lambda, \frac{\partial H(p)}{\partial p_k}\psi_{\lambda})_{\s^2_x(\R^n)},
\eeq
it is normal that several eigenvalues may emerge into one eigenvalue and one eigenvalue may split into several eigenvalues. This may cause that the minimum gap between two different eigenvalues goes to $0$ at $\eta=\eta_0$ for some $\eta_0\in \R^3$, which means when an eigenvalue $\lambda(|\eta|)$ is not an embedded eigenvalue, it is improper to use $F(|H^a(\eta) -\lambda(|\eta|)|<c), c>0$ to denote the projection on the space of all eigenfunctions of $H^a(\eta)$. 
\begin{definition}[overlapping eigenvalues]For a given $a=(jk)(l)\in L$ and $\eta_0\in \R^3$, we say that $\lambda(\eta_0)$ is an overlapping eigenvalue of $H^a(\eta_0)$ if its multiplicity, the dimension of the space of all eigenfunctions with an eigenvalue $\lambda(\eta_0)$, is greater than or equal to $2$. 
    
\end{definition}

\begin{definition}[collapse of overlapping eigenvalues]For a given $a=(jk)(l)\in L$ and $\eta_0\in \R^3$, assume that $H^a(\eta_0)$ has an overlapping eigenvalue $\lambda(\eta_0)$. Let $\tilde{N}(H^a(\eta),\lambda(\eta))$ be the multiplicity of $\lambda(\eta)$. We say that the overlapping eigenvalue $\lambda(\eta_0)$ \textbf{collapses} at $\eta_0$ if $\tilde{N}(H^a(\eta),\lambda(\eta))$ is not continuous at $\eta=\eta_0$.

\end{definition}
\begin{definition}[collapse of $H^a(\eta)$] For $a=(jk)(l)\in L$ and $\eta\in \R^3$, we say that $H^a(\eta)$ \textbf{collapses} at $\eta=\eta_0$ it meets any of the following conditions:
\begin{enumerate}
\item 
The eigenvalues of $H^a(\eta)$ collapse at $\eta=\eta_0$; 
\item An overlapping eigenvalue of $H^a(\eta)$ does the same at $\eta=\eta_0$.
\end{enumerate}
\end{definition}
According to Assumptions \ref{asp: short} and \ref{asp: critical b}, and Lemma \ref{Lem: dis}, we deduce that in our case, the eigenvalues of $H^a(\eta)$ will never collapse for all $\eta\in \R^3$.

\subsubsection{ The Threshold set of three-body Hamiltonians}
Now let us discuss the concept of a threshold point in a three-body Hamiltonian. Similar to the concept of a threshold point in a two-body Hamiltonian, we first fix the total quasi-momentum of the three particles, and then delve into the concept of a threshold point in a three-body Hamiltonian with a given total quasi-momentum.
\begin{definition}[Threshold points for $H_a$]\label{TPHa} For $a=(jk)(l)\in L$, let $I_a$ be as in Assumption \ref{asp: subH}. Given $v_T\in \mathbb{R}^3$ and $a=(jk)(l)\in L$, a threshold point $\lambda(v_T)$ of $H_a=\sum\limits_{m=1}^3 \omega_m(P_m)+V_{jk}(x_j-x_l)$ with a total quasi-momentum $v_T$ occurs if and only if $\lambda(v_T)\in \tau_a(v_T):=\bigcup\limits_{m=1}^3 \tau_{a,m}(v_T)$ where $\tau_{a,m}(v_T),m=1,2,3,$ are defined by: 
 with $\sum\limits_{m=1}^3\eta_m=v_T$,
\begin{enumerate}
\item $\tau_{a,1}(v_T):=\{\infty\}$.

\item $\lambda(v_T)\in \tau_{a,2}(v_T) $ if and only if $\lambda(v_T)=\sum\limits_{m=1}^3 E_m $, where $(E_1,E_2,E_3)$ satisfies that:
\begin{enumerate}
\item $E_l=\omega_l(v_T-\eta_l)$.
\item  \eq
\nabla_{\eta}[\omega_l(v_T-\eta)]+\nabla_\eta[\lambda_{jk}(\eta)]\vert_{\eta=\eta_l}=0.
\eeq

\end{enumerate}
\item Let $\{I_a(p)\}_{p\in \R^3}$ be as in Assumption \ref{asp: subH}. $\tau_{a,3}(v_T)= I_a(v_T)$.

\end{enumerate}

\end{definition} 
\begin{lemma}\label{Lem: tau finite}{If Assumptions \ref{asp: Hthreshold} and \ref{asp: subH} are valid, then there exists a positive integer $N$ such that $|\tau_a(v_T)|\leq N$ for all $v_T\in \R^3$ and all $a=(jk)(l)\in L$.}
    
\end{lemma}
\begin{proof}Take $v_T\in \R^3-\{0\}$. $|\tau_{a,1}(v_T)|=1$. By employing Assumption \ref{asp: Hthreshold}, there exists a positive integer $N_1$ such that $|\tau_{a,2}(v_T)|\leq N_1$ for all $v_T\in \R^3$. For $\tau_{a,3}(v_T)$, by employing Assumption \ref{asp: subH}, there exists a positive integer $N_2$ such that $|I_a|\leq N_2$ for all $a=(jk)(l)\in L$. Thus, for $v_T\in \R^3$, 
\begin{align*}
|\tau_a(v_T)|\leq& \sum\limits_{b=1}^3 |\tau_{a,b}(v_T)|\\
\leq& 1+N_1+N_2. 
\end{align*}
We finish the proof. 
    
\end{proof}

    
    
    

\begin{definition}[Threshold points for $H$]\label{TP}Given $v_T\in \mathbb{R}^3$, a threshold point $\lambda(v_T)$ of $H=H_0+V(x)$ with a total quasi-momentum $v_T$ occurs if and only if $\lambda(v_T)\in \tau(v_T):=\bigcup\limits_{m=1}^3 \tau_{m}(v_T)$ where $\tau_{1}(v_T)$ and $\tau_{2}(v_T)$ are defined by: 
\begin{enumerate}
\item $\tau_1(v_T):=\bigcup\limits_{a=(jk)(l)\in L} \tau_{a}(v_T).$
\item $\lambda(v_T)\in \tau_{2}(v_T) $ if and only if $\lambda(v_T)=\sum\limits_{j=1}^3\omega_j(\eta_j)$ where $\eta_j\in \mathbb{R}^3$ satisfy the following conditions:
\begin{enumerate}
\item $v_{j}(P_j)\vert_{P_j=\eta_j}=v_{l}(P_l)\vert_{P_l=\eta_l}=0\quad \text{ for some }1\leq j<l\leq 3$.\label{TPcon10}
\item $\sum\limits_{j=1}^3 \eta_j=v_T$.\label{TPcon20}
\end{enumerate}
\item $\lambda(v_T)\in \tau_{3}(v_T) $ if and only if $\lambda(v_T)=\sum\limits_{j=1}^3\omega_j(\eta_j)$ where $\eta_j\in \mathbb{R}^3$ satisfy the following conditions:
\begin{enumerate}
\item $v_{1}(\eta_1)=v_{2}(\eta_2)=v_{3}(\eta_3).$
\item $\sum\limits_{j=1}^3 \eta_j=v_T$.
\end{enumerate}
    \end{enumerate}
\end{definition}
Now let us prove that all non-embedded eigenvalues of $H^a(\eta)$ for $\eta\in \R^3$ and $a=(jk)(l)\in L$ will not disappear.
\begin{proof}[Proof of Lemma \ref{Lem: tauf}] By referring to Assumptions \ref{asp: Hthreshold} and \ref{asp: subH}, as well as Lemma \ref{Lem: tau finite}, it follows that $|\tau_1(v_T)|\leq N$ for some constant $N>0$. By using Assumption \ref{asp: critical}, we deduce: 
\begin{align*}
|\tau_2(v_T)|=& \sum\limits_{j=1}^2\sum\limits_{k=j+1}^{3} |\mathscr{C}_j|\times |\mathscr{C}_k|\\
\leq& N_1
\end{align*}
for some constant $N_1>0$. By employing Assumption \ref{asp: threshold}, it's also established that there exists a positive integer $N_2$ such that 
\eq
|\tau_{3}(v_T)|\leq N_2
\eeq
for all $v_T\in \R^3$. Consequently, 
\eq
|\tau(v_T)|\leq N+N_1+N_2
\eeq
and we finish the proof.
    
\end{proof}

Now let us discuss about smooth cut-off functions. Recall that given a finite set $A=\{a_1,\cdots,a_n\}\subseteq \R$ and $y\in \R$, we have
\eq
[y-A]_d:=(y-a_1,\cdots, y-a_n)
\eeq
and for $\epsilon>0$ and $a=(jk)(l)\in L$
\begin{align}
\bar{F}_{\tau}(H,\epsilon):=& (\Pi_{j=1}^{j=n} F(|H-\tau(P_1+P_2+P_3)|>\epsilon)) F(|H|\leq \frac{1}{\epsilon})F(|P_1+P_2+P_3|\leq \frac{1}{\epsilon}),
\end{align}
\eq
\bar{F}_{\tau^a}(H^a,\epsilon):=(\Pi_{j=1}^{j=n} F(|H^a-\tau^a(P_j+P_k)|>\epsilon))F(|H^a|\leq \frac{1}{\epsilon})F(|P_1+P_2+P_3|\leq \frac{1}{\epsilon}), 
\eeq
$F_\tau:=1-\bar{F}_\tau$ and $F_{\tau^a}:=1-\bar{F}_{\tau^a}$.
\begin{lemma}If Assumptions \ref{asp: critical}, \ref{asp: subH}, \ref{asp: Hthreshold} and \ref{asp: threshold} hold, then 
\eq
s\text{-}\lim\limits_{\epsilon \downarrow 0}F_\tau(H,\epsilon)P_{sc}=0,\quad \text{ on }\s^2_x(\R^9).
\eeq
    
\end{lemma}
\begin{proof}By using Assumptions \ref{asp: critical}, \ref{asp: subH}, \ref{asp: Hthreshold} and \ref{asp: threshold}, due to Lemma \ref{Lem: tau finite}, we have that for each $v_T\in \R^3$, $\tau(v_T)=\{a_1(v_T),\cdots,a_{N(v_T)}(v_T)\}$ for some $N(v_T)\in \mathbb{N}$ with $\sup\limits_{v_T\in \R^3} N(v_T)<\infty$. With $P_{sc}=P_c(P_1+P_2+P_3)$, we conclude that for all $v_T\in \R^3$,
\eq
s\text{-}\lim\limits_{\epsilon \downarrow 0}(1-\Pi_{d=1}^{d=N(v_T)}\bar{F}(|H-a_j(v_T)|>\epsilon) )P_c(v_T),\quad \s^2(X_{v_T}).
\eeq
Therefore, by using Dominated Convergence Theorem, we get 
\eq
s\text{-}\lim\limits_{\epsilon \downarrow0} F_\tau(H,\epsilon)P_{sc}=0,\quad \text{ on }\s^2_x(\R^9)
\eeq
and finish the proof.
    
\end{proof}

\begin{lemma}\label{Lem: subdecay}If Assumptions \ref{asp: critical}, \ref{asp: critical b}, \ref{asp: subH}, \ref{asp: Hthreshold} and \ref{asp: threshold} hold, then for all $a=(jk)(l)\in L$ and $\epsilon \in (0,1)$
    \eq
    \| \langle x_j-x_k\rangle^{-2} \bar{F}_\tau(H_a, \epsilon)e^{-itH_a}P_c(H^a)\langle x_j-x_k\rangle^{-2}\|_{\s^2_x(\R^9)\to \s^2_x(\R^9)}\leq \frac{C}{\langle t\rangle^{3/2}}
    \eeq
    for some $C=C(\epsilon, \sup\limits_{p\in \R^3} |\tau(p)|)>0$.
\end{lemma}
\begin{proof} Based on parts (a) of Assumptions \ref{asp: critical} and \ref{asp: critical b}, we can deduce the existence of $E_g\in \R$ such that $\omega_d(\eta)\geq E_g$ for all $d=1,2,3.$ Without loss of generality, we set
\eq
E_g\equiv \min\limits_{d\in \{1,2,,3\}} \min\limits_{\eta\in \R^3} \omega_d(\eta). \label{Eg}
\eeq
Consequently, when $H<E$, we have $H^a< E-E_g$. 

By using Assumptions \ref{asp: critical}, \ref{asp: subH}, \ref{asp: Hthreshold}, and \ref{asp: threshold}, as well as Lemma \ref{Lem: tauf}, we deduce that $\sup\limits_{v_T\in \R^3} |\tau (v_T)|<\infty.$ Let $\tau_c(v_T)\equiv \tau(v_T)-I_a(v_T)$ for $v_T\in \R^3$. Given a finite set $A=\{ a_1,\cdots, a_n\}\subset \R$ and $y\in \R$, we let
\eq
\bar{G}([y-A]_d, \epsilon)\equiv \Pi_{j=1}^{j=n} F(|y-a_j|>\epsilon), \quad \epsilon>0.
\eeq
Then we have
\eq
\bar{F}_\tau(H_a,\epsilon)=\bar{G}([H_a-\tau_c(\sum\limits_{d=1}^3 P_d)],\epsilon) \bar{F}_{I_a}(H_a,\epsilon).
\eeq
Given the continuity of $\sum\limits_{d=1}^3 \omega_d(\eta_d)$ with respect to the variable $(\eta_1,\eta_2,\eta_3)\in \R^9$, we deduce that 
\eq
E_{a,\epsilon}\equiv\sup\limits_{v_T\in \R^3, |v_T|\leq \frac{2}{\epsilon}}\sup\limits_{a\in \tau(v_T)} |a|\lesssim_\epsilon 1.\label{uniformbd}
\eeq
Our desired estimate follows if we can show that for all $b \in \R$
\eq
\|\langle x_j-x_k\rangle^{-2} F(|H_a-b|> \epsilon)\langle x_j-x_k\rangle^2 \|_{\s^2_{x}(\R^9)\to \s^2_x(\R^9)}\lesssim_\epsilon 1.\label{Sep.17.eq1}
\eeq
Now let us prove \eqref{Sep.17.eq1}. By employing inequality  
\begin{align}
    &\|\langle x_j-x_k\rangle^{-2} F(|H_a-b|> \epsilon)\langle x_j-x_k\rangle^2 \|_{\s^2_{x}(\R^9)\to \s^2_x(\R^9)}\nonumber\\
    \leq &1+\|\langle x_j-x_k\rangle^{-2} F(|H_a-b|\leq  \epsilon)\langle x_j-x_k\rangle^2 \|_{\s^2_{x}(\R^9)\to \s^2_x(\R^9)},
\end{align}
it suffices to show 
\eq
\|\langle x_j-x_k\rangle^{-2} F(|H_a-b|\leq \epsilon)\langle x_j-x_k\rangle^2 \|_{\s^2_{x}(\R^9)\to \s^2_x(\R^9)}\lesssim_{\epsilon} 1
\eeq
holds true for all $b\in \R$ with $|b|\leq E_\epsilon\equiv\sup\limits_{a=(jk)(l)\in L}E_{a,\epsilon}$ for some $E=E(\epsilon)>0$. By using Fourier representation to express $F(|H_a-b|\leq \epsilon)$, we deduce that for all $f\in C_0(\R^9)$, some constants $c>0$ and $E'=10E_\epsilon -E_g>0$
\begin{align}
    & \langle x_j-x_k\rangle^{-2} F(|H_a-b|\leq \epsilon)\langle x_j-x_k\rangle^2 f\nonumber\\
    =& \langle x_j-x_k\rangle^{-2} F(|H_a-b|\leq \epsilon)F(H^a\leq E')\langle x_j-x_k\rangle^2 f\nonumber\\
    =& c\int \hat{F}(w)dw \langle x_j-x_k\rangle^{-2}e^{\frac{iw(H_a-b)}{\epsilon}}F(H_a\leq E')\langle x_j-x_k\rangle^{2},
\end{align}
where $\hat{F}$ denotes the Fourier transform of $F(|y|\leq 1)$ in $y$ variable. By using Assumption \ref{asp: critical b}, we conclude: 
\begin{align}
    &\| \langle x_j-x_k\rangle^{-2} F(|H_a-b|\leq \epsilon)\langle x_j-x_k\rangle^2 f\|_{\s^2_x(\R^9)}\nonumber\\
    \lesssim_{E'}&\int dw |\hat{F}(w)| \times\langle \frac{w}{\epsilon}\rangle^2\nonumber\\
    \lesssim_{\epsilon} & 1,
\end{align}
where we use that $E'=E'(E_\epsilon)=E'(\epsilon)$ is a constant dependent on $\epsilon$. We finish the proof.

\end{proof}

\section{Existence of Channel wave operators and Channel projections}\label{Channel}
This section serves as an introduction to channel wave operators, and we present the proof of Theorem \ref{thm: channel} and Proposition \ref{Prop: Pab}.
\subsection{Channel wave operators}
To prove AC, we use induction. Specifically, we reduce the $3$-body problem to several $2$-body problems. We then further reduce each $2$-body system to a one-body system through translation. The one-body problem is well-studied, as demonstrated by  \cite{cycon2009schrodinger}.
\begin{definition}[Channel wave operators] The $H_a$ channel wave operator is defined by 
\eq
\Omega_a^*\psi(0):=s\text{-}\lim\limits_{t\to \infty} e^{itH_a}J_ae^{-itH}P_{sc}\psi(0)\quad \psi(0)\in \mathcal{H}=\s^2_x(\mathbb{R}^9)
\eeq
where $\{ H_a\}$ denotes a set of all some sub-Hamiltonians of $H$ and $\{J_a\}$ is a set of some smooth cut-off functions of $x, P, t$ satisfying
\eq
\sum\limits_{a} J_a=1.
\eeq
\end{definition}
Based on the definition of the channel wave operators, it is evident that the establishment of the existence of all channel wave operators leads to the proof of AC.  Choosing the right $J_a$ is pivotal, and assembling a suitable set of $J_a$ requires ingenuity. In the context that follows, we assume that $\alpha$ satisfies Assumption \ref{asp: omega} and we define $\eta_{jk}$ for $k=1,\cdots, N_{j}$ and $j=1,2,3,$ as given in \eqref{set: threshold}. In this note, we define
\eq
J_{\alpha,a}(t)=e^{-itH_a}F_{c,l,\alpha}(x_l,t,P_l)e^{itH_a}
\eeq
for $a=(jk)(l)\in L$, and 
\eq
J_{\alpha,free}(t):=e^{-itH_0}\left(\Pi_{l=1}^{l=3} F_{c,l,\alpha}(x_l,t,P_l)\right)e^{itH_0},
\eeq
for $a=(1)(2)(3)$. Here, it is worth noting that 
\eq
F_{c,j,\alpha}(x_j,t,P_j)=F_c(|x_j|\leq t^\alpha)\left(\Pi_{k=1}^{N_j}F_1(|P_j-\eta_{jk}|>\frac{1}{t^{\alpha/2}})\right)F_1(|P_j|\leq t^{\alpha/2}),
\eeq
where $\eta_{jk}$ for $k=1,\cdots, N_j,$ are defined in \eqref{set: threshold}. We define the new $3$-body channel wave operators as follows:
\begin{align}
\Omega_{a,\alpha}^{*}:=&s\text{-}\lim\limits_{t\to \infty} e^{itH_a}J_{\alpha,a}(t)e^{-itH}\nonumber\\
=&s\text{-}\lim\limits_{t\to \infty}F_{c,l,\alpha}(x_l,t,P_l) e^{itH_{a}}e^{-itH}\quad \text{ on }\s_x^2(\mathbb{R}^9)
\end{align}
for $a=(jk)(l)\in L$, and we define the new free channel wave operator as follows:
\begin{align}
\Omega_{free,\alpha}^{*}:=&s\text{-}\lim\limits_{t\to \infty} e^{itH_a}J_{\alpha,a}(t)e^{-itH}\nonumber\\
=&s\text{-}\lim\limits_{t\to \infty}\left(\Pi_{l=1}^{l=3} F_{c,l,\alpha}(x_l,t,P_l)\right) e^{itH_{0}}e^{-itH}\quad \text{ on }\s_x^2(\mathbb{R}^9).
\end{align}
Now we prove Theorem \ref{thm: channel}.
\begin{proof}[Proof of Theorem \ref{thm: channel}]We choose $\psi\in \s_x^2(\mathbb{R}^9)$. We will begin with $a=(12)(3)$ and prove the existence of $\Omega_{a,\alpha}^{*}$. Let
\eq
\Omega_{a,\alpha}^{*}(t)=e^{itH_{a}}J_{\alpha,a}(t)e^{-itH}=F_{c,l,\alpha}(x_l,t,P_l) e^{itH_{a}}e^{-itH}
\eeq
and
$$
\tilde{F}_1(P_3,t):=\left(\Pi_{k=1}^{N_3}F_1(|P_3-\eta_{3k}|>\frac{1}{t^{\alpha/2}})\right)F_1(|P_3|\leq t^{\alpha/2}).
$$
For $a=(12)(3)$, by using Cook's method to expand $\Omega_{a,\alpha}^{*}(t)\psi$, we obtain
\begin{align}
\Omega_{a,\alpha}^{*}(t)\psi=&\Omega_{a,\alpha}^{*}(1)\psi+\int_1^tds \partial_{s}[ F_c(|x_3|\leq s^\alpha)]F_1(P_{3},s)e^{isH_{a}}e^{-isH}\psi\nonumber\\
&+\int_1^tds  F_c(|x_3|\leq s^\alpha)\partial_{s}[F_1(P_{3},s)]e^{isH_{a}}e^{-isH}\psi\nonumber\\
&+(-i)\int_1^tds F_{c,3,\alpha}(x_3,s,P_3)e^{isH_{a}}V_{13}(x_1-x_3)e^{-isH}\psi\nonumber\\
&+(-i)\int_1^tdtF_{c,3,\alpha}(x_3,s,P_3)e^{isH_{a}}V_{23}(x_2-x_3)e^{-isH}\psi\nonumber\\
=:& \Omega_{a,\alpha}^{*}(1)\psi+\psi_{1}(t)+\psi_{2}(t)+\psi_{3}(t)+\psi_4(t).
\end{align}
The unitarity of $e^{-iH}$ and $e^{iH_{a}}$ implies that $\Omega_{a,\alpha}^{*}(1)\psi\in \s^2_x(\mathbb{R}^9)$. For $\psi_{3}(t)$ and $\psi_{4}(t)$, we use Assumptions \ref{asp: omega} and \ref{asp: short} and the unitarity of $e^{-isH}$ on $\s^2_x(\R^9)$ to obtain
\begin{align}
&\int_1^\infty \| F_{c,3,\alpha}(x_3,s,P_3)e^{isH_{a}}V_{13}(x_1-x_3)e^{-isH} \psi\|_{\s_x^2(\mathbb{R}^9)}ds\nonumber\\
\leq& \int_1^\infty\|  F_{c,3,\alpha}(x_3,s,P_3)e^{isH_{a}}\langle x_1-x_3\rangle^{-1-\epsilon}   \|_{\s^2_{x}(\mathbb{R}^9)\to\s^2_{x}(\mathbb{R}^9) } \| \langle \eta\rangle^{1+\epsilon}V_{13}(\eta)\|_{\s^\infty_\eta(\R^3)}\| e^{-isH}\psi\|_{\s^2_x(\R^9)}ds\nonumber\\
\leq&\left(\int_1^\infty\|  F_{c,3,\alpha}(x_3,s,P_3)e^{isH_{a}}\langle x_1-x_3\rangle^{-1-\epsilon}   \|_{\s^2_{x}(\mathbb{R}^9)\to\s^2_{x}(\mathbb{R}^9) }ds\right)\times\| \langle \eta\rangle^{1+\epsilon}V_{13}(\eta)\|_{\s^\infty_\eta(\R^3)}\| \psi\|_{\s^2_x(\R^9)}\nonumber\\
<&\infty,\label{Feb.7.1}
\end{align}
and
\begin{align}
&\int_1^\infty \| F_{c,3,\alpha}(x_3,s,P_3)e^{isH_{a}}V_{23}(x_2-x_3)e^{-isH} \psi\|_{\s_x^2(\mathbb{R}^9)}ds\nonumber\\
\leq& \int_1^\infty\|  F_{c,3,\alpha}(x_3,s,P_3)e^{isH_{a}}\langle x_2-x_3\rangle^{-1-\epsilon}   \|_{\s^2_{x}(\mathbb{R}^9)\to\s^2_{x}(\mathbb{R}^9) } \| \langle \eta\rangle^{1+\epsilon}V_{23}(\eta)\|_{\s^\infty_\eta(\R^3)}\| e^{-isH}\psi\|_{\s^2_x(\R^9)}ds\nonumber\\
\leq&\left(\int_1^\infty\|  F_{c,3,\alpha}(x_3,s,P_3)e^{isH_{a}}\langle x_2-x_3\rangle^{-1-\epsilon}   \|_{\s^2_{x}(\mathbb{R}^9)\to\s^2_{x}(\mathbb{R}^9) }ds\right)\times\| \langle \eta\rangle^{1+\epsilon}V_{23}(\eta)\|_{\s^\infty_\eta(\R^3)}\| \psi\|_{\s^2_x(\R^9)}\nonumber\\
<&\infty.\label{Feb.7.2}
\end{align}
Therefore, $\psi_{3}(\infty)$ and $\psi_{4}(\infty)$ exist in $\s_x^2(\mathbb{R}^9)$. As for $\psi_1(t)$, we utilize the propagation estimates first introduced by \cite{SW1} and \cite{SW2}: For a family of operators $\{B(t)\}_{t\in \R}$ on $\s^2_x(\R^9)$ and a family of data $\{\phi(t)\}_{t\in \R}$ in $\s^2_x(\R^9)$, we define
\eq
\langle B(t): \phi(t)\rangle_t:= (\phi(t), B(t)\phi(t))_{\s^2_x(\R^9)}.
\eeq
We observe $B(t)=e^{-itH_{a}}\tilde{F}_1(P_3, t)F_c(|x_3|\leq t^\alpha)\tilde{F}_1(P_3,t)e^{itH_{a}}$ with respect to $\psi(t)$. We compute $\langle B(t): \psi(t)\rangle_t$ as follows:
\begin{align}
&\langle B(t): \psi(t)\rangle_t=\langle B(t): \psi(t)\rangle_t\vert_{t=1}\nonumber\\
&+(-i)\int_1^tds (\psi(s), e^{-isH_{a}}\tilde{F}_1(P_3, t)F_c(|x_3|\leq s^\alpha)\tilde{F}_1(P_3,s)e^{isH_{a}}(V_{13}+V_{23})\psi(s) )_{\s^2_x(\mathbb{R}^9)}\nonumber\\
&+i\int_1^tds ((V_{13}+V_{23})\psi(s), e^{-isH_{a}}\tilde{F}_1(P_3, t)F_c(|x_3|\leq s^\alpha)\tilde{F}_1(P_3,s)e^{isH_{a}}\psi(s) )_{\s^2_x(\mathbb{R}^9)}\nonumber\\
&+\int_1^tds (\psi(s),   e^{-isH_{a}}\partial_s[\tilde{F}_1(P_3, s)F_c(|x_3|\leq s^\alpha)\tilde{F}_1(P_3, s)]e^{isH_{a}}\psi(s) )_{\s^2_x(\mathbb{R}^9)}\nonumber\\
=:& a_0+\int_1^tdsa_1(s)+\int_1^tdsa_2(s)+\int_1^tdsa_3(s).
\end{align}
Using the unitarity of $e^{\pm iH}$ on $\s^2_x(\R^9)$, we obtain
\begin{align}
|a_0|\leq &\|\psi(1)\|_{\s^2_x(\R^9)}\|e^{-iH_{a}}\tilde{F}_1(P_3, 1)F_c(|x_3|\leq 1)\tilde{F}_1(P_3,1)e^{iH_{a}} \|_{\s^2_x(\R^9)\to \s^2_x}\|\psi(1)\|_{\s^2_x(\R^9)}\nonumber\\
\leq & \|\psi(0)\|_{\s^2_x(\R^9}.
\end{align}
In addition, \eqref{Feb.7.1} and \eqref{Feb.7.2} imply that $a_1(s), a_2(s)\in \s^1_s[1,\infty)$. Regarding $a_3(s)$, we have
\begin{align}
a_3(s)=&(\psi(s),   e^{-isH_{a}}\tilde{F}_1(P_3, s)\partial_s[F_c(|x_3|\leq s^\alpha)]\tilde{F}_1(P_3, s)e^{isH_{a}}\psi(s) )_{\s^2_x(\mathbb{R}^9)}\nonumber\\
&+(\psi(s),   e^{-isH_{a}}\partial_s[\tilde{F}_1(P_3,s)]F_c(|x_3|\leq s^\alpha)\tilde{F}_1(P_3,s)e^{isH_{a}}\psi(s) )_{\s^2_x(\mathbb{R}^9)}\nonumber\\
&+(\psi(s),   e^{-isH_{a}}\tilde{F}_1(P_3,s)F_c(|x_3|\leq s^\alpha)\partial_s[\tilde{F}_1(P_3,s)]e^{isH_{a}}\psi(s) )_{\s^2_x(\mathbb{R}^9)}\nonumber\\
=:&a_{3,1}(s)+a_{3,2}(s)+a_{3,3}(s).
\end{align}
Here, we have $a_{3,1}(s)\geq 0$, 
\eq
a_{3,2}(s)=(\psi(s),   e^{-isH_{a}}\sqrt{F_c}\partial_s[\tilde{F}_1(P_3, s)]\tilde{F}_1(P_3, s)\sqrt{F_c}e^{isH_{a}}\psi(s) )_{\s^2_x(\mathbb{R}^9)}+a_{3,2,r}(s)\label{identity: a32}
\eeq
and 
\eq
a_{3,3}(s)=(\psi(s),   e^{-isH_{a}}\sqrt{F_c}\partial_s[\tilde{F}_1(P_3, s)]\tilde{F}_1(P_3, s)\sqrt{F_c}e^{isH_{a}}\psi(s) )_{\s^2_x(\mathbb{R}^9)}+a_{3,3,r}(s)
\eeq
for some reminder terms $a_{3,2,r}(s),a_{3,3,r}(s)$. By using Corollary 2.5 in \cite{SW2}, we prove that $$a_{3,2,r}(s),a_{3,3,r}(s)\in \s^1_s[1,\infty) :$$ 

When there is no confusion in the context, we use $F_c$ and $\tilde{F}_1$ to refer to $F_c(|x_3|\leq s^\alpha)$ and $\tilde{F}_1(P_3,s)$, respectively. By employing the identity
\begin{align}
&\partial_s[\tilde{F}_1]F_c\tilde{F}_1\nonumber\\
=&[\partial_s[\tilde{F}_1], \sqrt{F_c}] \sqrt{F_c}\tilde{F}_1+\sqrt{F_c}\partial_s[\tilde{F}_1]\sqrt{F_c}\tilde{F}_1\nonumber\\
=&[\partial_s[\tilde{F}_1], \sqrt{F_c}] \sqrt{F_c}\tilde{F}_1+\sqrt{F_c}\partial_s[\tilde{F}_1][\sqrt{F_c}, \tilde{F}_1]+\sqrt{F_c}\partial_s[\tilde{F}_1]F_1\sqrt{F_c},
\end{align}
we can show that \eqref{identity: a32} is valid by taking
\begin{align}
    a_{3,2,r}(s)=&(\psi(s),e^{-isH_a} [\partial_s[\tilde{F}_1], \sqrt{F_c}]\sqrt{F_c}\tilde{F}_1 e^{isH_a}\psi(s) )_{\s^2_x(\R^9)}\nonumber\\
    &+(\psi(s),e^{-isH_a} \sqrt{F_c}\partial_s[\tilde{F}_1][\sqrt{F_c}, \tilde{F}_1] e^{isH_a}\psi(s) )_{\s^2_x(\R^9)}.
\end{align}
 $a_{3,2,r}(s), a_{3,3,r}(s)\in \s^1_s[1,\infty)$ by using Corollary 2.5 in \cite{SW2}.

Furthermore, we have for all $q_3\in \R^3$ that
\eq
\partial_s[F_1(|q_3-\eta_{3k}|>\frac{1}{s^{\alpha}})]\geq 0, \quad k=1,\cdots, N_3,
\eeq
which implies
\eq
(\psi(s),   e^{-isH_{a}}\sqrt{F_c}\partial_s[\tilde{F}_1(P_3, s)]\tilde{F}_1(P_3, s)\sqrt{F_c}e^{isH_{a}}\psi(s) )_{\s^2_x(\mathbb{R}^9)}\geq 0,
\eeq
and 
\eq
(\psi(s),   e^{-isH_{a}}\sqrt{F_c}\partial_s[\tilde{F}_1(P_3, s)]\tilde{F}_1(P_3, s)\sqrt{F_c}e^{isH_{a}}\psi(s) )_{\s^2_x(\mathbb{R}^9)}\geq 0.
\eeq
Thus, via propagation estimates (see \cite{SW1}), we conclude that $a_{3,1}(s)\in \s^1_s[1,\infty)$, which implies
\eq
\psi_1(\infty)=\int_1^\infty ds \partial_{s}[ F_c(\frac{|x_3|}{s^\alpha}\leq 1)]\tilde{F}_1(P_{3}, s)e^{isH_{a}}e^{-isH}\psi \text{ exists in }\s^2_x(\mathbb{R}^9).\label{psiinfty1}
\eeq
Similarly, by observing $B(t)=e^{-itH_{a}}F_c(|x_2|/t^\alpha\leq 1)\tilde{F}_1(P_3, t)F_c(|x_2|/t^\alpha\leq 1)e^{itH_{a}}$ with respect to $\psi(t)$, we obtain $b(s)\in \s^1_s[1,\infty)$, where 
\eq
b(s):=(\psi(s),   e^{-isH_{a}}F_1(|x_3|\leq s^\alpha)\partial_s[\tilde{F}_1(P_3, s)]F_1(|x_3|\leq s^\alpha)e^{isH_{a}}\psi(s) )_{\s^2_x(\mathbb{R}^9)}.
\eeq
The fact that $b(s)\in \s^1_s[1,\infty)$ implies 
\eq
\int_1^\infty ds \partial_{s}[\tilde{F}_1(P_{3}, s)]  F_c(\frac{|x_3|}{s^\alpha}\leq 1)e^{isH_{a}}e^{-isH}\psi \text{ exists in }\s^2_x(\mathbb{R}^9).\label{Feb8.2}
\eeq
Furthermore, by applying Corollary 2.5 in \cite{SW2}, we obtain
\eq
\|[\partial_{s}[\tilde{F}_1(P_{3}, s)],  F_c(\frac{|x_3|}{s^\alpha}\leq 1)]e^{isH_{a}}e^{-isH}\psi \|_{\s^2_x(\mathbb{R}^9)}\in \s^1_s[1,\infty).\label{Feb8.3}
\eeq
Combining \eqref{Feb8.3} and \eqref{Feb8.2}, we deduce that
\eq
\psi_2(\infty)=\int_1^\infty ds F_c(\frac{|x_3|}{s^\alpha}\leq 1)\partial_{s}[\tilde{F}_1(P_{3},s)]e^{isH_{12}}e^{-isH}\psi \text{ exists in }\s^2_x(\mathbb{R}^9).\label{psiinfty2}
\eeq 
Using \eqref{psiinfty1} and \eqref{psiinfty2}, we conclude the existence of $\Omega_{a,\alpha}^{*}\psi$ in $\s^2_x(\mathbb{R}^9)$ which implies the existence of $\Omega_{a,\alpha}^{*}$ on $\s^2_x(\mathbb{R}^9)$ when $a=(12)(3)$. Analogously, we can obtain the existence of operators $\Omega_{a^2,\alpha}^{*}$ and $\Omega_{a^3,\alpha}^{*}$ where $a^2\equiv (13)(2)$ and $a^3\equiv (23)(1)$. To analyze the free channel wave operator, we use a similar argument. Let
\eq
\Omega_{free,\alpha}^{*}(t):=e^{itH_{0}}J_{\alpha,free}(t)e^{-itH}=\left(\Pi_{l=1}^{l=3} F_{c,l,\alpha}(x_l,t,P_l)\right)e^{itH_{0}}e^{-itH}.
\eeq
Expanding $\Omega_{free,\alpha}^{*}(t)\psi$ using Cook's method, we obtain
\begin{align}
\Omega_{free,\alpha}^{*}(t)\psi=&\Omega_{free,\alpha}^{*}(1)\psi+\int_1^tds \partial_{s}[ \left(\Pi_{l=1}^{l=3} F_{c,l,\alpha}(x_l,t,P_l)\right)]e^{isH_{0}}e^{-isH}\psi\nonumber\\
&+(-i)\int_1^tds \left(\Pi_{l=1}^{l=3} F_{c,l,\alpha}(x_l,t,P_l)\right)e^{isH_{0}}V(x)e^{-isH}\psi\nonumber\\
=:& \Omega_{free,\alpha}^{*}(1)\psi+\psi_{1}(t)+\psi_{2}(t).
\end{align}
Due to the unitarity of $e^{-iH}$ and $e^{iH_{0}}$, we have $ \Omega_{free,\alpha}^{*}(1)\psi\in \s^2_x(\mathbb{R}^9)$. For $\psi_{2}(t)$, by utilizing Assumptions \ref{asp: short} and \ref{asp: omega} along with the unitarity of $e^{-isH}$ on $\s^2_x(\R^9)$, we derive 
\begin{align}
&\|\left(\Pi_{l=1}^{l=3} F_{c,l,\alpha}(x_l,t,P_l)\right)e^{isH_{0}}V(x)e^{-isH} \psi\|_{\s_x^2(\mathbb{R}^9)}\nonumber\\
\leq&\sum\limits_{a=(jk)(l)\in L}\|   F_{c,j,\alpha}(x_j,t,P_j)e^{is\omega_j(P_j)}V_{jk}(x_j-x_k)e^{-isH} \psi   \|_{\s^2_{x}(\mathbb{R}^9)}\nonumber\\
\leq& \sum\limits_{a=(jk)(l)\in L} \| F_{c,j,\alpha}(x_j,t,P_j)e^{is\omega_j(P_j)}\langle x_j-x_k\rangle^{-4}\|_{\s^2_x(\R^9)\to \s^2_x(\R^9)}\|\langle \eta\rangle^{4}V_{jk}(\eta)\|_{\s^\infty_\eta(\R^3)}\|\psi\|_{\s^2_x(\R^9)}\nonumber\\
\in & \s^1_s[1,\infty).\label{Feb.7.1f}
\end{align}
As a result, $\psi_{2}(\infty)$ exists in $\s_x^2(\mathbb{R}^9)$. For $\psi_1(t)$, we use the propagation estimates, initiated by \cite{SW1} and \cite{SW2}, and estimate \eqref{Feb.7.1f}. Similarly, we obtain the existence of $\psi_1(\infty)$ in $\s^2_x(\mathbb{R}^9)$. This concludes the proof.

\end{proof}
\begin{lemma}\label{Lem: weak}For all $a=(jk)(l)\in L$, we have
\eq
w\text{-}\lim\limits_{t\to \infty}(1-F_{c,l,\alpha}(x_l,t,P_l))e^{itH_a}e^{-itH}=0\quad \text{ on }\s^2_x(\R^9)\label{weak: eq1}
\eeq
and
\eq
w\text{-}\lim\limits_{t\to \infty}\left[1-\left(\Pi_{l=1}^{l=3} F_{c,l,\alpha}(x_l,t,P_l)\right)\right]e^{itH_0}e^{-itH}=0\quad \text{ on }\s^2_x(\R^9).\label{weak: eq2}
\eeq
Therefore, if Assumptions \ref{asp: critical}, \eqref{asp: short}, and \ref{asp: omega} hold, and if $\alpha$ and $\alpha'$ satisfy Assumption \ref{asp: critical}, we have
\eq
\Omega_{free,\alpha}^{*}=\Omega_{free,\alpha'}^{*},\quad \text{ on }\s^2_x(\R^9),
\eeq
and
\eq
\Omega_{a,\alpha}^{*}=\Omega_{a,\alpha'}^{*},\quad \text{ on }\s^2_x(\R^9).
\eeq
\end{lemma}
\begin{proof}Choose $\psi,\phi\in \s^2_x(\R^9)$. Since
\begin{align*}
&1-F_{c,l,\alpha}(x_l,t,P_l)\\
=& F_c(|x_l|>t^\alpha)+F_c(|x_l|\leq t^\alpha )\left(1-\left(\Pi_{k=1}^{N_j}F_1(|P_j-\eta_{jk}|>\frac{1}{t^{\alpha/2}})\right)F_1(|P_j|\leq t^{\alpha/2})\right),
\end{align*}
using the unitarity of $e^{itH_a}$ and $e^{-itH}$ on $\s^2_x(\R^9)$ and using H\"older's inequality, we have, as $t\to \infty$,
\begin{align*}
   & \left|(\phi(x), (1-F_{c,l,\alpha}(x_l,t,P_l))e^{itH_a}e^{-itH}\psi(x))_{\s^2_x(\R^9)} \right|\\
   \leq & \| (1-F_{c,l,\alpha}(x_l,t,P_l))\phi(x) \|_{\s^2_x(\R^9)}\|\psi(x)\|_{\s^2_x(\R^9)}\\
   \to& 0,
\end{align*}
and 
\begin{align*}
   & \left|(\phi(x), \left[1-\left(\Pi_{l=1}^{l=3} F_{c,l,\alpha}(x_l,t,P_l)\right)\right]e^{itH_0}e^{-itH}\psi(x))_{\s^2_x(\R^9)} \right|\\
   \leq & \|\left[1-\left(\Pi_{l=1}^{l=3} F_{c,l,\alpha}(x_l,t,P_l)\right)\right]\phi(x) \|_{\s^2_x(\R^9)}\|\psi(x)\|_{\s^2_x(\R^9)}\\
   \to& 0.
\end{align*}
We get \eqref{weak: eq1} and \eqref{weak: eq2}. We finish the proof.
    
\end{proof}
The wave operators in the two-body channel, in the context of our three-body system, are given by
\eq
\Omega_{\alpha}^{a,*}:=s\text{-}\lim\limits_{t\to \infty} F_{c,j,\alpha}(x_j,t,P_j)e^{itH_0} e^{-itH_{a}}\quad \text{ on }\s^2_x(\mathbb{R}^9),
\eeq
for $a=(jk)(l)\in L$.
\begin{proposition}\label{thm:cW1}If Assumptions \ref{asp: short}, \ref{asp: critical} and \ref{asp: omega} are satisfied, then  $\Omega_{\alpha,jk}^{*}$ exists for all $a=(jk)(l)\in L,$ where $\alpha$ satisfies Assumption \ref{asp: omega}.
\end{proposition}
\begin{proof}
Choose $\psi\in \s_x^2(\mathbb{R}^9)$. We begin with $a=(12)(3)$ and the proof of the existence of $\Omega_{\alpha}^{a,*}$. Let
\eq
\Omega_{\alpha}^{a,*}(t)=F_{c,j,\alpha}(x_j,t,P_j)e^{itH_0} e^{-itH_{a}}.
\eeq
Using Cook's method to expand $\Omega_{\alpha}^{a,*}(t)\psi$, we have
\begin{align}
\Omega_{\alpha}^{a,*}(t)\psi=&\Omega_{\alpha}^{a,*}(1)\psi+\int_1^tds \partial_{s}[ F_{c,1,\alpha}(x_1,t,P_1)]e^{isH_{0}}e^{-isH_{a}}\psi\nonumber\\
&+(-i)\int_1^tds F_{c,1,\alpha}(x_1,t,P_1)e^{isH_{0}}V_{12}(x_1-x_2)e^{-isH_{a}}\psi\nonumber\\
=:& \Omega_{\alpha}^{a,*}(1)\psi+\psi_{1}(t)+\psi_{2}(t).
\end{align}
Thanks to the unitarity of $e^{-iH_0}$ and $e^{iH_a}$, we have $ \Omega_{\alpha}^{a,*}(1)\psi\in \s^2_x(\mathbb{R}^9)$. Regarding $\psi_{2}(t)$, using Assumption \ref{asp: omega} and the unitarity of $e^{-isH}$ on $\s^2_x(\R^9)$, we have
\begin{align}
&\| F_{c,1,\alpha}(x_1,t,P_1)e^{isH_{0}}V_{12}(x_1-x_2)e^{-isH_{a}}\psi\|_{\s_x^2(\mathbb{R}^9)}\nonumber\\
\leq &\| F_{c,1,\alpha}(x_1,t,P_1)e^{is\omega_1(P_1)}\langle x_1-x_2\rangle^{-4} \|_{\s^2_x(\R^9)\to\s^2_x(\R^9) }\|\langle \eta\rangle^{4}V_{12}(\eta)\|_{\s^\infty_\eta(\R^3)}\| \psi\|_{\s^2_x(\R^9)}\nonumber\\
\in & \s^1_s[1,\infty),\label{Feb8.32}
\end{align}
which implies that $\psi_2(\infty)$ exists in $\s^2_x(\mathbb{R}^9)$. Similarly, for $\psi_1(t)$, we can use the propagation estimate and estimate \ref{Feb8.32} to show that $\psi_1(\infty)$ exists in $\s^2_x(\mathbb{R}^9)$. Thus, $\Omega_{\alpha}^{a,*}\psi$ exists in $\s^2_x(\R^9)$. We can therefore conclude that $\Omega_{\alpha}^{a,*}$ exists on $\s^2_x(\mathbb{R}^9)$ when $a=(12)(3)$. Similarly, we can establish the existence of the operators $\Omega_{\alpha}^{a^2,*}$ and $\Omega_{\alpha}^{a^3,*}$, where $a^2=(13)(2)$ and $a^3=(1)(23)$. This completes the proof.
\end{proof}
\begin{lemma}\label{Lem: weak: 2}For all $a=(jk)(l)\in L$, we have
\eq
w\text{-}\lim\limits_{t\to \infty}(1-F_{c,j,\alpha}(x_l,t,P_l))e^{itH_0}e^{-itH_a}=0\quad \text{ on }\s^2_x(\R^9).\label{2: weak: eq1}
\eeq
Therefore, if Assumptions \ref{asp: critical}, \eqref{asp: short}, and \ref{asp: omega} hold, and if $\alpha$ and $\alpha'$ satisfy Assumption \ref{asp: critical}, we have
\eq
\Omega_{\alpha}^{a,*}=\Omega_{\alpha'}^{a,*},\quad \text{ on }\s^2_x(\R^9).
\eeq
\end{lemma}
\begin{proof}
    We can arrive at the same conclusion by using a similar argument to the one we employed in the proof of Lemma \ref{Lem: weak}.
\end{proof}

\subsection{Channel projections}\label{sec: channel projection}
Based on Theorem \ref{thm: channel}, we can define the channel projection in a new way. Recall that for $a=(jk)(l)\in L$, we have
\eq
J_{\alpha,a}(t)=e^{-itH_a} F_{c,l,\alpha}(x_l,t,P_l)e^{itH_a}.
\eeq
When $a=(jk)(l)\in L$, the $H_a$ channel projection is defined as
\eq
 P_{a}:=s\text{-}\lim\limits_{t\to \infty} e^{itH}e^{-itH_a}J_{\alpha,a}(t)e^{itH_a}e^{-itH}\quad \text{ on }\s^2_x(\mathbb{R}^9).
\eeq
When $a=(1)(2)(3)$, the $H_a$ channel is free channel and the free channel projection is defined as 
\eq
P_{free}:=s\text{-}\lim\limits_{t\to \infty} e^{itH}e^{-itH_0}J_{\alpha,free}(t)e^{itH_0}e^{-itH}\quad \text{ on }\s^2_x(\mathbb{R}^9),
\eeq
where we would like to remind the reader that 
\eq
J_{\alpha, free}(t)=e^{-itH_0} \left(\Pi_{l=1}^{l=3}F_{c,l,\alpha}(x_l,t,P_l)\right)e^{itH_0}.
\eeq 
We omit the subscript $\alpha$ here because, as Lemma \ref{Lem: Pab: weak} shows, $P_a$ and $P_{free}$ do not depend on the choice of $\alpha.$

The free channel projection of $H_a (a=(jk)(l)\in L)$ is defined by 
\eq
P_{sc}(H_a):=s\text{-}\lim\limits_{t\to \infty} e^{itH_a}e^{-itH_0}F_{c,j,\alpha}(x_j,t,P_j)e^{itH_0} e^{-itH_a},\quad \text{ on }\s^2_x(\R^9)
\eeq
and we let
\eq
P_{bs}(H_a):=1-P_{sc}(H_a).
\eeq
\begin{lemma}\label{Lem: Pab: weak}If $\alpha$ and $\alpha'$ satisfy Assumption \ref{asp: critical} and if Assumptions \ref{asp: critical}, \ref{asp: short} and \ref{asp: omega} hold, then we have
\eq
s\text{-}\lim\limits_{t\to \infty}e^{itH}e^{-itH_a}(J_{\alpha,a}(t)-J_{\alpha',a}(t))e^{itH_a}e^{-itH}=0,\quad \text{ on }\s^2_x(\R^9)\label{Lem: Pab: weak: eq1}
\eeq
and
\eq
s\text{-}\lim\limits_{t\to \infty}e^{itH}e^{-itH_0}(J_{\alpha,free}(t)-J_{\alpha',free}(t))e^{itH_0}e^{-itH}=0,\quad \text{ on }\s^2_x(\R^9).\label{Lem: Pab: weak: eq2}
\eeq
\end{lemma}
\begin{proof}Using Theorem \ref{thm: channel} and Lemma \ref{Lem: weak}, we have 
\eq
s\text{-}\lim\limits_{t\to \infty}(J_{\alpha,a}(t)-J_{\alpha',a}(t))e^{itH_a}e^{-itH}=0,\quad \text{ on }\s^2_x(\R^9)
\eeq
and
\eq
s\text{-}\lim\limits_{t\to \infty}(J_{\alpha,free}(t)-J_{\alpha',free}(t))e^{itH_0}e^{-itH}=0,\quad \text{ on }\s^2_x(\R^9).
\eeq
Using the unitarity of $e^{itH}$, $e^{-itH_a}$ and $e^{-itH_0}$ on $\s^2_x(\R^9)$, we get \eqref{Lem: Pab: weak: eq1} and \eqref{Lem: Pab: weak: eq2}. We finish the proof. 
    
\end{proof}
\begin{proposition}\label{cor 1}If Assumptions \ref{asp: short}, \ref{asp: critical} and \ref{asp: omega} are satisfied, then $P_a$, $P_{free},$ $P_{sc}(H_a)$ and $P_{bs}(H_a)$ exist for all $a=(jk)(l)\in L$.
\end{proposition}
\begin{proof}If we can demonstrate that for all $a=(jk)(l)\in L$, the limits
\eq
\Omega_{a}:=s\text{-}\lim\limits_{t\to \infty}e^{itH}e^{-itH_a},\quad \text{ on }\s^2_x(\R^9)\label{cor 1: 1},
\eeq
\eq
\Omega_{free}:=s\text{-}\lim\limits_{t\to \infty}e^{itH}e^{-itH_0},\quad \text{ on }\s^2_x(\R^9)\label{cor 1: 2}
\eeq
and 
\eq
\Omega^{a}:=s\text{-}\lim\limits_{t\to \infty}e^{itH_a}e^{-itH_0},\quad \text{ on }\s^2_x(\R^9)\label{cor 1: 3}
\eeq
exist, then we get the existence of $P_a,$ $P_{free}$ and $P_{sc}(H_a)$ by using Theorem \ref{thm: channel} and Proposition \ref{thm:cW1} and furthermore, the existence of $P_{bs}(H_a)$ can be deduced from the existence of $P_{sc}(H_a)$. Now let us prove the existence of \eqref{cor 1: 1}, \eqref{cor 1: 2} and \eqref{cor 1: 3}: Choose $f(x)\in \s^2_x(\R^9)$ and take $a=(12)(3)$. Let $\alpha$ be as in Assumption \ref{asp: omega}. Break $e^{itH}e^{-itH_a}f$ into $2$ pieces:
\begin{align}
    e^{itH}e^{-itH_a}f=&e^{itH}e^{-itH_a}F_{c,3,\alpha}(x_3,t,P_3)f+e^{itH}e^{-itH_a}\left(1-F_{c,3,\alpha}(x_3,t,P_3)\right)f\nonumber\\
    =:&f_1(t)+f_2(t).
\end{align}
For $f_1(t)$, we use Cook's method:
\begin{align}
    &e^{itH}e^{-itH_a}F_{c,3,\alpha}(x_3,t,P_3)f\nonumber\\
    =&e^{iH}e^{-iH_a}F_{c,3,\alpha}(x_3,1,P_3)f+i\int_1^t e^{isH}V_{13}(x_1-x_3)e^{-isH_a}F_{c,3,\alpha}(x_3,s,P_3)f(x)ds\nonumber\\
    &+i\int_1^t e^{isH}V_{23}(x_2-x_3)e^{-isH_a}F_{c,3,\alpha}(x_3,s,P_3)f(x)ds+\int_1^t e^{isH}e^{-isH_a} \p_s[F_{c,3,\alpha}(x_3,s,P_3)]f(x)\nonumber\\
    =:&f_1(1)+\int_1^tdsf_{1,in1}(s)+\int_1^tdsf_{1,in2}(s)+\int_1^tdsf_{1,+}(s).
\end{align}
$f_1(1)\in \s^2_x(\R^9)$. Assumption \ref{asp: omega} allows us to demonstrate that both $\int_1^\infty ds f_{1,in1}(s)$ and $\int_1^\infty ds f_{1,in2}(s)$ exist in $\s^2_x(\R^9)$, by utilizing the unitarity of $e^{-is(H_a-\omega_3(P_3))}$ and $e^{isH}$ on $\s^2_x(\R^9)$, as well as Assumption \ref{asp: omega}:
\begin{align}
   &\|f_{1,in1}(s)\|_{\s^2_x(\R^9)}\nonumber\\
   =&\| \|e^{isH}V_{13}(x_1-x_3)e^{-is\omega_3(P_3)}F_{c,3,\alpha}(x_3,s,P_3)e^{-is(H_a-\omega_3(P_3))}f(x)\|_{\s^2_x(\R^9)}\nonumber\\
   \leq &\|e^{isH}V_{13}(x_1-x_3)e^{-is\omega_3(P_3)}F_{c,3,\alpha}(x_3,s,P_3) \|_{\s^2_x(\R^9)\to\s^2_x(\R^9) }\|e^{-is(H_a-\omega_3(P_3))}f(x)\|_{\s^2_x(\R^9)}\nonumber\\
   \leq & \| \langle \eta\rangle^{4}V_{13}(\eta)\|_{\s^\infty_\eta(\R^3)}\|\langle x_1-x_3\rangle^{-4}e^{-is\omega_3(P_3)}F_{c,3,\alpha}(x_3,s,P_3)\|_{\s^2_x(\R^9)\to \s^2_x(\R^9)}\|f(x)\|_{\s^2_x(\R^9)}\nonumber\\
   \in& \s^1_s[1,\infty),\label{example  in}
   \end{align}
and
\begin{align}
   &\|f_{1,in2}(s)\|_{\s^2_x(\R^9)}\nonumber\\
   =&\| \|e^{isH}V_{23}(x_2-x_3)e^{-is\omega_3(P_3)}F_{c,3,\alpha}(x_3,s,P_3)e^{-is(H_a-\omega_3(P_3))}f(x)\|_{\s^2_x(\R^9)}\nonumber\\
   \leq &\|e^{isH}V_{23}(x_2-x_3)e^{-is\omega_3(P_3)}F_{c,3,\alpha}(x_3,s,P_3) \|_{\s^2_x(\R^9)\to\s^2_x(\R^9) }\|e^{-is(H_a-\omega_3(P_3))}f(x)\|_{\s^2_x(\R^9)}\nonumber\\
   \leq & \| \langle \eta\rangle^{4}V_{23}(\eta)\|_{\s^\infty_\eta(\R^3)}\|\langle x_2-x_3\rangle^{-4}e^{-is\omega_3(P_3)}F_{c,3,\alpha}(x_3,s,P_3)\|_{\s^2_x(\R^9)\to \s^2_x(\R^9)}\|f(x)\|_{\s^2_x(\R^9)}\nonumber\\
   \in& \s^1_s[1,\infty).
\end{align}
We also have that $\int_1^\infty ds f_{1,+}(s)$ exists in $\s^2_x(\R^9)$: For $t>T>1$, we can use Fubini's Theorem, the unitarity of $e^{isH}e^{-isH_a}$ on $\s^2_x(\R^9)$, and H\"older's inequality (in this order) to obtain the following:
\begin{align}
    &\|\int_T^t e^{isH}e^{-isH_a}\p_s[F_{c,3,\alpha}(x_3,s,P_3)]f(x)\|_{\s^2_x(\R^9)}\nonumber\\
=&\sup\limits_{\|g(x)\|_{\s^2_x(\R^9)}=1} \left| (g(x), \int_T^t e^{isH}e^{-isH_a}\p_s[F_{c,3,\alpha}(x_3,s,P_3)]f(x)ds)_{\s^2_x(\R^9)}\right|\nonumber\\
=&\sup\limits_{\|g(x)\|_{\s^2_x(\R^9)}=1} \left|\int_T^t (g(x),  e^{isH}e^{-isH_a}\p_s[F_{c,3,\alpha}(x_3,s,P_3)]f(x))_{\s^2_x(\R^9)}ds\right|\nonumber\\
=&\sup\limits_{\|g(x)\|_{\s^2_x(\R^9)}=1} \left|\int_T^t (e^{isH_a}e^{-isH}g(x),  \p_s[F_{c,3,\alpha}(x_3,s,P_3)]f(x))_{\s^2_x(\R^9)}ds\right|\nonumber\\
=&\sup\limits_{\|g(x)\|_{\s^2_x(\R^9)}=1} \left|\int_{\R^9} dx \int_T^t (e^{isH_a}e^{-isH}g(x))^*  \p_s[F_{c,3,\alpha}(x_3,s,P_3)]f(x)ds\right|\nonumber\\
\leq & \sup\limits_{\|g(x)\|_{\s^2_x(\R^9)}=1} \left|\int_{\R^9} dx \int_T^t (e^{isH_a}e^{-isH}g(x))^*  \p_s[F_{c,3,\alpha}(x_3,s,P_3)]f(x)ds\right|.\label{May.eq1}
\end{align}
We apply the Propagation estimates (as detailed in \cite{SW1}) to yield:
\begin{align}
  & \int_{\R^9} dx  \int_T^t \p_s[F_{c,3,\alpha}(x_3,s,P_3)]|e^{isH_a}e^{-isH}g(x)|^2ds \nonumber\\
  \leq & 2\sup\limits_{t\geq 1}(e^{itH_a}e^{-itH}g(x),F_{c,3,\alpha}(x_3,t,P_3)e^{itH_a}e^{-itH}g(x) )_{\s^2_x(\R^9)}\nonumber\\
  &+2E\|g(x)\|_{\s^2_x(\R^9)}^2\times \sup\limits_{1\leq j<k\leq 3}\|\langle \eta\rangle^{4}V_{jk}(\eta)\|_{\s^\infty_\eta(\R^3)} \nonumber\\
  \lesssim_E &\left(1+\sup\limits_{1\leq j<k\leq 3}\|\langle \eta\rangle^{4}V_{jk}(\eta)\|_{\s^\infty_\eta(\R^3)}\right)\|g(x)\|_{\s^2_x(\R^9)}^2.\label{g eq}
\end{align}
Based on \eqref{May.eq1} and \eqref{g eq}, we have
\begin{align}
&\|\int_T^t e^{isH}e^{-isH_a}\p_s[F_{c,3,\alpha}(x_3,s,P_3)]f(x)\|_{\s^2_x(\R^9)}\nonumber\\
    \leq & \sup\limits_{\|g(x)\|_{\s^2_x(\R^9)}=1} \int_{\R^9} dx \left( \int_T^t \p_s[F_{c,3,\alpha}(x_3,s,P_3)]|e^{isH_a}e^{-isH}g(x)|^2ds\right)^{1/2}\times\nonumber\\
    &\left( \int_T^t \p_s[F_{c,3,\alpha}(x_3,s,P_3)]|f(x)|^2ds\right)^{1/2}\nonumber\\
    \lesssim_{E} & \left(1+\sup\limits_{1\leq j<k\leq 3}\|\langle \eta\rangle^{1+\epsilon}V_{jk}(\eta)\|_{\s^\infty_\eta(\R^3)}\right)\times\left( \int_{\R^9} \int_T^t \p_s[F_{c,3,\alpha}(x_3,s,P_3)]|f(x)|^2dsdx \right)^{1/2}\nonumber\\
    \to & 0,\label{example  +}
\end{align}
as $T\to \infty$. As a result, we can infer that $\int_1^\infty f_{1,+}(s)ds$ exists in $\s^2_x(\R^9)$. Hence, $f_1(\infty)$ exists in $\s^2_x(\R^9)$. \par For $f_2(t)$, using the unitarity of $e^{itH}e^{-itH_a}$ on $\s^2_x(\R^9)$, we have 
\begin{align}
\lim\limits_{t\to \infty} \| f_2(t)\|_{\s^2_x(\R^9)}=&\lim\limits_{t\to \infty} \| e^{itH}e^{-itH_a}(1-F_{c,l,\alpha }(x_l,t,P_l))f(x)\|_{\s^2_x(\R^9)}\nonumber\\
=&\lim\limits_{t\to \infty} \| (1-F_{c,l,\alpha }(x_l,t,P_l))f(x)\|_{\s^2_x(\R^9)}\nonumber\\
=&0.
\end{align}
Hence, we can conclude that the limit
\eq
\lim\limits_{t\to \infty} e^{itH}e^{-itH_a}f=f_1(\infty)
\eeq
exists in $\s^2_x(\R^9)$. Thus, $\Omega_{a}$ exists on $\s^2_x(\R^9)$ when $a=(12)(3)$. Similarly, we get the existence of $\Omega_{a^2},$ $\Omega_{a^3}$ and  $\Omega_{free}$, where $a^2=(13)(2)$ and $a^3=(1)(23)$. Furthermore, we get the existence of $\Omega^a$, for all $a=(jk)(l)\in L$, on $\s^2_x(\R^9)$. As a result, we obtain \eqref{cor 1: 1}, \eqref{cor 1: 2} and \eqref{cor 1: 3}, which completes the proof.
\end{proof}
Let
\begin{align}
P_\mu:=& s\text{-}\lim\limits_{t\to \infty} e^{itH}\left(\Pi_{l=1}^{l=3} e^{-itH_a}(1-F_{c,l,\alpha}(x_l,t,P_l))e^{itH_a}\right)e^{-itH}
\nonumber\\
=& s\text{-}\lim\limits_{t\to \infty} e^{itH}\left(\Pi_{l=1}^{l=3} e^{-it\omega_l(P_l)}(1-F_{c,l,\alpha}(x_l,t,P_l))e^{it\omega_l(P_l)}\right)e^{-itH}
\end{align}
on $\s^2_x(\R^9)$ for $a=(jk)(l)\in L$. Here, it is worth noting that for all $a=(jk)(l)\in L$, we have 
\eq
e^{-itH_a}(1-F_{c,l,\alpha}(x_l,t,P_l))e^{itH_a}=e^{-it\omega_l(P_l)}(1-F_{c,l,\alpha}(x_l,t,P_l))e^{it\omega_l(P_l)}.\label{intro: id10}
\eeq 
\begin{proposition}\label{Prop: Pmu: exist}If Assumptions \ref{asp: short}, \ref{asp: critical} and \ref{asp: omega} are satisfied, then $P_\mu$ exists.
    
\end{proposition}
\begin{proof}
    Based on Proposition \ref{cor 1}, we have
    \eq
    s\text{-}\lim\limits_{t\to \infty} e^{itH}e^{-it\omega_l(P_l)}F_{c,l,\alpha}e^{it\omega_l(P_l)}e^{-itH}=s\text{-}\lim\limits_{t\to \infty} e^{itH}e^{-itH_a}F_{c,l,\alpha}e^{itH_a}e^{-itH}\label{E. eq1}
    \eeq
    exists on $\s^2_x(\R^9)$ for all $a=(jk)(l)\in L$. \eqref{E. eq1} implies that 
    \eq
 s\text{-}\lim\limits_{t\to \infty} e^{itH}e^{-it\omega_l(P_l)}\left(1-F_{c,l,\alpha}(x_l,t,P_l)\right)e^{it\omega_l(P_l)}e^{-itH}
    \eeq
    exists on $\s^2_x(\R^9)$. Therefore, $P_\mu$ exists on $\s^2_x(\R^9)$ since
    \begin{align}
     &e^{itH}\left(\Pi_{l=1}^{l=3} e^{-it\omega_l(P_l)}(1-F_{c,l,\alpha}(x_l,t,P_l))e^{it\omega_l(P_l)}\right)e^{-itH}\nonumber\\
        =&e^{itH} e^{-it\omega_1(P_1)}(1-F_{c,1,\alpha}(x_1,t,P_1))e^{it\omega_1(P_1)}e^{-itH}\times \nonumber\\
        &e^{itH} e^{-it\omega_2(P_2)}(1-F_{c,2,\alpha}(x_2,t,P_2))e^{it\omega_2(P_2)}e^{-itH}\times\nonumber\\
        &e^{itH} e^{-it\omega_3(P_3)}(1-F_{c,3,\alpha}(x_3,t,P_3))e^{it\omega_3(P_3)}e^{-itH},
    \end{align}
   and since
    \begin{align}
     &s\text{-}\lim\limits_{t\to \infty}e^{itH} e^{-it\omega_1(P_1)}(1-F_{c,1,\alpha}(x_1,t,P_1))e^{it\omega_1(P_1)}e^{-itH}\times \nonumber\\
        &s\text{-}\lim\limits_{t\to \infty}e^{itH} e^{-it\omega_2(P_2)}(1-F_{c,2,\alpha}(x_2,t,P_2))e^{it\omega_2(P_2)}e^{-itH}\times\nonumber\\
        &s\text{-}\lim\limits_{t\to \infty}e^{itH} e^{-it\omega_3(P_3)}(1-F_{c,3,\alpha}(x_3,t,P_3))e^{it\omega_3(P_3)}e^{-itH}
    \end{align}
    exists on $\s^2_x(\R^9)$.
\end{proof}
Given $a=(jk)(l), b=(j'k')(l')\in L$, we define the following:
\eq
\Pp^{l'}_l:=s\text{-}\lim\limits_{t\to \infty}e^{itH_b} e^{-itH_a}F_{c,l,\alpha}(x_l,t,P_l)e^{itH_a}e^{-itH_b},\quad \text{ on }\s^2_x(\R^9).
\eeq 
\begin{proposition}\label{cor two-body}If Assumptions \ref{asp: short}, \ref{asp: critical} and \ref{asp: omega} are satisfied, then $\Pp^{l'}_l$ exists for all $a=(jk)(l)$, $b=(j'k')(l')\in L$.
 \end{proposition}
\begin{proof}When $b=a$, we have
\eq
\Pp_{l}^{l'}=s\text{-}\lim\limits_{t\to \infty} F_{c,l,\alpha}(x_l,t,P_l)=1
\eeq
exists on $\s^2_x(\R^9)$. When $a\neq b$, let
\eq
\Pp^{l'}_l(t):=e^{itH_b} e^{-itH_a}F_{c,l,\alpha}(x_l,t,P_l)e^{itH_a}e^{-itH_b}.
\eeq
Recall that 
\eq
V_a(x)=V_{jk}(x_j-x_k).
\eeq
Using Cook's method to expand $\Pp^{l'}_l(t)$, we have 
\begin{align}
    \Pp^{l'}_l(t)=&\Pp^{l'}_l(1)+\int_1^t \p_s[\Pp^{l'}_l(s)]ds\nonumber\\
    =&\Pp^{l'}_l(1)+(-i)\int_1^t e^{isH_b}e^{-isH_a}F_{c,l,\alpha}(x_l,s,P_l)e^{isH_a}V_b(x)e^{-isH_b}ds\nonumber\\
    &+i\int_1^t e^{isH_b}e^{-isH_a}F_{c,l,\alpha}(x_l,s,P_l)e^{isH_a}V_a(x)e^{-isH_b}ds\nonumber\\
    &+(-i)\int_1^t e^{isH_b}V_a(x)e^{-isH_a}F_{c,l,\alpha}(x_l,s,P_l)e^{isH_a}e^{-isH_b}ds\nonumber\\
    &+i\int_1^t e^{isH_b}V_b(x)e^{-isH_a}F_{c,l,\alpha}(x_l,s,P_l)e^{isH_a}e^{-isH_b}ds\nonumber\\
    &+\int_1^t e^{isH_b}e^{-isH_a}\p_s[F_{c,l,\alpha}(x_l,s,P_l)]e^{isH_a}e^{-isH_b}ds\nonumber\\
    =:&\Pp^{l'}_l(1)+\left(\sum\limits_{j=1}^4\int_1^tds \Pp_{l,j}^{l'}(s)\right)+\int_1^tds \Pp_{l,+}^{l'}(s).
\end{align}
By employing a similar argument to the one used for $f_{1,in1}(s)$ in \eqref{example in}, we have that
\eq
\int_1^\infty \Pp_{l,j}^{l'}(s)ds=s\text{-}\lim\limits_{t\to \infty} \Pp_{l,j}^{l'}(s)ds\quad \text{ on }\s^2_x(\R^9)
\eeq
exists for all $j=1,2,3,4,$ and all $a,b\in L$ where $a\neq b$
. By employing a similar argument to the one used for $f_{1,+}(s)$ in \eqref{May.eq1} and \eqref{example  +}, we have that
\eq
\int_1^\infty \Pp_{l,+}^{l'}(s)ds=s\text{-}\lim\limits_{t\to \infty} \Pp_{l,+}^{l'}(s)ds\quad \text{ on }\s^2_x(\R^9)
\eeq
exists for all $a,b\in L$ where $a\neq b$
. Therefore, we can conclude that
\eq
\Pp_{l}^{l'}=s\text{-}\lim\limits_{t\to \infty}\Pp_{l}^{l'}(t)\quad \text{ on }\s^2_x(\R^9)
\eeq
exists for all $a,b\in L$ where $a\neq b$. We complete the proof.

\end{proof}
Now we prove Proposition \ref{Prop: Pab}.
\begin{proof}[Proof of Proposition \ref{Prop: Pab}]It follows from Propositions \ref{thm:cW1}, \ref{cor 1} and \ref{cor two-body}.
    
\end{proof}


 \section{Proof Outline for Proposition \ref{Prop: Pmu} and Forward/Backward Waves}\label{sec 3}
 In this section, we prove Proposition \ref{Prop: Pmu} by establishing some a priori estimates for $e^{-itH}P_{sc}f$ $( f\in \s^2_x(\R^9))$. Our discussion in this section is based on Assumptions \ref{asp: short} through \ref{asp: subHH}. Throughout this paper, for all $a=(jk)(l)\in L$, we adopt $\Omega_a^*$, $\Omega^{a,*}$ instead of $\Omega_{a,\alpha}^*$, $\Omega_\alpha^{a,*}$, respectively.
\subsection{Outline for Proposition \ref{Prop: Pmu}}\label{sec:outline}
Now let us explain how to prove Proposition \ref{Prop: Pmu}. Choose a scattering state, $\psi\in \s^2_x(\R^9)$. Then $\psi=\psi_{sc}$ with $\|\psi\|_{\s^2_x(\R^9)}=1$. By using Proposition \ref{Rwhole}, we have that for any $M\geq 1$, there exists a sequence of time $\{t_n\}_{n=1}^{n=\infty}$ with $t_n=t_n(M,\psi)\uparrow \infty$ as $n\to \infty$, such that 
\eq
\| \chi(|\tilde{x}|\leq M)e^{-it_nH}\psi\|_{\s^2_x(\R^9)}< \frac{1}{n}. \label{Req1}
\eeq
If one can show that for each $\epsilon>0$, there exist $M=M(\epsilon,\psi)\geq 1$ and $T=T(\epsilon, M, \psi)>0$ such that when $t\geq T$, 
\eq
\| \chi(|\tilde{x}|>M) e^{-itH}\psi- e^{-itH_0}\psi_{free}-\sum\limits_{a=(jk)(l)\in L}e^{-itH_a}\psi_a\|_{\s^2_x(\R^9)}<\epsilon,\label{sec 4main:eq}
\eeq
then based on \eqref{Req1} and \eqref{sec 4main:eq}, we get a sequential AC. The sequential AC implies that 
\eq
e^{it_nH}\left( \Pi_{l=1}^{l=3} e^{-it_n\omega_l(P_l)} (1-F_{c,l,\alpha}(x_l,t,P_l))e^{it_n\omega_l(P_l)}\right)e^{-it_nH}\psi\to 0
\eeq
in $\s^2_x(\R^9)$ for a sequence of time $\{t_n\}_{n=1}^{n=\infty}$. Due to Proposition \ref{Prop: Pmu: exist}, we get $P_\mu P_{sc}=0$. So it suffices to prove \eqref{sec 4main:eq}. 

Now let us explain how to prove \eqref{sec 4main:eq}. Choose $\epsilon>0$. We will specify the choice of $M$ later. For $\chi(|\tilde{x}|\geq M)$, break it into three pieces:
\eq
\chi(|\tilde{x}|\geq M)=\sum\limits_{a=(jk)(l)\in L} g_{M,a}(\tilde{x}),
\eeq
where $ g_{M,a}(\tilde{x}):=\chi(|\tilde{x}|\geq M)g_a(\tilde{x})$ with $\{g_a\}_{a=(jk)(l)\in L}$ a smooth partition of unity which satisfies that
\eq
 |x_j-x_k| \leq 10\min\limits_{a'=(j'k')(l')\in L} |x_{j'}-x_{k'}|.\label{gMa}
\eeq 
Based on the setup of $g_{M,a}$, we break $ \chi(|\tilde{x}|>M) e^{-itH}\psi$ into three pieces accordingly: 
\eq
\chi(|\tilde{x}|>M) e^{-itH}\psi=\sum\limits_{a=(jk)(l)\in L} g_{M,a}(\tilde{x})e^{-itH}\psi.
\eeq
For $ g_{M,a}(\tilde{x})e^{-itH}\psi$ with $a=(jk)(l)\in L$, we claim that there exists $M_a=M_a(\epsilon, \psi)\geq 1$ such that when $M\geq M_a$ and $t\geq T_a\geq 0$ for some $T_a=T_a(\epsilon,M_a,\psi)>0$, we have 
\eq
\| g_{M,a}(\tilde{x})e^{-itH}\psi-g_a(\tilde{x})e^{-itH_0}\psi_{free}-e^{-itH_a}P_b(H_a)\psi_a\|_{\s^2_x(\R^9)}< \epsilon/3\label{claim1}
\eeq
where $ \psi_{free}=\Omega_{free}^*\psi$ and $ \psi_a=\Omega_a^{*}\psi$:
\begin{claim}\label{Claim1}If Assumptions \ref{asp: short}-\ref{asp: subHH} hold, then for all $a=(jk)(l)\in L$, there exists $M_a=M_a(\epsilon, \psi)\geq 1$ such that when $M\geq M_a$ and $t\geq T_a\geq 0$ for some $T_a=T_a(\epsilon,M,\psi)>0$, \eqref{claim1} is valid.
\end{claim}
If Claim \ref{Claim1} is true, then \eqref{sec 4main:eq} is valid by taking $M=\max\limits_{a=(jk)(l)\in L}\{ M_a\}$ and $T=\max\limits_{a=(jk)(l)\in L}\{ T_a\}$. Now let us explain how to prove Claim \ref{Claim1}. Break $g_{M,a}(\tilde{x})e^{-itH}\psi$ into three pieces:
\begin{align}
    g_{M,a}(\tilde{x})e^{-itH}\psi=& g_{M,a}(\tilde{x})\bar{F}_\tau(H_a,\epsilon_1)e^{-itH}\psi+g_{M,a}(\tilde{x})F_\tau(H_a,\epsilon_1)\bar{F}_\tau(H, 10\epsilon_1)e^{-itH}\psi\nonumber\\
    &+ g_{M,a}(\tilde{x})F_\tau(H_a,\epsilon_1)F_\tau(H, 10\epsilon_1)e^{-itH}\psi\nonumber\\
    =:&\psi_{M,a}(x,t)+\psi_{M,a,r1}(x,t)+\psi_{M,a,r2}(x,t),
\end{align}
where $\epsilon_1=\epsilon_1(\epsilon)>0$ is small enough such that 
\eq
\| F_\tau(H,10\epsilon_1)\psi\|_{\s^2_x(\R^9)}< \frac{\epsilon}{10000},\label{B1epsilon}
\eeq
which implies 
\begin{align}
    \| \psi_{M,a,r2}(x,t)\|_{\s^2_x(\R^9)}\leq & \| F_\tau(H,10\epsilon_1)\psi\|_{\s^2_x(\R^9)}\nonumber\\
    < & \frac{\epsilon}{10000}. 
\end{align}
By using Assumption \ref{asp: subHH} part b) and the unitarity of $e^{-itH}$ on $\s^2_x(\R^9)$, due to Lemma \ref{Lem: tauf}, we have 
\begin{align}
    &\| \psi_{M,a,r1}(x,t)\|_{\s^2_x(\R^9)}\nonumber\\
    \leq & \left( \| \langle x_j-x_l\rangle^{-3} \chi(|x_j-x_k|\geq M^{\frac{1}{100}}) g_{M,a}(\tilde{x})\|_{\s^2_x(\R^9)\to \s^2_x(\R^9)}\right.\nonumber\\
    &+\left. \| \langle x_k-x_l\rangle^{-3} \chi(|x_j-x_k|\geq M^{\frac{1}{100}}) g_{M,a}(\tilde{x})\|_{\s^2_x(\R^9)\to \s^2_x(\R^9)}\right)\nonumber\\
    &\times \| F_\tau(H_a,\epsilon_1)\bar{F}_\tau(H, 10\epsilon_1)\|_{\s^2_x(\R^9)\to \s^2_a(\R^9)}\|e^{-itH}\psi\|_{\s^2_x(\R^9)}\nonumber\\
    \leq & \frac{C_0}{M^3}\label{psiMar1}
\end{align}
for some $C_0=C_0(\epsilon_1)>0$. Then take $M_0= \left(10000(C_0+1)/\epsilon\right)^{1/3}$ and we have that when $M\geq M_0$,
\eq
\| \psi_{M,a,r1}(x,t)\|_{\s^2_x(\R^9)}< \frac{\epsilon}{10000}.
\eeq
For $\psi_{M,a}(x,t)$, break it into three pieces
\begin{align}
\psi_{M,a}(x,t)=&\chi(|x_j-x_k|\geq M^{\frac{1}{100}}) g_{M,a}(\tilde{x})\bar{F}_\tau(H_a,\epsilon_1)e^{-itH}\psi\nonumber\\
&+\chi(|x_j-x_k|< M^{\frac{1}{100}}) g_{M,a}(\tilde{x})P_c(H^a)\bar{F}_\tau(H_a,\epsilon_1)e^{-itH}\psi\nonumber\\
&+\chi(|x_j-x_k|< M^{\frac{1}{100}}) g_{M,a}(\tilde{x})P_b(H^a)\bar{F}_\tau(H_a,\epsilon_1)e^{-itH}\psi\nonumber\\
=&\psi_{M,a,f}(x,t)+\psi_{M,a,\text{pre-f}}(x,t)+\psi_{M,a,b}(x,t),
\end{align}
where 
\eq
\psi_{M,a,f}(x,t):=\chi(|x_j-x_k|\geq M^{\frac{1}{100}}) g_{M,a}(\tilde{x})\bar{F}_\tau(H_a,\epsilon_1)e^{-itH}\psi,
\eeq
\eq
\psi_{M,a,\text{pre-f}}(x,t):=\chi(|x_j-x_k|< M^{\frac{1}{100}}) g_{M,a}(\tilde{x})P_c(H^a)\bar{F}_\tau(H_a,\epsilon_1)e^{-itH}\psi,
\eeq
and 
\eq
\psi_{M,a,b}(x,t):=\chi(|x_j-x_k|< M^{\frac{1}{100}}) g_{M,a}(\tilde{x})P_b(H^a)\bar{F}_\tau(H_a,\epsilon_1)e^{-itH}\psi.
\eeq
Here, $\frac{1}{100}$ represents a small positive number.  
\begin{lemma}\label{Lem: free}For each $R\geq 1$, we let $\chi_R(\tilde{x})$ be some cut-off function satisfying that for all $(j'k')(l')\in L$,  
\eq
|x_{j'}-x_{k'}|\geq R \quad \text{ in the support of }\chi_R(\tilde{x}).
\eeq
Assume that $\psi\in \s^2_x(\R^9)$ with $\| \psi\|_{\s^2_x(\R^9)}=1$ and $\psi=\bar{F}_\tau(H,c_0)\psi$ for some $c_0>0$. If Assumptions \ref{asp: short}-\ref{asp: subHH} hold, then given $\epsilon>0$ there exists a $R_{f}=R_f(\epsilon,\epsilon_1, \psi)\geq 1$ such that when $R\geq R_f$ and $t\geq T_f$  for some $T_f=T_f(\epsilon,\epsilon_1, R_f,\psi)>0$, 
\eq
\|  \chi_R( \tilde{x})\bar{F}_\tau(H_a,\epsilon_1)e^{-itH}\psi-e^{-itH_0}\psi_{free}  \|_{\s^2_x(\R^9)}<\frac{\epsilon}{10},
\eeq
holds true, where $\psi_{free}=\Omega_{free}^*\psi\in \s^2_x(\R^9)$. 

\end{lemma}
We defer the proof for Lemma \ref{Lem: free} to section \ref{section: proof free}. By using Lemma \ref{Lem: free}, we have that by taking $M_1=M^{100}_f$ and $T_1=T_f$, when $ M\geq M_1=M_f$ and $t\geq T_1=T_f$, 
\eq
\|\psi_{M,a,f}(x,t)-  g_a(\tilde{x})e^{-itH_0}\psi_{free}  \|_{\s^2_x(\R^9)}<\frac{\epsilon}{100}.\label{maineq: psif}
\eeq
For $\psi_{M,a,\text{pre-f}}(x,t)$, we let $\tilde{g}_{M,a,s}(\tilde{x})$ and $\tilde{g}_{M,a,l}(\tilde{x})$ be two cut-off functions which satisfy that 
\eq
\tilde{g}_{M,a,s}(\tilde{x}):=\chi(|x_j-x_k|\leq M^{1/100}) \chi(|x_j-x_l|\geq M^{1/4})\chi(|x_k-x_l|\geq M^{1/4})
\eeq
and
\eq
\tilde{g}_{M,a,l}(\tilde{x}):=\chi(|x_j-x_k|> M^{1/100}) \chi(|x_j-x_l|\geq M^{1/4})\chi(|x_k-x_l|\geq M^{1/4}).
\eeq
Let 
\eq
\bar{\chi}_a(\tilde{x}, M):=1-\chi(|x_j-x_l|\geq M^{1/4})\chi(|x_k-x_l|\geq M^{1/4}),\quad M\geq 1, a=(jk)(l)\in L.
\eeq
We break $\psi_{M,a,\text{pre-f}}(x,t)$ into three pieces:
\begin{align}
&\psi_{M,a,\text{pre-f}}(x,t)\nonumber\\
=& \chi(|x_j-x_k|< M^{\frac{1}{100}}) g_{M,a}(\tilde{x}) \bar{F}_\tau(H_a,\epsilon_1)P_c(H^a)e^{i\sqrt{M}H}\tilde{g}_{M,a,s}(\tilde{x})e^{-i(t-\sqrt{M})H}\psi\nonumber\\
&+\chi(|x_j-x_k|< M^{\frac{1}{100}}) g_{M,a}(\tilde{x})\bar{F}_\tau(H_a,\epsilon_1) P_c(H^a)e^{i\sqrt{M}H}\tilde{g}_{M,a,l}(\tilde{x})e^{-i(t-\sqrt{M})H}\psi\nonumber\\
&+ \chi(|x_j-x_k|< M^{\frac{1}{100}}) g_{M,a}(\tilde{x})\bar{F}_\tau(H_a,\epsilon_1) P_c(H^a)e^{i\sqrt{M}H}\bar{\chi}_a(\tilde{x},M)e^{-i(t-\sqrt{M})H}\psi \nonumber\\
=:&\psi_{M,a,\text{pre-f},1}(x,t)+\psi_{M,a,\text{pre-f},2}(x,t)+\psi_{M,a,\text{pre-f},3}(x,t).
\end{align}
\begin{lemma}\label{Lem: pref1}If Assumptions \ref{asp: short}, \ref{asp: critical}, \ref{asp: subH} and \ref{asp: subHH} hold, then for all $M\geq 1$, all $\epsilon>0$ and all $a=(jk)(l)\in L$, we have 
\begin{align}
&\| \chi(|x_j-x_k|<M^{\frac{1}{100}})g_{M,a}(\tilde{x})\bar{F}_\tau(H_a,\epsilon) P_c(H^a)e^{-i\sqrt{M}H}\chi(|x_j-x_k|<M^{\frac{1}{100}})\|_{\s^2_x(\R^9)\to \s^2_x(\R^9)}\nonumber\\
\leq& \frac{C}{\sqrt{M}}
\end{align}
for some $C=C(\epsilon, \opnorm{V})>0$.
\end{lemma}
\begin{proof}Let 
\eq
A_{a,M,\epsilon}:=\chi(|x_j-x_k|<M^{\frac{1}{100}})g_{M,a}(\tilde{x})F_\tau(H_a,\epsilon) P_c(H^a)e^{-i\sqrt{M}H}\chi(|x_j-x_k|<M^{\frac{1}{100}}).
\eeq
By using Duhamel's formula to expand $e^{-i\sqrt{M}H}$ in $\psi_{a,M,\epsilon}(x,t)$ with respect to $e^{-itH_a}$, we have
\begin{align}
    &A_{a,M,\epsilon}\nonumber\\
    =& \chi(|x_j-x_k|<M^{\frac{1}{100}})g_{M,a}(\tilde{x})\bar{F}_\tau(H_a,\epsilon) P_c(H^a)e^{-i\sqrt{M}H_a}\chi(|x_j-x_k|<M^{\frac{1}{100}})\nonumber\\
    &+(-i)\int_0^{\sqrt{M}}ds\chi(|x_j-x_k|<M^{\frac{1}{100}})g_{M,a}(\tilde{x})\bar{F}_\tau(H_a,\epsilon) P_c(H^a)e^{-i(\sqrt{M}-s)H_a}V_{jl}(x_j-x_l)\nonumber\\
    &\times e^{-isH}\chi(|x_j-x_k|<M^{\frac{1}{100}})\nonumber\\
    &+(-i)\int_0^{\sqrt{M}}ds\chi(|x_j-x_k|<M^{\frac{1}{100}})g_{M,a}(\tilde{x})\bar{F}_\tau(H_a,\epsilon) P_c(H^a)e^{-i(\sqrt{M}-s)H_a}V_{kl}(x_k-x_l)\nonumber\\
    &\times e^{-isH}\chi(|x_j-x_k|<M^{\frac{1}{100}})\nonumber\\
    =:&A_{a,M,\epsilon,1}+A_{a,M,\epsilon,2}+A_{a,M,\epsilon,3}. 
\end{align}
By using Lemma \ref{Lem: subdecay}, we have 
\eq
\| A_{a,M,\epsilon,1}\|_{\s^2_x(\R^9)\to \s^2_x(\R^9)}\leq \frac{C}{M^{3/4}}\times M^{\frac{1}{25}}
\eeq
for some $C=C(\epsilon, \sup\limits_{v_T\in \R^3} |\tau(v_T)|)>0 $. Additionally, due to Lemma \ref{Lem: tauf}, we have 
\eq
\sup\limits_{v_T\in \R^3} |\tau(v_T)|<\infty.
\eeq
Therefore, we have
\eq
\| A_{a,M,\epsilon,1}\|_{\s^2_x(\R^9)\to \s^2_x(\R^9)}\leq \frac{C}{\sqrt{M}}\label{estimate: Aa1}
\eeq
for some $C=C(\epsilon)>0 $. For $A_{a,M,\epsilon,2}$ and $A_{a,M,\epsilon,3}$, by using Assumption \ref{asp: subHH} part c), Assumption \ref{asp: subH} part 2 and the unitarity of $e^{-isH}$ on $\s^2_x(\R^9)$, we have 
\begin{align}
    &\|A_{a,M,\epsilon,2} \|_{\s^2_x(\R^9)\to \s^2_x(\R^9)}\nonumber\\
    \leq & \| \langle x_j-x_l\rangle^{-4}g_{M,a}(\tilde{x})\|_{\s^2_x(\R^9)\to \s^2_x(\R^9)}\int_0^{\sqrt{M}}ds \| \langle x_j-x_l\rangle^4 P_c(H^a)\langle x_j-x_l\rangle^{-4}\|_{\s^2_x(\R^9)\to \s^2_x(\R^9)}\nonumber\\
    & \times \| \langle x_j-x_l\rangle^4 F_\tau(H_a,\epsilon)e^{-i(\sqrt{M}-s)H_a}\langle x_j-x_l\rangle^{-4}\|_{\s^2_x(\R^9)\to \s^2_x(\R^9)}\opnorm{V}\nonumber\\
    \leq & \frac{C}{M^4} \times \int_0^{\sqrt{M}}ds \langle \sqrt{M}-s\rangle^4 \opnorm{V}\nonumber\\
    \leq & \frac{C}{M^{3/2}}\label{estimate: Aa2}
\end{align}
for some $C=C(\epsilon, \opnorm{V})>0$, and similarly, 
\begin{align}
    &\|A_{a,M,\epsilon,3} \|_{\s^2_x(\R^9)\to \s^2_x(\R^9)}\nonumber\\
    \leq & \| \langle x_k-x_l\rangle^{-4}g_{M,a}(\tilde{x})\|_{\s^2_x(\R^9)\to \s^2_x(\R^9)}\int_0^{\sqrt{M}}ds \| \langle x_k-x_l\rangle^4 P_c(H^a)\langle x_k-x_l\rangle^{-4}\|_{\s^2_x(\R^9)\to \s^2_x(\R^9)}\nonumber\\
    & \times \| \langle x_k-x_l\rangle^4 F_\tau(H_a,\epsilon)e^{-i(\sqrt{M}-s)H_a}\langle x_k-x_l\rangle^{-4}\|_{\s^2_x(\R^9)\to \s^2_x(\R^9)}\opnorm{V}\nonumber\\
    \leq & \frac{C}{M^4} \times \int_0^{\sqrt{M}}ds \langle \sqrt{M}-s\rangle^4 \opnorm{V}\nonumber\\
    \leq & \frac{C}{M^{3/2}}\label{estimate: Aa3}
\end{align}
for some $C=C(\epsilon, \opnorm{V})>0$. According to \eqref{estimate: Aa1}, \eqref{estimate: Aa2} and \eqref{estimate: Aa3}, we conclude that 
\eq
\| A_{a,M,\epsilon}\|_{\s^2_x(\R^9)\to 
s^2_x(\R^9)}\leq \frac{C}{\sqrt{M}}
\eeq
for some $C=C(\epsilon, \opnorm{V})>0$. We finish the proof.
    
\end{proof}
For $\psi_{M,a,\text{pre-f},1}(x,t)$, according to Lemma \ref{Lem: pref1} and the unitarity of $e^{-i(t-\sqrt{M})H}$ on $\s^2_x(\R^9)$, we have that
\begin{align}
&\| \psi_{M,a,\text{pre-f},1}(x,t)\|_{\s^2_x(\R^9)}\nonumber\\
\leq & \| \chi(|x_j-x_k|< M^{\frac{1}{100}})g_{M,a}(\tilde{x})P_c(H^a)\bar{F}_\tau(H_a,\epsilon_1)e^{i\sqrt{M}H}\chi(|x_j-x_k|<M^{\frac{1}{100}})\|_{\s^2_x(\R^9)\to \s^2_x(\R^9)}\nonumber\\
&\times \| e^{-i(t-\sqrt{M})H}\psi\|_{\s^2_x(\R^9)}\nonumber\\
\leq & \frac{C_{\epsilon_1}}{\sqrt{M}} \label{Sep1eq1}
\end{align}
for some $C_{\epsilon_1}=C_{\epsilon_1}(\epsilon_1)>0$. By taking $M_2=\frac{10^6(C_{\epsilon_1}+1)^2}{\epsilon^2}$, we conclude that for all $M\geq M_2$, according to \eqref{Sep1eq1}, 
\eq
\| \psi_{M,a,\text{pre-f},1}(x,t)\|_{\s^2_x(\R^9)} < \frac{\epsilon}{1000}. 
\eeq
For $\psi_{M,a,\text{pre-f},2}(x,t)$, we use Lemma \ref{Lem: free} and following fact: given $M\geq 1$ and $\psi_{free}\in \s^2_x(\R^9)$, we have
\eq
\lim\limits_{t\to \infty} \| \chi(|x_j-x_k|<M^{1/100})g_{M,a}(\tilde{x})\bar{F}_\tau(H_a,\epsilon_1)P_c(H^a)e^{-i\sqrt{M}H}e^{-i(t-\sqrt{M}H_0)}\psi_{free}\|_{\s^2_x(\R^9)}=0\label{Sep:ineq2}
\eeq
and 
\eq
\lim\limits_{t\to \infty}\| \tilde{g}_{M,a,s}(\tilde{x})e^{-i(t-\sqrt{M}H_0)}\psi_{free}\|_{\s^2_x(\R^9)}=0.\label{Sep:ineq3}
\eeq
According to \eqref{Sep:ineq2} and \eqref{Sep:ineq3}, we have that there exists $T_{31}>0$ such that when $t\geq T_{31}$, 
\begin{align}
&\| \chi(|x_j-x_k|< M^{\frac{1}{100}}) g_{M,a}(\tilde{x})\bar{F}_\tau(H_a,\epsilon_1) P_c(H^a)e^{i\sqrt{M}H}\tilde{g}_{M,a,l}(\tilde{x})e^{-i(t-\sqrt{M}H_0)}\psi_{free}\|_{\s^2_x(\R^9)}\nonumber\\
\leq & \| \tilde{g}_{M,a,s}(\tilde{x})e^{-i(t-\sqrt{M}H_0)}\psi_{free}\|_{\s^2_x(\R^9)}\nonumber\\
&+\| \chi(|x_j-x_k|<M^{1/100})g_{M,a}(\tilde{x})\bar{F}_\tau(H_a,\epsilon_1) P_c(H^a)e^{-i\sqrt{M}H}e^{-i(t-\sqrt{M}H_0)}\psi_{free}\|_{\s^2_x(\R^9)}\nonumber\\
<& \frac{\epsilon}{100}.
\end{align}
Using Lemma \ref{Lem: free}  by taking $ M_3^{1/100}=R_f$ and $T_{32}=T_f$, we have that when $M\geq M_3$ and $t\geq \max\{T_{32}+\sqrt{M}, T_{31}\}$, 
\begin{align}
    &\| \psi_{M,a,\text{pre-f},2}(x,t)\|_{\s^2_x(\R^9)}\nonumber\\
    \leq & \| \tilde{g}_{M,a,l}(\tilde{x})P_c(H^a)e^{-i(t-\sqrt{M})H}\psi- \tilde{g}_{M,a,l}(\tilde{x})e^{-i(t-\sqrt{M}H_0)}\psi_{free}\|_{\s^2_x(\R^9)}\nonumber\\
    &+\| \chi(|x_j-x_k|<M^{1/100})g_{M,a}(\tilde{x})P_c(H^a)e^{-i\sqrt{M}H}\tilde{g}_{M,a,l}(\tilde{x})e^{-i(t-\sqrt{M}H_0)}\psi_{free}\|_{\s^2_x(\R^9)}\nonumber\\
    <& \frac{\epsilon}{50}.\label{Sep1eq2}
\end{align}
For $\psi_{M,a,\text{pre-f},3}(x,t)$, we use Assumption \ref{asp: subHH} part c). Based on the definition of $g_{M,a}$, we have: for $d=1,2,3,$
\eq
\| \frac{1}{\langle x_j-x_l\rangle^d} \chi(|x_j-x_k|< M^{\frac{1}{100}}) g_{M,a}(\tilde{x})\|_{\s^2_x(\R^9)\to \s^2_x(\R^9)}\leq \frac{C_0^d}{M^d}\label{C0eq1}
\eeq
and 
\eq
\| \frac{1}{\langle x_k-x_l\rangle^d} \chi(|x_j-x_k|< M^{\frac{1}{100}}) g_{M,a}(\tilde{x})\|_{\s^2_x(\R^9)\to \s^2_x(\R^9)}\leq \frac{C_0^d}{M^d}\label{C0eq2}
\eeq
for some constant $C_0>0$. Let $C_{c,N}$ be as in Assumption \ref{asp: subHH} part c). According to Lemma \ref{Lem: tauf}, by taking $c=\epsilon_1 $, $N=\max\limits_{v_T\in \R^3}|\tau(v_T)|$ and 
\eq
M_4=\left(\frac{2000(1+C_{c,N})(C_0+1)}{\epsilon}\right)^4, 
\eeq
where $C_0$ is as in \eqref{C0eq1} and \eqref{C0eq2}, we have that when $M\geq M_4$, by using the unitarity of $e^{-i(t-\sqrt{M}H)}$ on $\s^2_x(\R^9)$, we have
\begin{align}
    &\| \psi_{M,a,\text{pre-f},3}(x,t)\|_{\s^2_x(\R^9)}\nonumber\\
    \leq & \frac{C_0}{M} \| \langle x_j-x_l\rangle P_c(H^a)e^{-i\sqrt{M}H} \chi(|x_j-x_l|<M^{1/4} )\|_{\s^2_x(\R^9)\to \s^2_x(\R^9)}\|\psi\|_{\s^2_x(\R^9)}\nonumber\\
    & +\frac{C_0}{M} \| \langle x_k-x_l\rangle P_c(H^a)e^{-i\sqrt{M}H} \chi(|x_k-x_l|<M^{1/4} )\|_{\s^2_x(\R^9)\to \s^2_x(\R^9)}\|\psi\|_{\s^2_x(\R^9)}\nonumber\\
    \leq & \frac{2C_0 C_{c,N} \sqrt{M}\langle M^{1/4}\rangle}{M}\nonumber\\
    <& \frac{\epsilon}{100}. \label{Sep1eq3}
\end{align}
According to \eqref{Sep1eq1}, \eqref{Sep1eq2} and \eqref{Sep1eq3}, we have that by taking $M_5=\max\limits_{d=2,3,4}M_d$, when $M\geq M_5$ and $t\geq T_5=\max\{T_{32}+\sqrt{M}, T_{31}\}$,
\eq
\| \psi_{M,a,\text{pre-f}}(x,t)\|_{\s^2_x(\R^9)}< \frac{\epsilon}{10}.\label{maineq: psipref}
\eeq
\begin{lemma}\label{mainLem: Pba}Let $N^a$ be as in \eqref{def: Na}. If Assumptions \ref{asp: short}-\ref{asp: subHH} hold, then given $\epsilon>0$ there exist $M_f=M_f(\epsilon,\psi)\geq 1$ and $T_f=T_f(\epsilon, M_f, \psi)$ such that when $M\geq M_f$ and $t\geq T_f$, we have 
\eq
\| \psi_{M,a,b}(x,t)-P_b(H^a)e^{-itH_a} \psi_a(x)\|_{\s^2_x(\R^9)}<\frac{N^a\epsilon}{100}
\eeq
for all $a=(jk)(l)\in L$, where $\psi_a(x)$ is given by $\psi_a(x)=\Omega_{a}^{*}\psi$. 
    
\end{lemma}
We defer the proof of Lemma \ref{mainLem: Pba} to Section \ref{sec: Pb}. By using Lemma \ref{mainLem: Pba} with $M_b=M_f$ and $T_b=T_f$, we obtain that when $M\geq M_b$ and $t\geq T_b$, we have
\eq
\| \psi_{M,a,b}(x,t)-P_b(H^a)e^{-itH_a} \psi_a(x)\|_{\s^2_x(\R^9)}<\frac{\epsilon}{100}.\label{maineq: psib}
\eeq
Thus, according to \eqref{maineq: psif}, \eqref{maineq: psipref} and \eqref{maineq: psib}, we obtain that when $M\geq M_a:=\max\{M_1, M_5, M_b\}$ and $t\geq T_a:=\max\{T_1,T_5, T_b \}$, 
\eq
\| g_{M,a}(\tilde{x})e^{-itH}\psi-e^{-itH_0}\psi_{free}-P_b(H^a)e^{-itH_a}\psi_{a}\|_{\s^2_x(\R^9)}< \frac{\epsilon}{3},
\eeq
with $\psi_{free}=\Omega_{free}^{*}\psi$ and $\psi_{a}=\Omega_{a}^{*}\psi$. Hence, we get Claim \ref{Claim1}. Therefore, once we prove Lemma \ref{Lem: free} and Lemma \ref{mainLem: Pba}, we get $P_\mu P_{sc}=0$ and finish the proof of Proposition \ref{Prop: Pmu}.

\subsection{Forward/backward Propagation waves} 

Before proving Lemmas \ref{Lem: free} and \ref{mainLem: Pba}, we begin by introducing the concept of forward/backward propagation waves. The forward/backward propagation waves are similar to the incoming/outgoing waves which are initiated by Mourre \cite{M1979}. Let $S^2$ be the unit sphere in $\R^3$. Assume that  a class of functions on $S^2$, $\{F^{\hat{h}}(\xi)\}_{\hat{h}\in I}$, is a smooth partition of unity with an index set 
 \eq
 I=\{ \hat{h}_1,\cdots, \hat{h}_N\}\subseteq S^2\label{indexI}
 \eeq
 for some $N\in \N^+$. We also assume that for all $\hat{h}\in I$, 
 \eq
 F^{\hat{h}}(\xi)=\begin{cases}1 & \text{ when }|\xi-\hat{h}|<c\\ 0 & \text{ when }|\xi-\hat{h}|>2c\end{cases}, \quad \xi\in S^2\label{ceq1}
 \eeq
 for some $c\in (0,\frac{1}{200})$. For each $\hat{h}\in I$, we also let 
 \eq
 \tilde{F}^{\hat{h}}: S^2\to \R, \text{ be a smooth cut-off function}, 
 \eeq
 satisfying 
 \eq
 \tilde{F}^{\hat{h}}(\xi)=\begin{cases}1 & \text{ when }|\xi-\hat{h}|<4c\\ 0 & \text{ when }|\xi-\hat{h}|>8c\end{cases},\quad \xi\in S^2.\label{ceq2}
 \eeq
We assume that $c>0$, defined in \eqref{ceq1} and \eqref{ceq2}, is proper enough such that for all $x,q\in \R^3$, 
\eq
F^{\hat{h}}(\hat{x})\tilde{F}^{\hat{h}}(\hat{q})|x+q|\geq \frac{1}{10}(|x|+|q|),\label{Feq1}
\eeq
and
\eq
F^{\hat{h}}(\hat{x})(1-\tilde{F}^{\hat{h}}(\hat{q}))|x-q|\geq \frac{1}{10^6}(|x|+|q|).\label{Feq2}
\eeq
Throughout this note, we write $\hat{h}$ as $h/|h|$ for all $h\in \R^3-\{0\}$. Now let us define the projection on forward/backward propagation set, the set of forward/backward propagation waves, with respect to the phase-space $(r,v)\in \R^{3+3}$: 
\begin{definition}[Projection on the forward/backward propagation set]The projection on the forward propagation set, with respect to $(r,v)$, is defined by
\eq
P^+(r,v):=\sum\limits_{b=1}^N F^{\hat{h}_b}(\hat{r})\tilde{F}^{\hat{h}_b}(\hat{v}),\label{Prv+}
\eeq
and the projection on the backward 
propagation set, by 
\eq
P^-(r,v):=1-P^+(r,v).\label{Prv-}
\eeq
 
\end{definition}
\begin{lemma}\label{Lem: Pprop}Let $\omega: \R^3\to \R, \xi\mapsto \omega(\xi)$, with $\nabla_{\xi}[\omega(\xi)]=v(\xi)$. Let $x=r\in \R^3$. For all $\sigma>0$, all $ s\geq 0$, all $M\geq1$ and all $m\in (0,1)$, $P^\pm(r,v)$, defined in \eqref{Prv+} and \eqref{Prv-}, respectively, satisfy that for all $g\in C_0^{[\sigma]+1}(\R^3)$ with $\text{supp}(g(\xi))\subseteq \{ q_k\in \R^3: dist(q_k, \{q_k: v(q_k)=0 \})\geq m\}$, and for all $G\in C^{[\sigma]+1}_0(\R^9)$,
    \begin{align}
    &\| \bar{F}(|r|\geq M) P^\pm(r,v(P)) e^{\pm is\omega(P)}g(v(P))\langle r\rangle^{-\sigma}\|_{\s^2_r(\R^3)\to\s^2_r(\R^3)}\nonumber\leq\frac{C}{\langle M+ms\rangle^\sigma},\label{Lemmain:goal}
    \end{align}
    for some constant $C=C(m, \| g\|_{C^{[\sigma]+1}_0(\R^3)})>0$. 
    
\end{lemma}
\begin{proof}Choose $f\in \s^2_r(\R^3)$ with $f=\chi(|r|\leq \frac{1}{10^7}(M/2+sm))f$. Then by using Fourier representation, we have 
\begin{align}
    & \bar{F}(|r|\geq M) P^\pm(r,v(P)) e^{\pm is\omega(P)}g(v(P))\langle r\rangle^{-\sigma}f\nonumber\\
    =& c \int_{\R^6} d^3qd^3y\bar{F}(|r|\geq M) P^\pm (r, v(q))e^{\pm i s\omega(q)}g(v(q)) e^{-iq\cdot y} \langle y\rangle^{-\sigma} f(y).
\end{align}
In the support of $P^\pm(r,v(q)) g(v(q))$, based on the condition on $\text{supp}(g(\xi))$, when $|r|\geq M/2$ we have that for all $s\geq 1$, 
\begin{align}
| r\pm sv(q)-y  |\geq& \frac{1}{10^6}(|r|+ s|v(q)| )-\frac{1}{10^7}(M/2+sm)\nonumber\\
\geq & \frac{1}{10^7}(M/2+sm).
\end{align}
Thus, by using the method of non-stationary phase(by taking integartion by parts for $[\sigma]+1$ times), we have 
\begin{align}
\| \bar{F}(|r|\geq M) P^\pm(r,v(P)) e^{\pm is\omega(P)}g(v(P))\langle r\rangle^{-\sigma}f\|_{\s^2_r(\R^3)}\leq &\frac{C\|f\|_{\s^2_r(\R^3)}}{\langle M/2+ms\rangle^\sigma}\nonumber\\
\leq & \frac{C\|f\|_{\s^2_r(\R^3)}}{\langle M+ms\rangle^\sigma}
\end{align}
for some constant $C=C(m,\| g\|_{C_0^{[\sigma]+1}}(\R^3))$.

For $f\in \s^2_r(\R^3)$ with $f=\chi(|r|> \frac{1}{10^7}(M/2+sm))f$, we have 
\begin{align}
    &\|\bar{F}(|r|\geq M) P^\pm(r,v(P)) e^{\pm is\omega(P)}g(v(P))\langle r\rangle^{-\sigma}f \|_{\s^2_r(\R^3)}\nonumber\\
    \leq & \|\langle r\rangle^{-\sigma}f \|_{\s^2_r(\R^3)}\nonumber\\
    \leq & \frac{C}{\langle M/2+sm\rangle^\sigma}\|f\|_{\s^2_r(\R^3)}\nonumber\\
    \leq & \frac{C}{\langle M+sm\rangle^\sigma}\|f\|_{\s^2_r(\R^3)}
\end{align}
for some constant $C>0$. Therefore, for all $f\in \s^2_r(\R^3)$, we have 
\begin{align}
\| \bar{F}(|r|\geq M) P^\pm(r,v(P)) e^{\pm is\omega(P)}g(v(P))\langle r\rangle^{-\sigma}f\|_{\s^2_r(\R^3)}\leq & \frac{C\|f\|_{\s^2_r(\R^3)}}{\langle M+ms\rangle^\sigma}
\end{align}
holds true for some constant $C=C(m,\| g\|_{C_0^{[\sigma]+1}}(\R^3))$. We finish the proof.

\end{proof}
\begin{lemma}\label{Lem: Pprop2}Let $\omega: \R^3\to \R, \xi\mapsto \omega(\xi)$, with $\nabla_{\xi}[\omega(\xi)]=v(\xi)$. Let $x=r\in \R^3$. For all $\sigma>0$, all $ s\geq 0$, all $M\geq1$ and all $m\in (0,1)$, $P^\pm(r,v)$, defined in \eqref{Prv+} and \eqref{Prv-}, respectively, satisfy that for all $g\in C_0^{[\sigma]+1}(\R^3)$ with $\text{supp}(g(\xi))\subseteq \{ q_k\in \R^3: dist(q_k, \{q_k: v(q_k)=0 \})\geq m\}$, and for all $G\in C^{[\sigma]+1}_0(\R^9)$,
    \begin{align}
    &\| \bar{F}(|r|\geq M) P^\pm(r,v(P)) e^{\pm is\omega(P)}g(v(P))F(|r|<\frac{1}{10^{10}}(M+sm))\|_{\s^2_r(\R^3)\to\s^2_r(\R^3)}\nonumber\leq\frac{C}{\langle M+ms\rangle^\sigma},\label{Lemmain:goal0}
    \end{align}
    for some constant $C=C(m, \| g\|_{C^{[\sigma]+1}_0(\R^3)})>0$. 
    
\end{lemma}
\begin{proof}Choose $f\in \s^2_r(\R^3)$ with $f=F(|r|\leq \frac{1}{10^{10}}(M+sm))f$. Then by using Fourier representation, we have 
\begin{align}
    & \bar{F}(|r|\geq M) P^\pm(r,v(P)) e^{\pm is\omega(P)}g(v(P))\langle r\rangle^{-\sigma}f\nonumber\\
    =& c \int_{\R^6} d^3qd^3y\bar{F}(|r|\geq M) P^\pm (r, v(q))e^{\pm i s\omega(q)}g(v(q)) e^{-iq\cdot y}  F(|y|\leq \frac{1}{10^{10}}(M+sm))f(y).
\end{align}
In the support of $P^\pm(r,v(q)) g(v(q))$, based on the condition on $\text{supp}(g(\xi))$, when $|r|\geq M/2$ we have that for all $s\geq 1$, 
\begin{align}
| r\pm sv(q)-y  |\geq& \frac{1}{10^6}(|r|+ s|v(q)| )-\frac{1}{10^{10}}(M+sm)\nonumber\\
\geq & \frac{1}{10^7}(M+sm).
\end{align}
Thus, by using the method of non-stationary phase(by taking integartion by parts for $[\sigma]+1$ times), we have 
\begin{align}
\| \bar{F}(|r|\geq M) P^\pm(r,v(P)) e^{\pm is\omega(P)}g(v(P))f\|_{\s^2_r(\R^3)}\leq &\frac{C\|f\|_{\s^2_r(\R^3)}}{\langle M+ms\rangle^\sigma}
\end{align}
for some constant $C=C(m,\| g\|_{C_0^{[\sigma]+1}}(\R^3))$. We finish the proof.

\end{proof}
\section{Proof of Lemma \ref{mainLem: Pba}}\label{sec: Pb}
Now let us explain how to prove Lemma \ref{mainLem: Pba}. Fix $a=(jk)(l)\in L$ and fix $\epsilon>0$. Recall that
\eq
\psi_{M,a,b}(x,t)=\chi(|x_j-x_k|< M^{\frac{1}{100}}) g_{M,a}(\tilde{x})P_b(H^a)\bar{F}_\tau(H_a,\epsilon_1)e^{-itH}\psi.
\eeq
By using Assumption \ref{asp: Pb}, $\psi_{M,a,b}(x,t)$ can be rewritten as 
\begin{align}
    \psi_{M,a,b}(x,t)=\sum\limits_{d=1}^{N^a}\psi_{M,a,b,d}(x,t),\label{psiMab}
\end{align}
where 
\eq
\psi_{M,a,b,d}(x,t):=\chi(|x_j-x_k|< M^{\frac{1}{100}}) g_{M,a}(\tilde{x})P_{b,d}(H^a)\bar{F}_\tau(H_a,\epsilon_1)e^{-itH}\psi.
\eeq
For each $d\in \{1,\cdots, N^a\}$, we define the forward and backward projections with respect to the flow $e^{-it(\omega_l(P_l)+\lambda_{a,d}(P_j+P_k))}$ as:
\eq
P_{a,d}^\pm :=P^\pm (r,v)
\eeq
where $r=x_l-x_j$ and $v=v_l(P_l)-\nabla_{P_j+P_k}[\lambda_{a,d}(P_j+P_k)]$. We then decompose $\psi_{M,a,b,d}(x,t)$ into two components:
\begin{align}
    \psi_{M,a,b,d}(x,t)=& \chi(|x_j-x_k|< M^{\frac{1}{100}}) g_{M,a}(\tilde{x})P^+_{a,d}P_{b,d}(H^a)\bar{F}_\tau(H_a,\epsilon_1)e^{-itH}\psi\nonumber\\
    &+\chi(|x_j-x_k|< M^{\frac{1}{100}}) g_{M,a}(\tilde{x})P^-_{a,d}P_{b,d}(H^a)\bar{F}_\tau(H_a,\epsilon_1)e^{-itH}\psi\nonumber\\
    =:&\psi_{M,a,b,d}^+(x,t)+\psi_{M,a,b,d}^-(x,t).\label{psiMabd}
\end{align}
For $\psi_{M,a,b,d}^+(x,t)$, we approximate it using $P_b(H^a) e^{-itH_a} \Omega_{a}^{*}\psi$. For $\psi_{M,a,b,d}^-(x,t)$, we use $P_b(H^a) e^{-itH_a} \psi$ as its approximation. Their errors decay in $M$ by using forward/backward propagation estimates with respect to flow $e^{-it(\omega_l(P_l)+\lambda_{a,d}(P_j+P_k))}$. 
\begin{proof}[Proof of Lemma \ref{mainLem: Pba}]
Fix $a=(jk)(l)\in L$ and fix $\epsilon>0$. Choose $\psi\in \s^2_x(\R^9)$ with $\|\psi\|_{\s^2_x(\R^9)}=1$. Following \eqref{psiMab} and \eqref{psiMabd}, it suffices to control $\psi_{M,a,b,d}^\pm(x,t)$. For $\psi_{M,a,b,d}^+(x,t)$, by using Duhamel's formula to expand $e^{-itH}$ in $\psi_{M,a,b,d}^+(x,t)$ with respect to $e^{-itH_a}$, we obtain
    \begin{align}
        &\psi_{M,a,b,d}^+(x,t)\nonumber\\
        =& \chi(|x_j-x_k|< M^{\frac{1}{100}}) g_{M,a}(\tilde{x})P^+_{a,d}P_{b,d}(H^a)\bar{F}_\tau(H_a,\epsilon_1)e^{-itH_a}\Omega_{a}^{*}\psi\nonumber\\
        &+i\int_t^\infty ds \chi(|x_j-x_k|< M^{\frac{1}{100}}) g_{M,a}(\tilde{x})P^+_{a,d}P_{b,d}(H^a)\bar{F}_\tau(H_a,\epsilon_1)e^{-i(t-s)H_a}V_{jl}(x_j-x_l)e^{-isH}\psi\nonumber\\
        &+i\int_t^\infty ds \chi(|x_j-x_k|< M^{\frac{1}{100}}) g_{M,a}(\tilde{x})P^+_{a,d}P_{b,d}(H^a)\bar{F}_\tau(H_a,\epsilon_1)e^{-i(t-s)H_a}V_{kl}(x_k-x_l)e^{-isH}\psi\nonumber\\
        =:& \psi_{M,a,b,d,0}^+(x,t)+\psi_{M,a,b,d,1}^+(x,t)+\psi_{M,a,b,d,2}^+.
    \end{align}
Recall that for $j=1,2,3$, $\tilde{H}(P_j,j):=\tilde{H}_0+V_{jk}(-x_k)$ where $\tilde{H}_0=\omega_l(P_l)+\omega_k(P_k)+\omega_{j}(P_j-P_l-P_k) $. Let
    \begin{align}
    \bar{F}_{\tau,j}(\tilde{H}(P_j,j), \epsilon_1):=&F(|\tilde{H}(P_j,j)-\tau(P_j) |>\epsilon_1)F(|P_j|\leq \frac{1}{\epsilon_1})F(\tilde{H}(P_j,j)<\frac{1}{\epsilon_1}).
    \end{align}
    Due to the definition of $\tau(q_3), q_3\in \R^3$ and the continuity of $v_l(q_l)-\nabla_{\eta-q_l} \lambda_{a,d}(\eta-q_l)$, we have that for each $\epsilon_1>0$, there exists $\epsilon_2=\epsilon_2(\epsilon_1)>0$ such that  
    \begin{align}
    \bar{F}_{\tau,j}(\tilde{H}(P_j,j), \epsilon_1)=&F(|v_l(P_l)-\nabla_{P_j-P_l}[\lambda_{a,d}(P_j-P_l)]|>\epsilon_2)F_{\tau,j}(\tilde{H}(P_j,j), \epsilon_1)\nonumber\\
    &\times F(|v_l(P_l)-\nabla_{P_j-P_l}[\lambda_{a,d}(P_j-P_l)]|\leq \frac{1}{\epsilon_2}).\label{id: Ftauj}
    \end{align}
By using 
    \eq
    \chi(|x_j-x_k|<M^{1/100})g_{M,a}(\tilde{x})=\chi(|x_j-x_k|<M^{1/100})g_{M,a}(\tilde{x})\chi(|x_l-x_j|>\frac{M}{100}),
    \eeq
    and by using
    \begin{align}
    &\chi(|x_l-x_j|>\frac{M}{100})P_{a,d}^+P_{b,d}(H^a)\bar{F}_\tau(H_a,\epsilon_1)e^{-i(t-s)H_a}V_{jl}(x_j-x_l)\nonumber\\
    =&e^{-iP_l\cdot x_j}e^{-iP_l\cdot x_k} \chi(|x_l|>\frac{M}{100})P_{a,j,d}^+e^{-i(t-s)(\omega_l(P_l)+\lambda_{a,d}(P_j-P_l))}P_{b,d}(H^a(P_j-P_l))\bar{F}_\tau(\tilde{H}(P_j,j),\epsilon_1)\nonumber\\
    &\times V_{jl}(-x_l)e^{iP_l\cdot x_j}e^{iP_l\cdot x_k},
    \end{align}
    where $H_{a,j}:=\omega_l(P_l)+\omega_k(P_k)+\omega_j(P_j-P_l)+V_{jk}(x_j-x_k)$ and 
    \eq
    P_{a,j,d}^+:=P^+(r,v)
    \eeq
    with $r=x_l$ and $v=v_l(P_l)-\nabla_{P_j-P_l}[\lambda_{a,d}(P_j-P_l)]$, according to \eqref{id: Ftauj} and Lemma \ref{Lem: Pprop}, we have that for $u\geq 0$
    \begin{multline}\label{Pbeq1}
    \| \chi(|x_l|> \frac{M}{100})P^\pm_{a,j,d}e^{\pm iu(\omega_l(P_l)+ \lambda_{a,d}(P_j-P_l))}\\
   \times F(|v_l(P_l)-\nabla_{P_j-P_l}[\lambda_{a,d}(P_j-P_l)]|>\epsilon_2)\langle x_l\rangle^{-2}  \|_{\s^2_x(\R^9)\to \s^2_x(\R^9)}\lesssim_{\epsilon_2} \frac{1}{\langle M+u\rangle^2 },
    \end{multline}
which also implies that 
\begin{multline}\label{Pbeq2}
    \| \chi(|x_l-x_j|> \frac{M}{100})P^\pm_{a,d}e^{\pm iu(\omega_l(P_l)+ \lambda_{a,d}(P_j+P_k))}\\
   \times F(|v_l(P_l)-\nabla_{P_j+P_k}[\lambda_{a,d}(P_j+P_k)]|>\epsilon_2)\langle x_l-x_j\rangle^{-2}  \|_{\s^2_x(\R^9)\to \s^2_x(\R^9)}\lesssim_{\epsilon_2} \frac{1}{\langle M+u\rangle^2 },
    \end{multline}
By using \eqref{Pbeq1}, \eqref{Pb: eq}, Assumption \ref{asp: short} and the unitarity of $e^{-isH}$ on $\s^2_x(\R^9)$, we have 
    \begin{multline}
        \| \psi_{M,a,b,d,1}^+(x,t)\|_{\s^2_x(\R^9)}\leq \int_t^\infty ds\|  \chi(|x_l|> \frac{M}{100})P^\pm_{a,j,d}\\
        \times e^{ i(s-t)(\omega_l(P_l)+ \lambda_{a,d}(P_j-P_l))} F(|v_l(P_l)-\nabla_{P_j-P_l}[\lambda_{a,d}(P_j-P_l)]|>\epsilon_2)\langle x_l\rangle^{-2}\|_{\s^2_x(\R^9)\to \s^2_x(\R^9)}\\
        \times \| \langle x_l\rangle^2 P_{a,d}(H^a(P_j-P_l))\langle x_l\rangle^{-2} \|_{\s^2_x(\R^9)\to \s^2_x(\R^9)}\| \langle x_l\rangle^2 V_{jl}(-x_l)\|_{\s^\infty_{x_l}(\R^3)}\| \psi\|_{\s^2_x(\R^9)}\\
        \leq \int_t^\infty ds\frac{C}{\langle M+s-t\rangle^2}\|\langle \eta\rangle^4 V_{jl}(\eta)\|_{\s^\infty_\eta(\R^3)}\\
        \leq \frac{C_1}{M} \opnorm{V}\label{0psiMabd1}
    \end{multline}
for some $C_1=C_1(\epsilon_2)>0$. Since $\epsilon_2=\epsilon_2(\epsilon_1)>0$ and $\epsilon_1=\epsilon_1(\epsilon)>0$, we have $\epsilon_2=\epsilon_2(\epsilon)>0$ and therefore $C_1$ in \eqref{psiMabd1} is some positive constant which only depends on $\epsilon$. Similarly, by using Assumption \ref{asp: subH} part 3, we have 
\begin{multline}
        \| \psi_{M,a,b,d,2}^+(x,t)\|_{\s^2_x(\R^9)}\leq \int_t^\infty ds\|  \chi(|x_l|> \frac{M}{100})P^\pm_{a,j,d}\\
        \times e^{ i(s-t)(\omega_l(P_l)+ \lambda_{a,d}(P_j-P_l))} F(|v_l(P_l)-\nabla_{P_j-P_l}[\lambda_{a,d}(P_j-P_l)]|>\epsilon_2)\langle x_l\rangle^{-2}\|_{\s^2_x(\R^9)\to \s^2_x(\R^9)}\\
        \times \| \langle x_l\rangle^2 P_{a,d}(H^a(P_j-P_l))\langle x_l-x_k\rangle^{-2} \|_{\s^2_x(\R^9)\to \s^2_x(\R^9)}\| \langle \eta\rangle^2 V_{kl}(\eta)\|_{\s^\infty_{\eta}(\R^3)}\| \psi\|_{\s^2_x(\R^9)}\\
        \leq \int_t^\infty ds\frac{C}{\langle M+s-t\rangle^2}\opnorm{V}\\
        \leq \frac{C_2}{M} \opnorm{V}\label{psiMabd1}
    \end{multline}
for some $C_2=C_2(\epsilon)>0$. Additionally, due to the definition of $g_{M,a}(\tilde{x})$ in \eqref{gMa}, we have 
\eq
1-g_{M,a}(\tilde{x})=(1-g_{M,a}(\tilde{x}))\chi(|x_j-x_k|\geq \frac{M}{1000}).\label{gMajk}
\eeq
Therefore, by using intertwining identity 
\eq
\bar{F}_\tau(H_a,\epsilon_1)\Omega_{a}^{*}\psi=\Omega_{a}^{*}\bar{F}_\tau(H,\epsilon_1)\psi,
\eeq
unitarity of $e^{-itH_a}$ on $\s^2_x(\R^9)$, Assumption \ref{asp: subH}, \eqref{B1epsilon} and \eqref{gMajk}, we have
\begin{align}
    &\| \psi_{M,a,b,d,0}^+(x)-P_{b,d}(H^a)e^{-itH_a}\Omega_{a}^{*}\psi \|_{\s^2_x(\R^9)}\nonumber\\
    \leq & \|(1-\chi(|x_j-x_k|<M^{\frac{1}{100}})g_{M,a}(\tilde{x}))P^+_{a,d}P_{b,d}(H^a)\bar{F}_\tau(H_a,\epsilon_1)e^{-itH_a}\Omega_{a}^{*}\psi \|_{\s^2_x(\R^9)}+\nonumber\\
    & \|P^+_{a,d}P_{b,d}(H^a)F_\tau(H_a,\epsilon_1)e^{-itH_a}\Omega_{a}^{*}\psi \|_{\s^2_x(\R^9)}+\|P^-_{a,d}P_{b,d}(H^a)e^{-itH_a}\Omega_{a}^{*}\psi \|_{\s^2_x(\R^9)}\nonumber\\
    \leq & \frac{C}{M^4} \| \langle x_j-x_k \rangle^4 P_{b,d}(H^a) \|_{\s^2_x(\R^9)\to \s^2_x(\R^9)}\|\psi\|_{\s^2_x(\R^9)}+\frac{\epsilon}{10000}+\|P^-_{a,d}P_{b,d}(H^a)e^{-itH_a}\Omega_{a}^{*}\psi \|_{\s^2_x(\R^9)}\nonumber\\
    \leq & \frac{C_0}{M^4}+\frac{\epsilon}{10000}+\|P^-_{a,d}P_{b,d}(H^a)e^{-itH_a}\Omega_{a}^{*}\psi \|_{\s^2_x(\R^9)}
\end{align}
for some constant $C_0>0$. Using \eqref{Pbeq2}, we obtain 
\eq
\lim\limits_{t\to \infty}\| P_{a,d}^-P_{b,d}(H^a) e^{-itH_a}\Omega_{a}^{*}\psi\|_{\s^2_x(\R^9)} =0.
\eeq
Therefore, there exists $T_{ad1}=T_{ad1}(M,\epsilon)>0$ such that when $t\geq T_{ad0}$, 
\eq
\| P_{a,d}^-P_{b,d}(H^a) e^{-itH_a}\Omega_{a}^{*}\psi\|_{\s^2_x(\R^9)} <\frac{\epsilon}{10000}.
\eeq
Hence, by taking 
\eq
M_{ad1}=\max\limits_{d=1,2,3} M_d
\eeq
where 
\eq
M_1:=\frac{10000C_1(1+\opnorm{V})}{\epsilon },
\eeq
\eq
M_2:=\frac{10000C_2(1+\opnorm{V})}{\epsilon }
\eeq
and 
\eq
M_3:= \left(\frac{10000C_0}{\epsilon}\right)^{1/4},
\eeq
when $M\geq M_{ad1}$ and $t\geq T_{ad1}(M,\epsilon)$, we have
\eq
\| \psi^+_{M,a,b,d}(x)-P_{b,d}(H^a)e^{-itH_a}\Omega_{a}^{*}\psi\|_{\s^2_x(\R^9)}< \frac{\epsilon}{2500}.\label{psiMabd+}
\eeq
For $\psi_{M,a,b,d}^-(x,t)$, by using Duhamel's formula to expand $e^{-itH}$ in $\psi_{M,a,b,d}^-(x,t)$, we obtain
    \begin{align}
        &\psi_{M,a,b,d}^-(x,t)\nonumber\\
        =& \chi(|x_j-x_k|< M^{\frac{1}{100}}) g_{M,a}(\tilde{x})P^-_{a,d}P_{b,d}(H^a)\bar{F}_\tau(H_a,\epsilon_1)e^{-itH_a}\psi\nonumber\\
        &+(-i)\int_0^t ds \chi(|x_j-x_k|< M^{\frac{1}{100}}) g_{M,a}(\tilde{x})P^-_{a,d}P_{b,d}(H^a)\bar{F}_\tau(H_a,\epsilon_1)e^{-i(t-s)H_a}V_{jl}(x_j-x_l)e^{-isH}\psi\nonumber\\
        &+(-i)\int_0^t ds \chi(|x_j-x_k|< M^{\frac{1}{100}}) g_{M,a}(\tilde{x})P^-_{a,d}P_{b,d}(H^a)\bar{F}_\tau(H_a,\epsilon_1)e^{-i(t-s)H_a}V_{kl}(x_k-x_l)e^{-isH}\psi\nonumber\\
        =:& \psi_{M,a,b,d,0}^-(x,t)+\psi_{M,a,b,d,1}^-(x,t)+\psi_{M,a,b,d,2}^-.
    \end{align}
    By using \eqref{Pbeq2}, we obtain, similarly,
    \eq
    \| \psi^-_{M,a,b,d,1}(x,t)\|_{\s^2_x(\R^9)}\leq \frac{C_3}{M}\opnorm{V}
    \eeq
    and 
    \eq
    \| \psi^-_{M,a,b,d,2}(x,t)\|_{\s^2_x(\R^9)}\leq \frac{C_4}{M}\opnorm{V}
    \eeq
    for some constants $C_3=C_3(\epsilon)>0$ and $C_4=C_4(\epsilon)>0$. Using \eqref{Pbeq2} again, we obtain 
\eq
\lim\limits_{t\to \infty}\| P_{a,d}^-P_{b,d}(H^a) e^{-itH_a}\psi\|_{\s^2_x(\R^9)} =0.
\eeq
Therefore, there exists $T_{ad2}=T_{ad2}(M,\epsilon)>0$ such that when $t\geq T_{ad2}$, 
\eq
\| P_{a,d}^-P_{b,d}(H^a) e^{-itH_a}\psi\|_{\s^2_x(\R^9)} <\frac{\epsilon}{10000}.
\eeq
Hence, by taking 
\eq
M_{ad2}=\max\limits_{d=4,5} M_d
\eeq
where 
\eq
M_4:=\frac{10000C_3(1+\opnorm{V})}{\epsilon },
\eeq
and
\eq
M_5:=\frac{10000C_4(1+\opnorm{V})}{\epsilon },
\eeq
when $M\geq M_{ad2}$ and $t\geq T_{ad2}(M,\epsilon)$, we have 
\eq
\| \psi^-_{M,a,b,d}(x,t)\|_{\s^2_x(\R^9)}< \frac{3\epsilon}{10000}.\label{psiMabd-}
\eeq
According to \eqref{psiMabd+} and \eqref{psiMabd-}, we have that when $M\geq M_{ad}:=\max\{M_{ad1}, M_{ad2}\}$ and $t\geq T_{ad}:=\max\{ T_{ad1}, T_{ad2}\}$, we have 
\eq
\| \psi_{M,a,b,d}(x)-P_{b,d}(H^a)e^{-itH_a}\Omega_{a}^{*}\psi\|_{\s^2_x(\R^9)}< \frac{\epsilon}{1250}.\label{eq: psiMabd}
\eeq
Based on \eqref{eq: psiMabd}, we have that by taking $M_{ab}:=\max\limits_{d=1,\cdots,N^a}M_{ad}$ and $t\geq T_{ad}:=\max\limits_{d=1,\cdots,N^a} T_{ad}$, 
\begin{align}
\| \psi_{M,a,b}(x)-P_b(H^a)e^{-itH_a}\Omega_{a}^{*}\psi\|_{\s^2_x(\R^9)}\leq& \sum\limits_{d=1}^{N^a} \| \psi_{M,a,b,d}(x)-P_{b,d}(H^a)e^{-itH_a}\Omega_{a}^{*}\psi\|_{\s^2_x(\R^9)}\nonumber\\
< & \frac{N^a\epsilon}{1250}
\end{align}
holds true for all $a=(jk)(l)\in L$ and $\psi=\psi_{sc}\in \s^2_x(\R^9)$ with $\| \psi\|_{\s^2_x(\R^9)}=1$. We finish the proof. 
\end{proof}

\section{Proof of Lemma \ref{Lem: free}}\label{section: proof free}
In this section, we prove Lemma \ref{Lem: free}.
\subsection{Outline of the proof} \label{sec: outline}Fix $a=(jk)(l)\in L$ and choose $\psi=\psi_{sc}\in \s^2_x(\R^9)$ with $\| \psi\|_{\s^2_x(\R^9)}=1$ and $\psi=\bar{F}_\tau(H,c_0)\psi$ for some $c_0>0$. Fix $\epsilon>0$. Given $\epsilon>0$, we assume that $\epsilon_1=\epsilon_1(\epsilon)>0$ is smaller enough such that \eqref{B1epsilon} is valid. Let
\eq
\psi_{R,\epsilon_1}(x,t)=\chi_R(\tilde{x})\bar{F}_\tau(H_a,\epsilon_1)e^{-itH}\psi.
\eeq
Break $\psi_{R,\epsilon_1}(x,t)$ into three pieces:
\begin{align}
    \psi_{R,\epsilon_1}(x,t)=&\chi_R(\tilde{x})\bar{F}_\tau(H_0,\epsilon_1/10)e^{-itH}\psi\nonumber\\
    &+\chi_R(\tilde{x})F_\tau(H_0,\epsilon_1/10)\bar{F}_\tau(H_a,\epsilon_1)e^{-itH}\psi+\chi_R(\tilde{x})\bar{F}_\tau(H_0,\epsilon_1/10)F_\tau(H_a,\epsilon_1)e^{-itH}\psi\nonumber\\
    =:&\psi_{R,\epsilon_1,m}(x,t)+\psi_{R,\epsilon_1,r1}(x,t)+\psi_{R,\epsilon_1,r_2}.
\end{align}
For $\psi_{R,\epsilon_1,r_1}(x,t)$, by using Assumption \ref{asp: subHH} {part a)} and the unitarity of $e^{-itH}$ on $\s^2_x(\R^9)$, we have 
\begin{align}
    \| \psi_{R,\epsilon_1,r_1}(x,t)\|_{\s^2_x(\R^9)}\leq & \| \langle x_j-x_k\rangle^{-3}\chi_R(\tilde{x})\|_{\s^2_x(\R^9)\to \s^2_x(\R^9)}\nonumber\\
    &\times \|\langle x_j-x_k\rangle^3 F_\tau(H_0,\epsilon_1/10)\bar{F}_\tau(H_a,\epsilon_1) \|_{\s^2_x(\R^9)\to \s^2_x(\R^9
    )}\|e^{-itH}\psi\|_{\s^2_x(\R^9)}\nonumber\\
    \leq & \frac{C}{R^3}
\end{align}
for some $C=C(\epsilon_1)=C(\epsilon)>0$. Therefore, by taking $R_{r1}=\left(\frac{10000C}{\epsilon}\right)^{1/3}$, when $R\geq R_{r1}$, 
\eq
\| \psi_{R,\epsilon_1,r_1}(x,t)\|_{\s^2_x(\R^9)}< \frac{\epsilon}{10000}.\label{psiRr1}
\eeq
For $\psi_{R,\epsilon_1,r_2}(x,t)$, by using \eqref{B1epsilon} and some process which is similar to what we did for $\psi_{M,a,r1}(x,t)$ in \eqref{psiMar1}, we have that there exists $R_{r2}>0$ such that when $R\geq R_{r2}$, 
\eq
\| \psi_{R,\epsilon_1,r_2}(x,t)\|_{\s^2_x(\R^9)}< \frac{\epsilon}{10000}.\label{psiRr2}
\eeq
Now let us look at $\psi_{R,\epsilon_1,m}(x,t)$. $\psi_{R,\epsilon_1,m}(x,t)$ can be rewritten as 
\eq
\psi_{R,\epsilon_1,m}(x,t)=\sum\limits_{l=1}^3 \psi_{R,\epsilon_1,m,l}(x,t),
\eeq
where $\{\tilde{g}_a(q)\}_{a=(jk)(l)\in L}$ is a smooth partition of unity, which satisfies 
\eq
 |q_j-q_k| \geq \frac{1}{10}\max\limits_{a'=(j'k')(l')\in L} |q_{j'}-q_{k'}|,\label{tgMa}
\eeq 
\eq
F_{\tau,l}(H_0,\epsilon_1/10):= \tilde{g}_a(v_j(P_j)-v_k(P_k))F_\tau(H_0, \epsilon_1/10)
\eeq
and
\eq
\psi_{R,\epsilon_1,m,l}(x,t):=\chi_R(\tilde{x})F_{\tau,l}(H_0,\epsilon_1/10)e^{-itH}\psi.
\eeq
Due to the continuity of $\sum\limits_{d=1}^3 \omega_d(\eta_d)$ in $(\eta_1,\eta_2,\eta_3)$ and due to the definition of $\tau(v_T)$ for $v_T\in \R^3$, we have that there exists $\epsilon_2=\epsilon_2(\epsilon_1)>0$ such that 
\eq
F_{\tau,l}(H_0,\epsilon_1/10)=F(|v_j(P_j)-v_k(P_k)|>\epsilon_2)F_{\tau,l}(H_0,\epsilon_1/10).\label{epsilon2}
\eeq
For $\psi_{R,\epsilon_1,m,l}(x,t)$, we need some projections on new forward/backward propagation waves. For $1\leq j<k\leq 3(j,k\in \N^+)$, taking $(r,v)=(x_j-x_k,v_j(P_j)-v_k(P_k)),$ we define 
\eq
P^\pm_{jk}:=P^\pm(x_j-x_k,v_j(P_j)-v_k(P_k)).
\eeq
Based on the definition of $P^\pm_{jk}, 1\leq j<k\leq 3$, we break $\psi_{R,\epsilon_1,m,l}(x,t)$ into two pieces: 
\begin{align}
    \psi_{R,\epsilon_1,m,l}(x,t)=&\chi_R(\tilde{x})P^{+}_{jk}F_{\tau,l}(H_0,\epsilon_1/10)e^{-itH}\psi+\chi_R(\tilde{x})P^{-}_{jk}F_{\tau,l}(H_0,\epsilon_1/10)e^{-itH}\psi\nonumber\\
    =:&\psi^+_{R,\epsilon_1,m,l}(x,t)+\psi^-_{R,\epsilon_1,m,l}(x,t).
\end{align}
For $\psi^+_{R,\epsilon_1,m,l}(x,t)$, we use $\chi_R(\tilde{x})P_{jk}^+ F_{\tau,l}(H_0,\epsilon_1/10)e^{-itH_0}\Omega_{free}^*\psi$ to approximate it. To be precise, by using Duhamel's formula to expand $e^{-itH}$ in $\psi^+_{R,\epsilon_1,m,l}(x,t)$, we obtain
\begin{align}
    \psi_{R,\epsilon_1,m,l}^+(x,t)=&\chi_R(\tilde{x})P^+_{jk}F_{\tau,l}(H_0,\epsilon_1/10)e^{-itH_0}\Omega_{free}^*\psi\nonumber\\
    &+i\chi_R(\tilde{x})P^+_{jk}F_{\tau,l}(H_0,\epsilon_1/10)\int_t^\infty ds e^{-i(t-s)H_0}V_{jk}(x_j-x_k)e^{-isH}\psi\nonumber\\
    &+i\chi_R(\tilde{x})P^+_{jk}F_{\tau,l}(H_0,\epsilon_1/10)\int_t^\infty ds e^{-i(t-s)H_0}V_{jl}(x_j-x_l)e^{-isH}\psi\nonumber\\
    &+i\chi_R(\tilde{x})P^+_{jk}F_{\tau,l}(H_0,\epsilon_1/10)\int_t^\infty ds e^{-i(t-s)H_0}V_{kl}(x_k-x_l)e^{-isH}\psi\nonumber\\
    =:&\psi_{free,m,l}^+(x,t)+\psi_{m,l,jk}(x,t)+\psi_{m,l,jl}(x,t)+\psi_{m,l,kl}(x,t).
\end{align}
Here, we define $V_{kj}(x_k-x_j)$ as $V_{jk}(x_j-x_k)$ when $k<j$. For $\psi^+_{free,m,l}(x,j)$, since
\eq
\lim\limits_{t\to \infty} \| \chi_R(\tilde{x})P^+_{jk} e^{-itH_0}\Omega_{ free}^*\psi-\tilde{g}_a(v_j(P_j)-v_k(P_k))e^{-itH_0}\Omega_{free}^*\psi\|_{\s^2_x(\R^9)}=0
\eeq
and since 
\begin{align}
    & \|\tilde{g}_a(v_j(P_j)-v_k(P_k)) F_\tau(H_0,\epsilon_1/10) e^{-itH_0}\Omega_{free}^*\psi\|_{\s^2_x(\R^9)}\nonumber\\
    \leq & \| \Omega_{ free}^* F_\tau(H,\epsilon_1/10)\psi\|_{\s^2_x(\R^9)}\nonumber\\
    \leq & \| F_\tau(H,10\epsilon_1)\psi\|_{\s^2_x(\R^9)}\nonumber\\
    < & \frac{\epsilon}{10000},
\end{align}
we conclude that there exists $T_0\geq 1$ such that when $t\geq T_0$, we have 
\eq
\| \psi^+_{free,m,l}(x,t)-\tilde{g}_a(v_j(P_j)-v_k(P_k))e^{-itH_0}\Omega_{free}^*\psi\|_{\s^2_x(\R^9)}< \frac{\epsilon}{1000}.\label{psifreel}
\eeq
Now recall that $\epsilon_2$ is defined in \eqref{epsilon2}. For $\psi_{m,l,jk}(x,t)$, by using Lemma \ref{Lem: Pprop} and the unitarity of $e^{-isH}$ on $\s^2_x(\R^9)$, we have 
\begin{align}
    \| \psi_{m,l,jk}(x,t)\|_{\s^2_x(\R^9)}\leq & \int_t^\infty ds \| \chi_R(\tilde{x})P^+_{jk}\bar{F}_{\tau,l}(H_0, \epsilon_1/10)e^{-i(t-s)H_0}\langle x_j-x_k\rangle^{-2}\|_{\s^2_x(\R^9)\to \s^2_x(\R^9)}\nonumber\\
    &\times\| \langle \eta\rangle^2 V_{jk}(\eta)\|_{\s^\infty_\eta(\R^3)}\| e^{-isH}\psi\|_{\s^2_x(\R^9)}\nonumber\\
    \lesssim_{\epsilon_1,\epsilon_2} &\int_t^\infty \frac{ds}{\langle R+(t-s)\epsilon_2\rangle^2} \| \langle \eta\rangle^2 V_{jk}(\eta)\|_{\s^\infty_\eta(\R^3)}\|\psi\|_{\s^2_x(\R^9)}\nonumber\\
    \lesssim_\epsilon& \frac{1}{R}\| \langle \eta\rangle^2 V_{jk}(\eta)\|_{\s^\infty_\eta(\R^3)}.\label{psifreejk}
\end{align}
Here, we use the fact that $\epsilon_2$ depends on $\epsilon_1$ and $\epsilon_1$ depends on $\epsilon$. For $\psi_{m,l,jl}(x,t)$, we break it into two pieces: 
\begin{align}
   &\psi_{m,l,jl}(x,t)\nonumber\\
   =& i\chi_R(\tilde{x})P^+_{jk}\bar{F}_{\tau,l}(H_0,\epsilon_1/10)\int_t^\infty ds e^{-i(t-s)H_0}F(|x_j-x_k|>\frac{1}{10^{11}}(R+|t-s|\epsilon_2))V_{jl}(x_j-x_l)e^{-isH}\psi\nonumber\\
   &+i\chi_R(\tilde{x})P^+_{jk}\bar{F}_{\tau,l}(H_0,\epsilon_1/10)\int_t^\infty ds e^{-i(t-s)H_0}F(|x_j-x_k|\leq \frac{1}{10^{11}}(R+|t-s|\epsilon_2))V_{jl}(x_j-x_l)e^{-isH}\psi\nonumber\\
   =:&\psi_{m,l,jl,1}(x,t)+\psi_{m,l,jl,2}(x,t).
\end{align}
For $\psi_{m,l,jl,2}(x,t)$, by using Lemma \ref{Lem: Pprop2} and unitarity of $e^{-isH}$ on $\s^2_x(\R^9)$, we have 
\begin{align}
   & \| \psi_{m,l,jl,2}(x,t)\|_{\s^2_x(\R^9)}\nonumber\\
   \leq & \int_t^\infty ds \| \chi_R(\tilde{x})P^+_{jk}\bar{F}_{\tau,l}(H_0, \epsilon_1/10)e^{-i(t-s)H_0}F(|x_j-x_k|\leq \frac{1}{10^{11}}(R+|t-s|\epsilon_2))\|_{\s^2_x(\R^9)\to \s^2_x(\R^9)}\nonumber\\
    &\times\|  V_{jl}(\eta)\|_{\s^\infty_\eta(\R^3)}\| e^{-isH}\psi\|_{\s^2_x(\R^9)}\nonumber\\
    \lesssim_{\epsilon_1,\epsilon_2} &\int_t^\infty \frac{ds}{\langle R+(t-s)\epsilon_2\rangle^2} \|  V_{jl}(\eta)\|_{\s^\infty_\eta(\R^3)}\|\psi\|_{\s^2_x(\R^9)}\nonumber\\
    \lesssim_\epsilon& \frac{1}{R}\|  V_{jk}(\eta)\|_{\s^\infty_\eta(\R^3)}.\label{psimljl2}
\end{align}
For $\psi_{m,l,jl,1}(x,t)$, we need the notion of the projections on following forward/backward propagation waves. For $j,k\in\{1,2,3\}$ with $j\neq k$($j$ is not necessary to be smaller than $k$), taking $(r,v)=(x_j-x_k, v_j(P_j))$, we define 
\eq
\tilde{P}_{jk,j}^\pm:=P^\pm(x_j-x_k,v_j(P_j)).
\eeq
For $\psi_{m,l,jl,1}(x,t,s)$, let $a':=(jl)(k)\in L$,
\begin{align}
   \psi_{m,l,jl,1}^0(x,t,s):=&F(|x_j-x_k|>\frac{1}{10^{11}}(R+|t-s|\epsilon_2)) V_{jl}(x_j-x_l)e^{-isH}\psi
\end{align}
\begin{align}
   \psi_{m,l,jl,1}^+(x,t,s):=&F(|x_j-x_k|>\frac{1}{10^{11}}(R+|t-s|\epsilon_2)) V_{jl}(x_j-x_l)\tilde{P}_{kj,k}^+e^{-isH_{a'}}\Omega_{a'}^{*}\psi
\end{align}
and
\begin{align}
   \psi_{m,l,jl,1}^-(x,t,s):=&F(|x_j-x_k|>\frac{1}{10^{11}}(R+|t-s|\epsilon_2)) V_{jl}(x_j-x_l)\tilde{P}_{kj,k}^-e^{-isH_{a'}}\psi.
\end{align}
\begin{lemma}\label{Lem: two-body}Let $\psi$ be as in Lemma \ref{Lem: free}. If Assumptions \ref{asp: short} and \ref{asp: subH} hold, then for all $\psi\in \s^2_x(\R^9)$ and all $R\geq 1$
\eq
\lim\limits_{t\to \infty}\| \int_t^\infty ds \bar{F}_{\tau,l}(H_0,\epsilon_1/10)e^{-i(t-s)H_0}\psi_{m,l,jl,1}^\pm(x,t,s)\|_{\s^2_{x}(\R^9)}=0,\label{9.11end: eq1}
\eeq
and 
\eq
\lim\limits_{t\to \infty} \| \chi_R(\tilde{x})P^-_{jk}\bar{F}_{\tau,l}(H_0,\epsilon_1/10)\int_0^t ds e^{-i(t-s)H_0}\psi_{m,l,jl,1}^\pm(x,t,s)\|_{\s^2_{x}(\R^9)}=0.\label{9.11end: eq2}
\eeq
\end{lemma}
\begin{proof}Let $\mathcal{S}_x(\R^9)$ be the Schwartz space. Since $\tilde{P}_{kj,k}^\pm $ commutes with $V_{jl}$, by using Assumption \ref{asp: subH} part 2, Fubini's theorem and the unitarity of $e^{-i(t-s)H_0}$ on $\s^2_x(\R^9)$, we have that for all $\phi\in \mathcal{S}_x(\R^9)$
\begin{align}
    &\left|(\phi(x), \int_t^\infty ds \bar{F}_{\tau,l}(H_0,\epsilon_1/10)e^{-i(t-s)H_0}\psi^\pm_{m,l,jl,1}(x,t,s))_{\s^2_x(\R^9)}\right|\nonumber\\
    \leq & \int_t^\infty ds \left| (\phi(x), \bar{F}_{\tau,l}(H_0,\epsilon_1/10)e^{-i(t-s)H_0}\psi^\pm_{m,l,jl,1}(x,t,s))_{\s^2_x(\R^9)}\right|\nonumber\\
    \leq & \int_t^\infty ds \| \langle x_j-x_l \rangle^{-2}\bar{F}_{\tau,l}(H_0,\epsilon_1/10)e^{i(t-s)H_0}\phi \|_{\s^2_x(\R^9)}\| \langle \eta\rangle^4 V_{jl}(\eta)\|_{\s^\infty_\eta(\R^3)} \nonumber\\
    &\times \| \langle x_j-x_l\rangle^{-2} e^{-isH_{a'}}\bar{F}_{\tau}(H_{a'}, c_0)\Omega_{a'}^{*}\psi\|_{\s^2_x(\R^9)}\nonumber\\
    \leq &\opnorm{V} \| \langle x_j-x_l \rangle^{-2}\bar{F}_{\tau,l}(H_0,\epsilon_1/10)e^{i(t-s)H_0}\phi \|_{\s^2_{x,t}(\R^{10})}\nonumber\\
    &\times\left( \int_t^\infty ds \| \langle x_j-x_l\rangle^{-2} e^{-isH_{a'}}\bar{F}_{\tau}(H_{a'}, c_0)\Omega_{a'}^{*}\psi\|_{\s^2_x(\R^9)}^2  \right)^{1/2}\nonumber\\
    \leq & C\opnorm{V}\|\phi\|_{\s^2_x(\R^9)}\times\left( \int_t^\infty ds \| \langle x_j-x_l\rangle^{-2} e^{-isH_{a'}}\bar{F}_{\tau}(H_{a'}, c_0)\Omega_{a'}^{*}\psi\|_{\s^2_x(\R^9)}^2  \right)^{1/2}
\end{align}
for some constant $C=C(c_0, \epsilon_1,\sup\limits_{p\in \R^3}|\tau(p)|)>0$. Therefore, by using the B.L.T. Theorem (Theorem I.7 on page $9$ in \cite{RS4}), we have that for some constant $C=C(c_0, \epsilon_1, \sup\limits_{p\in \R^3}|\tau(p)|)>0$
\begin{align}
    &\| \int_t^\infty ds e^{-i(t-s)H_0}\psi_{m,l,jl,1}^\pm(x,t,s)\|_{\s^2_x(\R^9)}\nonumber\\
    \leq & C\opnorm{V} \left( \int_t^\infty ds \| \langle x_j-x_l\rangle^{-2} e^{-isH_{a'}}\bar{F}_{\tau}(H_{a'}, c_0)\Omega_{a'}^{*}\psi\|_{\s^2_x(\R^9)}^2  \right)^{1/2}\to 0
\end{align}
as $t\to \infty.$ Hence, we get \eqref{9.11end: eq1}. By using \eqref{9.11end: eq1}, we conclude that 
\eq
\int_0^\infty ds \bar{F}_{\tau,l}(H_0,\epsilon_1/10)e^{-i(t-s)H_0}\psi_{m,l,jl,1}^\pm(x,t,s)
\eeq
exists in $\s^2_x(\R^9)$. By using \eqref{9.11end: eq1} again, we obtain that 
\begin{align}
    &\| \chi_R(\tilde{x})P^-_{jk}\bar{F}_{\tau,l}(H_0,\epsilon_1/10)\int_0^t ds e^{-i(t-s)H_0}\psi_{m,l,jl,1}^\pm(x,t,s)\|_{\s^2_{x}(\R^9)}\nonumber\\
    \leq & \| \chi_R(\tilde{x})P^-_{jk}\bar{F}_{\tau,l}(H_0,\epsilon_1/10)\int_0^\infty ds e^{-i(t-s)H_0}\psi_{m,l,jl,1}^\pm(x,t,s)\|_{\s^2_{x}(\R^9)}\nonumber\\
    &+ \| \chi_R(\tilde{x})P^-_{jk}\bar{F}_{\tau,l}(H_0,\epsilon_1/10)\int_t^\infty ds e^{-i(t-s)H_0}\psi_{m,l,jl,1}^\pm(x,t,s)\|_{\s^2_{x}(\R^9)}\nonumber\\
    \to & 0
\end{align}
as $t\to \infty$. Here, we also use the fact that for all $g\in \s^2_x(\R^9)$
\eq
\lim\limits_{t\to \infty} \| P_{jk}^- e^{-itH_0}g\|_{\s^2_x(\R^9)}=0.
\eeq
Therefore, we get \eqref{9.11end: eq2} and finish the proof.
\end{proof}
To demonstrate that $\| \psi_{m,l,jl,1}\|_{\s^2_x(\R^9)}$ can be made small by wisely selecting values for $R$ and $t$, we require the following proposition:
\begin{proposition}\label{Prop: end}For all $M\geq M_0$ $c>0$ and $a=(jk)(l)\in L$, we let
\begin{align}
\psi_{a,M,c}(x,t)=&\langle x_j-x_k\rangle^{-5/2}F(|x_j-x_k|\leq c_1M)F(|x_j-x_l|>M) \bar{F}_{\tau}(H_a,c)e^{-itH}\psi
\end{align}
\begin{align}
   \psi_{a,M,c}^+(x,t)=&\langle x_j-x_k\rangle^{-5/2}F(|x_j-x_k|\leq c_1M)\bar{F}(|x_j-x_l|>M) \tilde{P}_{lj,l}^+\bar{F}_{\tau}(H_a,c)e^{-itH_a}\Omega_{a}^{*}\psi
\end{align}
and
\begin{align}
   \psi_{a,M,c}^-(x,t)=&\langle x_j-x_k\rangle^{-5/2}F(|x_j-x_k|\leq c_1M)\bar{F}(|x_j-x_l|>M) \tilde{P}_{lj,l}^-\bar{F}_{\tau}(H_a,c)e^{-itH_a}\psi,
\end{align}
where $c_1, \frac{1}{M_0}>0$ denote two sufficiently small numbers such that when $|r_1|\geq \frac{M_0}{2}$ and  $|r_2|\leq 2 cM_0$,
\eq
|r_1\pm r_2+q|\geq \frac{1}{100}(|r_1\pm r_2|+|q|),\label{c21}
\eeq
and
\eq
|r_1\pm r_2-q|\geq \frac{1}{10^7}(|r_1\pm r_2|+|q|),\label{c22}
\eeq
if $\hat{r}_1\in \text{supp}(F^{\hat{h}})$ $q\in \R^3- \{0\}$ and $\hat{q}\in \text{supp}(\tilde{F}^{\hat{h}})$ for all $h\in I$. If Assumptions \ref{asp: short}- \ref{asp: subHH} hold, then for all $\epsilon>0$, there exists a $M_1=M_1(M_0, \epsilon,c,c_1, \opnorm{V})\geq M_0$ and $T=T(M_0, \epsilon,c,c_1, \opnorm{V})\geq 1$ such that
\eq
\|\psi_{a,M,c}(x,t)-\psi_{a,M,c}^+(x,t)-\psi_{a,M,c}^-(x,t)\|_{\s^2_x(\R^9)}< \frac{\epsilon}{100000}
\eeq
holds true for all $M\geq M_1, c>0$ and $t\geq T$.    
\end{proposition}
\begin{remark}
    Throughout this note, we regard $c_1$ as a constant.
\end{remark}
We defer the proof of Proposition \ref{Prop: end} to subsection \ref{subsection: prop4.1}. Given $a=(jk)(l)\in L$, we let 
\eq
v_{a}:=\Omega^{a}v_a\Omega^{a,*},
\eeq
\eq
P_{c,\beta}^M(H_{a}):=P_c(H_{a})\bar{F}(|v_a|> \frac{1}{M^{\beta}}),
\eeq
and
\eq
\bar{P}_{c,\beta}^M(H_a):=P_c(H_{a})F(|v_a|\leq \frac{1}{M^{\beta}})
\eeq
for $\beta\in (0,\frac{1}{2})$. Here, $1/100$ standards for some sufficiently small positive number. The proof of Proposition \ref{Prop: end} requires the $\s^2_x(\R^9)$ estimates of, for $t\geq 0$, $a=(jk)(l)\in L$ and $M\geq 1$, 
\eq
\A^{\pm,\beta}_{a,c,c_1,j}(t,M):=\langle x_j-x_k\rangle^{-5/2} F(|x_l-x_j|\geq M) \tilde{P}_{lj,l}^\pm P_{c,\beta}^M(H_{a})\bar{F}_\tau(H_a,c)e^{\pm itH_a}\langle x_l-x_j\rangle^{-2} 
\eeq
and  
\begin{align}
&\bar{\A}^{\pm,\beta}_{a,c,c_1,k}(t,M):=\langle x_j-x_k\rangle^{-5/2}\nonumber\\
&\times F(|x_j-x_k|\leq c_1M) F(|x_l-x_j|\geq M) \tilde{P}_{lj,l}^\pm P_{c,\beta}^M(H_{a})\bar{F}_\tau(H_a,c)e^{\pm itH_a}\langle x_l-x_k\rangle^{-2}. 
\end{align}
To be precise, we have to show that: when $t\geq 2 \sqrt{M}$ and $M\geq 1$, 
\eq
\| \A^{\pm,\beta}_{a,c,c_1,j}(t,M)\|_{\s^2_x(\R^9)\to \s^2_x(\R^9)}\lesssim_{ \epsilon, c,c_1} \frac{1}{\langle t+M\rangle^{3/2(1-\epsilon)}}\label{main:eq}
\eeq
and 
\eq
\| \bar{\A}^{\pm,\beta}_{a,c,c_1,k}(t,M)\|_{\s^2_x(\R^9)\to \s^2_x(\R^9)}\lesssim_{ \epsilon, c,c_1} \frac{1}{\langle t+M\rangle^{3/2(1-\epsilon)}}\label{main:eq2}
\eeq
for some $\beta=\beta(\epsilon)\in (0,\frac{1}{2})$ and all $\epsilon \in (0,1)$ and $c>0$. Here, $c_2$ is some fixed number which satisfies \eqref{c21} and \eqref{c22}.

\begin{proposition}\label{main:prop}Let $c_2, M_0$ be as in Proposition \ref{Prop: end}. If Assumptions \ref{asp: short} to \ref{asp: subHH} hold, then when $t\geq 0$ and $M\geq M_0$, \eqref{main:eq} and \eqref{main:eq2} hold true for all $\epsilon \in (0,1)$ and $c>0$.
    
\end{proposition}
We defer the proof of Proposition \ref{main:prop} to subsection \ref{subsection: Prop2}.

By using Proposition \ref{Prop: end}, we have 
\eq
\| \psi_{m,l,jl,1}^0(x,t,s)-\psi_{m,l,jl,1}^+(x,t,s)-\psi_{m,l,jl,1}^-(x,t,s)\|_{\s^2_x(\R^9)}\lesssim_\epsilon \frac{1}{( R+|t-s|)^{\frac{1}{2}+\delta}}\label{end: eq}
\eeq
for some $\delta>0$. Consequently, there exists a sufficiently large number $R_{free,1}\geq 1$ such that when $R\geq R_{free,1}$, we have that there exists $T_{free,1}=T_{free,1}(R)\geq 1$ such that when $R\geq R_{free,1}$ and $t\geq T_{free,1}$, 
\begin{align}
\| \psi_{m,l,jl,1}(x,t)\|_{\s^2_x(\R^9)}< \frac{\epsilon}{10000}\label{psimljl1}
\end{align}
holds true by using Lemma \ref{Lem: two-body} and Proposition \ref{Prop: end}. Due to \eqref{psimljl1} and \eqref{psimljl2}, we have 
\begin{align}
\| \psi_{m,l,jl}(x,t)\|_{\s^2_x(\R^9)}< \frac{\epsilon}{5000}\label{psimljl}.
\end{align}
Similarly, there exists a sufficiently large number $R_{free,2}\geq 1$ such that when $R\geq R_{free,2}$, there exists $T_{free,2}=T_{free,2}(R)\geq 1$ such that when $R\geq R_{free,2}$ and $t\geq T_{free,2}$, we have
\begin{align}
\| \psi_{m,l,kl}(x,t)\|_{\s^2_x(\R^9)}< \frac{\epsilon}{5000}\label{psimlkl}.
\end{align}
Therefore, due to \eqref{psimljl}, \eqref{psimlkl}, \eqref{psifreejk} and \eqref{psifreel}, we conclude that there exists a sufficiently large number $R^{l,+}_{free}\geq 1$ such that when $R\geq R^{l,+}_{free}$, there exists $T^{l,+}_{free}=T^{l,+}_{free}(R)\geq 1$ such that when $R\geq R^{l,+}_{free}$ and $t\geq T^{l,+}_{free}$, we have
\eq
\| \psi_{R,\epsilon_1,m,l}^+(x,t)-\tilde{g}_a(v_j(P_j)-v_k(P_k))e^{-itH_0}\Omega_{free}^*\psi\|_{\s^2_x(\R^9)}< \frac{\epsilon}{500}\label{psiml+}.
\eeq
Similarly, we obtain that there exists a sufficiently large number $R^{l,-}_{free}\geq 1$ such that when $R\geq R^{l,-}_{free}$, there exists $T^{l,-}_{free}=T^{l,-}_{free}(R)\geq 1$ such that when $R\geq R^{l,-}_{free}$ and $t\geq T^{l,-}_{free}$, we have
\eq
\| \psi_{R,\epsilon_1,m,l}^-(x,t)\|_{\s^2_x(\R^9)}< \frac{\epsilon}{500}\label{psiml-}.
\eeq
Due to \eqref{psiml+} and \eqref{psiml-}, we have that there exists a sufficiently large number $R^{l}_{free}\geq 1$ such that when $R\geq R^{l}_{free}$, there exists $T^{l}_{free}=T^{l}_{free}(R)\geq 1$ such that when $R\geq R^{l}_{free}$ and $t\geq T^{l}_{free}$, we have
\eq
\| \psi_{R,\epsilon_1,m,l}(x,t)-\tilde{g}_a(v_j(P_j)-v_k(P_k))e^{-itH_0}\Omega_{free}^*\psi\|_{\s^2_x(\R^9)}< \frac{\epsilon}{250}\label{psiml}.
\eeq
Hence, we have that there exists a sufficiently large number $R_{free}\geq 1$ such that when $R\geq R_{free}$, there exists $T_{free}=T_{free}(R)\geq 1$ such that when $R\geq R_{free}$ and $t\geq T_{free}$, we have
\begin{align}
    \| \psi_{R,\epsilon_1,m}-e^{-itH_0}\Omega_{free}^*\psi\|_{\s^2_x(\R^9)}\leq &\sum\limits_{l=1}^3  \| \psi_{R,\epsilon_1,m,l}-\tilde{g}_a(v_j(P_j)-v_k(P_k))e^{-itH_0}\Omega_{free}^*\psi\|_{\s^2_x(\R^9)}\nonumber\\
    <& \frac{\epsilon}{60}.\label{psiRm}
\end{align}
By using \eqref{psiRr1}, \eqref{psiRr2} and \eqref{psiRm}, we have that there exists a sufficiently large number $R_{free}'\geq 1$ such that when $R\geq R_{free}'$, there exists $T_{free}'=T_{free}'(R)\geq 1$ such that when $R\geq R_{free}'$ and $t\geq T_{free}'$, we have
\begin{align}
    \| \psi_{R,\epsilon_1}-e^{-itH_0}\Omega_{free}^*\psi\|_{\s^2_x(\R^9)} <& \frac{\epsilon}{10}.
\end{align}
We finish the proof of Lemma \ref{Lem: free} if we have Propositions \ref{Prop: end} and \ref{main:prop}.

\subsection{Proof of Proposition \ref{main:prop}}\label{subsection: Prop2}
The proof of Proposition \ref{main:prop} requires following lemmas. Let $\tau_j, j=1,2,3$, denote the set of all critical points of $\omega_j$. 

\begin{lemma}\label{Lem: main: 2} If Assumptions \ref{asp: short} to \ref{asp: subHH} hold true, then for all $a=(jk)(l)\in L$, $c,c_1>0$ and $\beta\in (0,1/2)$, we have
   \begin{align}
 \| \langle x_j-x_k\rangle^{-2}&P_{c}(H_a)\bar{F}_\tau(H_a,c)e^{itH_a}\bar{F}_{\tau_j}(v_j(P_j),c_1)\bar{F}_\tau(H_0,\frac{c}{10})\langle x_j-x_l\rangle^{-2}\|_{\s^2_x(\R^9)\to \s^2_x(\R^9)}\nonumber\\
\lesssim_{c,c_1,\beta}&\frac{1}{\langle t\rangle^{3/2}} (1+\|\langle \eta\rangle^{4}V_{jk}(\eta)\|_{\s^\infty_\eta(\R^3)}),\label{9.11: eq1}
\end{align}
and
\begin{align}
 \| \langle x_j-x_k\rangle^{-2}&P_{c,\beta}^M(H_a)\bar{F}_\tau(H_a,c)e^{itH_a}\bar{F}_{\tau_j}(v_j(P_j),c_1)\bar{F}_\tau(H_0,\frac{c}{10})\langle x_j-x_l\rangle^{-2}\|_{\s^2_x(\R^9)\to \s^2_x(\R^9)}\nonumber\\
\lesssim_{c,c_1,\beta}&\frac{1}{\langle t\rangle^{3/2}} (1+\|\langle \eta\rangle^{4}V_{jk}(\eta)\|_{\s^\infty_\eta(\R^3)}).\label{9.11: eq2}
\end{align}
 
\end{lemma}
\begin{proof}
Let
\begin{align}
    Q_{M}^a(t):=& \langle x_j-x_k\rangle^{-2}P_c(H_a)\bar{F}_\tau(H_a,c)e^{itH_a}\bar{F}_{\tau_j}(v_j(P_j),c_1)\bar{F}_\tau(H_0,\frac{c}{10})\nonumber\\
    &\times\langle x_j-x_l\rangle^{-2}.
\end{align}
The estimate for $Q_M^a(t)$ requires the estimate for $e^{iu\omega_j(P_j)}(u\in \R)$: Based on the definition of $\tau_j$, by using the method of stationary phase, we have 
    \eq
    \| e^{iu\omega_j(P_j)}\bar{F}_{\tau_j}(v_j(P_j),c_1)\|_{\s^1_{x_j}(\R^3)\to \s^\infty_{x_j}(\R^3)}\lesssim_{c_1} \frac{1}{|u|^{3/2}},\quad \forall|u|\geq 1.\label{free: stationary}
    \eeq
Now let us estimate $Q_M^a(t)$ on $\s^2_x(\R^9)$. Let
\begin{align}
    Q_{M1}^a(t):=&\langle x_j-x_k\rangle^{-2}P_c(H_a)\bar{F}_\tau(H_a,c)e^{itH_0}\bar{F}_{\tau_j}(v_j(P_j),c_1)\bar{F}_\tau(H_0,\frac{c}{10})\nonumber\\
    &\times\langle x_j-x_l\rangle^{-2},
\end{align}
and 
\begin{align}
    Q_{M2}^a(t):=&i\int_0^tds\langle x_j-x_k\rangle^{-2}P_c(H_a)\bar{F}_\tau(H_a,c)e^{isH_a}V_{jk}(x_j-x_k)e^{i(t-s)H_0}\bar{F}_{\tau_j}(v_j(P_j),c_1)\nonumber\\
    &\times\bar{F}_\tau(H_0,\frac{c}{10})\langle x_j-x_l\rangle^{-2}.
\end{align}
Using Duhamel's formula to expand $e^{itH_a}$ in $Q_M^a(t)$, we have 
\begin{align}
    Q_M^a(t)=Q_{M1}^a(t)+Q_{M_2}^a(t).
\end{align}
For $Q_{M1}^a(t)$, by using \eqref{free: stationary} and Assumptions \ref{asp: subH} and \ref{asp: subHH}, we have 
\begin{align}
    \|Q_{M1}^a(t)\|_{\s^2_x(\R^9)\to \s^2(\R^9)}\leq & \| \langle x_j-x_k\rangle^{-2}P_c(H_a)\bar{F}_\tau(H_a,c)\langle x_j-x_k\rangle^{2}\|_{\s^2_x(\R^9)\to \s^2_x(\R^9)}\nonumber\\
    &\times\| \langle x_j-x_k\rangle^{-2}e^{itH_0}\bar{F}_{\tau_j}(v_j(P_j),c_1)\langle x_j-x_l\rangle^{-2}\|_{\s^2_x(\R^9)]\to \s^2_x(\R^9)} \nonumber\\
    &\times\| \langle x_j-x_l\rangle^{2}\bar{F}_\tau(H_0,\frac{c}{10})\langle x_j-x_l\rangle^{-2}\|_{\s^2_x(\R^9)\to \s^2_x(\R^9)}\nonumber\\
    \lesssim_{c_1,c} & \frac{1}{\langle t\rangle^{3/2}},\label{Lem4.3: QM1}
\end{align}
where we also use that when $|t|\geq 1$, by using the unitarity of $e^{it\omega_{l}(P_l)}$ and $e^{it\omega_k(P_k)}$ on $\s^2_x(\R^9)$,
\begin{align}
  &  \| \langle x_j-x_k\rangle^{-2}e^{itH_0}\bar{F}_{\tau_j}(v_j(P_j),c_1)\langle x_j-x_l\rangle^{-2}\|_{\s^2_x(\R^9)]\to \s^2_x(\R^9)}\nonumber\\
  =& \| \langle x_j-x_k\rangle^{-2}e^{it\omega_j(P_j)}\bar{F}_{\tau_j}(v_j(P_j),c_1)\langle x_j-x_l\rangle^{-2}\|_{\s^2_x(\R^9)]\to \s^2_x(\R^9)} \nonumber\\
  \leq & \|\|\langle x_j-x_k\rangle^{-2 }\|_{\s^2_{x_j}(\R^3)}\| e^{it\omega_j(P_j)}\bar{F}_{\tau_j}(v_j(P_j),c_1)\|_{\s^1_{x_j}(\R^3)\to \s^\infty_{x_j}(\R^3)}\nonumber\\
  &\times \|\langle x_j-x_l\rangle^{-2 }\|_{\s^2_{x_j}(\R^3)}\|_{\s^2_{x_k,x_l}(\R^6)\to \s^2_{x_k,x_l}(\R^6)}\nonumber\\
  \lesssim_{ c_1}& \frac{1}{|t|^{3/2}}.
\end{align}
From the estimate for $Q_{M1}^a(t)$ on $\s^2_x(\R^9)$, we have
\eq
\| \langle x_j-x_k \rangle^{-2}e^{itH_0}\bar{F}_\tau (H_0,\frac{c}{10})\bar{F}_{\tau_j}(\omega_j(P_j),c_1)\langle x_j-x_l\rangle^{-2} \|_{\s^2_x(\R^9)\to \s^2_x(\R^9)}\lesssim_{ c_1}\frac{1}{|t|^{3/2}}\label{Lem4.3: eq1}
\eeq
when $|t|\geq 1$. Now we use \eqref{Lem4.3: eq1} to estimate $Q_{M2}^a(t)$ on $\s^2_x(\R^9)$. By employing \eqref{Lem4.3: eq1} and Assumption \ref{asp: subHH}, we have 
\begin{align}
    &\|Q_{M2}^a(t)\|_{\s^2_x(\R^9)\to \s^2_x(\R^9)}\nonumber\\
    \leq & \int_0^tds \| \langle x_j-x_k\rangle^{-2}P_c(H_a)\bar{F}_\tau(H_a,c)e^{isH_a}\langle x_j-x_k\rangle^{-2}\|_{\s^2_x(\R^9)\to \s^2_x(\R^9)}\nonumber\\
    \times& \|\langle \eta\rangle^{4} V_{jk}(\eta)\|_{\s^\infty_\eta(\R^3)}\| \langle x_j-x_k \rangle^{-2}e^{i(t-s)H_0}\bar{F}_\tau (H_0,\frac{c}{10})\bar{F}_{\tau_j}(\omega_j(P_j),c_1)\langle x_j-x_l\rangle^{-2} \|_{\s^2_x(\R^9)\to \s^2_x(\R^9)}\nonumber\\
    \lesssim_{c,c_1}& \int_0^t ds \frac{1}{\langle s\rangle^{3/2}}\times \|\langle \eta\rangle^4 V_{jk}(\eta)\|_{\s^\infty_\eta(\R^3)}\times\frac{1}{\langle t-s\rangle^{3/2}}\nonumber\\
    \lesssim_{c,c_1}& \frac{1}{\langle t\rangle^{3/2}} \|\langle \eta\rangle^{4} V_{jk}(\eta)\|_{\s^\infty_\eta(\R^3)}.\label{Lem4.3: QM2}
\end{align}
According to \eqref{Lem4.3: QM1} and \eqref{Lem4.3: QM2}, we have 
\eq
\|Q_M^a(t)\|_{\s^2_x(\R^9)\to \s^2_x(\R^9)}\lesssim_{c,c_1}\frac{1}{\langle t\rangle^{3/2}} (1+\|\langle \eta\rangle^{4}V_{jk}(\eta)\|_{\s^\infty_\eta(\R^3)}).
\eeq
This completes the proof of \eqref{9.11: eq1}. Similarly, we have \eqref{9.11: eq2}.

\end{proof}
\begin{lemma}\label{Lem: main: 3}If Assumptions \ref{asp: short} to \ref{asp: subHH} hold true, then for all $a=(jk)(l)\in L,$ all $ c,c_1>0$, all $M\geq 1$, all $\alpha\in (0,1)$ and all $t\geq 2M^\alpha$, we have

\begin{align}
    \| &\langle x_j-x_k\rangle^{-5/2}P_c(H_a) e^{i(t^\alpha+M^\alpha) H_a} \bar{F}_\tau(H_a,c)\bar{F}_{\tau^a}(H^a,c_1)(e^{i(t-(t^\alpha+M^\alpha)))H_a}\nonumber\\
    &-e^{i(t-(t^\alpha+M^\alpha)))H_0})\times \bar{F}_\tau (H_0,\frac{c}{10})F_{\tau_j}(v_j(P_j),c_1)\|_{\s^2_x(\R^9)\to \s^2_x(\R^9)}\nonumber\\
    &\lesssim_{  c,c_1}\frac{1}{\langle M^\alpha +t^\alpha\rangle^{3/2} }\|\langle \eta\rangle^4 V_{jk}(\eta)\|_{\s^\infty_\eta(\R^3)}.
\end{align}
    
\end{lemma}
\begin{proof}
    Let
    \begin{align}
        Q_M^a(t):=&\langle x_j-x_k\rangle^{-5/2}P_c(H_a) e^{i(t^\alpha+M^\alpha) H_a} \bar{F}_\tau(H_a,c)\bar{F}_{\tau^a}(H^a,c_1)(e^{i(t-(t^\alpha+M^\alpha)))H_a}\nonumber\\
    &-e^{i(t-(t^\alpha+M^\alpha)))H_0})\times \bar{F}_\tau (H_0,\frac{c}{10})F_{\tau_j}(v_j(P_j),c_1).
    \end{align}
    By using Duhamel's formula to expand $e^{i(t-t^\alpha-M^\alpha)H_a}$ in $Q_{M}^a(t)$, we have 
    \begin{align}
        Q_M^a(t):=&i\int_0^{t-t^\alpha-M^\alpha}ds\langle x_j-x_k\rangle^{-5/2}P_c(H_a) e^{i(t^\alpha+M^\alpha+s) H_a} \bar{F}_\tau(H_a,c)\bar{F}_{\tau^a}(H^a,c_1)\nonumber\\
        &V_{jk}(x_j-x_k)e^{i(t-(t^\alpha+M^\alpha)-s)H_0}\times \bar{F}_\tau (H_0,\frac{c}{10})F_{\tau_j}(v_j(P_j),c_1).
    \end{align}
    By employing Assumptions \ref{asp: subH} and \ref{asp: subHH}, we have that for all $u\in \R$,
    \begin{align}
    &\| \langle x_j-x_k\rangle^{-5/2}e^{iuH_a}\bar{F}_\tau(H_a,c)\bar{F}_{\tau^a}(H^a,c_1) \langle x_j-x_k\rangle^{-5/2}  \|_{\s^2_x(\R^9)\to\s^2_x(\R^9)}\nonumber\\
    \leq &\| \langle x_j-x_k\rangle^{-5/2}\bar{F}_\tau(H_a,c)\langle x_j-x_k\rangle^{5/2}\|_{\s^2_x(\R^9)\to \s^2_x(\R^9)} \nonumber\\
    &\times\| \langle x_j-x_k\rangle^{-5/2} P_c(H_a) e^{iuH_a}\bar{F}_{\tau^a}(H^a,c_1) \langle x_j-x_k\rangle^{-5/2}\|_{\s^2_x(\R^9)\to \s^2_x(\R^9)} \\
    \lesssim_{c,c_1}&\frac{1}{\langle u\rangle^{5/2}},
    \end{align}
    and by using the unitarity of $e^{-iuH_0}(\forall u\in \R)$ on $\s^2_x(\R^9)$,
    we have 
    \begin{align}
        \|Q_M^a(t)\|_{\s^2_x(\R^9)\to \s^2_x(\R^9)}\lesssim_{c,c_1}& \int_0^{t-t^\alpha-M^\alpha }ds \frac{1}{\langle t^\alpha+M^\alpha+s\rangle^{5/2}}\|\langle \eta\rangle^{5/2} V_{jk}(\eta)\|_{\s^\infty_\eta(\R^3)}\nonumber\\
        \lesssim_{c,c_1}& \frac{1}{\langle t^\alpha+M^\alpha \rangle^{3/2}}\|\langle \eta\rangle^{4} V_{jk}(\eta)\|_{\s^\infty_\eta(\R^3)}.
    \end{align}
    This completes the proof.
\end{proof}
\begin{lemma}\label{Lem: main: 4}  For all $a=(jk)(l)\in L,$ all $ c,c_1>0$, all $M\geq 1$, all $\alpha\in (0,1), h\in S^2$($S^2$ denotes the unit sphere in $\R^3$) and all $t\geq 2M^\alpha$,  \begin{align}
\| &F(|x_l-x_j|\geq M) \tilde{P}_{lj,l}^\pm P_c(H_a)e^{\pm i(t^\alpha+M^\alpha)H_a}\bar{F}_\tau(H_a,c)\bar{F}_{\tau^a}(H^a,c_1)\bar{F}_{\tau_l}(v_l(P_l),10^{10}c_1)\nonumber\\
&\times e^{\pm i(t-t^\alpha-M^\alpha)H_0}\bar{F}_\tau(H_0,\frac{c}{10})F_{\tau_j}(v_j(P_j),c_1)\langle x_l-x_j\rangle^{-2} \|_{\s^2_x(\R^9)\to \s^2_x(\R^9)}\nonumber\\
&\lesssim_{ \alpha,c,c_1}\frac{1}{\langle t+M\rangle^2}. 
\end{align}
\end{lemma}
\begin{proof}Let 
\begin{align}
Q_{M}^{a,\pm}(t):=   & F(|x_l-x_j|\geq M) \tilde{P}_{lj,l}^\pm P_c(H_a)e^{\pm i(t^\alpha+M^\alpha)H_a}\bar{F}_\tau(H_a,c)\bar{F}_{\tau^a}(H^a,c_1)\nonumber\\
&\times \bar{F}_{\tau_l}(v_l(P_l),10^{10}c_1)e^{\pm i(t-t^\alpha-M^\alpha)H_0}\bar{F}_\tau(H_0,\frac{c}{10})F_{\tau_j}(v_j(P_j),c_1)\langle x_l-x_j\rangle^{-2}.
\end{align}
Write $e^{\pm i(t^\alpha+M^\alpha)H_a}e^{\pm i(t-t^\alpha-M^\alpha)H_0}$ as 
\eq
e^{\pm i(t^\alpha+M^\alpha)H_a}e^{\pm i(t-t^\alpha-M^\alpha)H_0}=e^{\pm it\omega_l(P_l)}e^{\pm i(t^\alpha+M^\alpha)H^a}e^{\pm i(t-t^\alpha-M^\alpha)H^a_0},
\eeq
where $H^a_0:=\omega_j(P_j)+\omega_k(P_k)$. Let
\begin{align}
    Q_{M1}^{a,\pm}(t):=&F(|x_j-x_l|\geq M)\tilde{P}_{lj,l}^\pm e^{\pm it\omega_l(P_l)}\bar{F}_{\tau_l}(v_l(P_l),10^{10}c_1)\nonumber\\
    &\times F(|x_l-x_j|<\frac{1}{100}(M+10^{10}tc_1))e^{\pm i(t^\alpha+M^\alpha)H^a}\bar{F}_{\tau}(H_a,c)\bar{F}_{\tau^a}(H^a,c_1)e^{\pm i(t-t^\alpha-M^\alpha)H_0^a}\nonumber\\
    &\times \bar{F}_\tau(H_0,\frac{c}{10})F_{\tau_j}(v_j(P_j),c_1)\langle x_j-x_l\rangle^{-2},
\end{align}
    \begin{align}
    Q_{M2}^{a,\pm}(t):=&F(|x_j-x_l|\geq M)\tilde{P}_{lj,l}^\pm e^{\pm it\omega_l(P_l)}\bar{F}_{\tau_l}(v_l(P_l),10^{10}c_1)\nonumber\\
    &\times F(|x_l-x_j|\geq\frac{1}{100}(M+10^{10}tc_1))e^{\pm i(t^\alpha+M^\alpha)H^a}\bar{F}_{\tau}(H_a,c)\bar{F}_{\tau^a}(H^a,c_1)\nonumber\\
    F(|x_l-x_j|&<\frac{1}{1000}(M+10^{10}tc_1))e^{\pm i(t-t^\alpha-M^\alpha)H_0^a} \bar{F}_\tau(H_0,\frac{c}{10})F_{\tau_j}(v_j(P_j),c_1)\langle x_j-x_l\rangle^{-2},
\end{align}
and
\begin{align}
    Q_{M3}^{a,\pm}(t):=&F(|x_j-x_l|\geq M)\tilde{P}_{lj,l}^\pm e^{\pm it\omega_l(P_l)}\bar{F}_{\tau_l}(v_l(P_l),10^{10}c_1)\nonumber\\
    &\times F(|x_l-x_j|\geq\frac{1}{100}(M+10^{10}tc_1))e^{\pm i(t^\alpha+M^\alpha)H^a}\bar{F}_{\tau}(H_a,c)\bar{F}_{\tau^a}(H^a,c_1)\nonumber\\
    F(|x_l-x_j|&\geq\frac{1}{1000}(M+10^{10}tc_1))e^{\pm i(t-t^\alpha-M^\alpha)H_0^a} \bar{F}_\tau(H_0,\frac{c}{10})F_{\tau_j}(v_j(P_j),c_1)\langle x_j-x_l\rangle^{-2}.
\end{align}
We break $Q_M^{a,\pm}(t)$ into three pieces: 
\eq
Q_{M}^{a,\pm}(t)=Q_{M1}^{a,\pm}(t)+Q_{M2}^{a,\pm}(t)+Q_{M3}^{a,\pm}(t).
\eeq
For $Q_{M1}^{a,\pm}(t)$, thanks to the method of non-stationary phase, we have that for all $M\geq 1$ and all $ t\geq 0$,
\begin{align}
A_M^{a,\pm}(t):=&\| F(|x_j-x_l|\geq M)\tilde{P}_{lj,l}^\pm e^{\pm it\omega_l(P_l)}\bar{F}_{\tau_l}(v_l(P_l),10^{10}c_1)\nonumber\\
&\times F(|x_l-x_j|\leq \frac{1}{100}(M+10^{10}tc_1))\|_{\s^2_x(\R^9)\to \s^2_x(\R^9)}\nonumber\\
\lesssim_{c_1,N} & \frac{1}{\langle M+10^{10}tc_1\rangle^N} .\label{Lem4.5: eq1}
\end{align}
By using \eqref{Lem4.5: eq1}, we have 
\begin{align}
    \| Q_{M1}^{a,\pm}(t)\|_{\s^2_x(\R^9)\to \s^2_x(\R^9)}\leq & A_M^{a,\pm}(t)\| e^{\pm i(t^\alpha+M^\alpha)H^a}\bar{F}_\tau(H_a,c)\bar{F}_{\tau^a}(H^a,c_1)e^{\pm i(t-t^\alpha-M^\alpha)H_0^a}\nonumber\\
    &\times \bar{F}_\tau(H_0,\frac{c}{10})F_{\tau_j}(v_j(P_j),c_1)\|_{\s^2_x(\R^9)\to \s^2_x(\R^9)}\|\langle x_l-x_j\rangle^{-2}\|_{\s^2_x(\R^9)\to \s^2_x(\R^9)}\nonumber\\
    \lesssim_{c_1}& \frac{1}{\langle M+10000tc_1\rangle^2}\nonumber\\
    \lesssim_{c_1}& \frac{1}{\langle M+t\rangle^2}.\label{Lem5: QM1}
\end{align}
For $Q_{M2}^{a,\pm}(t)$, since for all $\alpha\in (0,1)$,
{\begin{align}
   B^{a,\pm}_M(t):= &\| F(|x_l-x_j|\geq \frac{1}{100}(M+10^{10}tc_1))e^{\pm i(t^\alpha+M^\alpha)H^a}\bar{F}_\tau(H_a,c)\bar{F}_{\tau^a}(H^a,c_1)\nonumber\\
    &\times F(|x_l-x_j|<\frac{1}{1000}(M+10^{10}tc_1))\|_{\s^2_x(\R^9)\to \s^2_x(\R^9)}\nonumber\\
    \lesssim_{c,c_1,N,\alpha}&\frac{1}{\langle M+10000c_1t\rangle^N }
\end{align}
holds true, we have}
\begin{align}
   & \|Q_{M2}^{a,\pm}(t)\|_{\s^2_x(\R^9)\to \s^2_x(\R^9)}\nonumber\\
   \leq & \| F(|x_j-x_l|\geq M)\tilde{P}_{lj}^\pm e^{\pm it\omega_l(P_l)}\bar{F}_{\tau_l}(v_l(P_l),10^{10}c_1)\|_{\s^2_x(\R^9)\to \s^2_x(\R^9)}\nonumber\\
    &\times B_M^a(t)\|e^{\pm i(t-t^\alpha-M^\alpha)H_0^a} \bar{F}_\tau(H_0,\frac{c}{10})F_{\tau_j}(v_j(P_j),c_1)\langle x_j-x_l\rangle^{-2} \|_{\s^2_x(\R^9)\to \s^2_x(\R^9)}\nonumber\\
    \lesssim_{c,c_1,\alpha}& \frac{1}{\langle M+10^{10}c_1t\rangle^2}\nonumber\\
    \lesssim_{c,c_1,\alpha}&\frac{1}{\langle M+t\rangle^2}.\label{Lem5: QM2}
\end{align}
For $Q_{M3}^{a,\pm}(t)$, by using 
\begin{align}
  C_M^{a,\pm}(t):= & \| F(|x_l-x_j|\geq\frac{1}{1000}(M+10^{10}tc_1))e^{\pm i(t-t^\alpha-M^\alpha)H_0^a} \bar{F}_\tau(H_0,\frac{c}{10})F_{\tau_j}(v_j(P_j),c_1)\nonumber\\
   &\times \langle x_j-x_l\rangle^{-2}\|_{\s^2_x(\R^9)\to \s^2_x(\R^9)}\nonumber\\
   \lesssim_{c,c_1}& \frac{1}{\langle M+10^{10}tc_1\rangle^2},\label{Lem4.5: eq 1}
\end{align}
we have 
\begin{align}
   &\|Q_{M3}^{a,\pm}(t)\|_{\s^2_x(\R^9)\to \s^2_x(\R^9)}\nonumber\\
   \leq & \| F(|x_j-x_l|\geq M)\tilde{P}_{lj,l}^\pm e^{\pm it\omega_l(P_l)}\bar{F}_{\tau_l}(v_l(P_l),10^{10}c_1)\nonumber\\
    &\times F(|x_l-x_j|\geq\frac{1}{100}(M+10^{10}tc_1))e^{\pm i(t^\alpha+M^\alpha)H^a}\bar{F}_{\tau}(H_a,c)\bar{F}_{\tau^a}(H^a,c_1)\|_{\s^2_x(\R^9)\to \s^2_x(\R^9)}\times C_M^a(t) \nonumber\\
    \lesssim_{c,c_1} &\frac{1}{\langle M+10^{10}tc_1\rangle^2}\nonumber\\
    \lesssim_{c,c_1}& \frac{1}{\langle M+t\rangle^2}.\label{Lem5: QM3}
\end{align}
Here, \eqref{Lem4.5: eq 1} follows by using the method of non-stationary phase. Based on \eqref{Lem5: QM1}, \eqref{Lem5: QM2} and \eqref{Lem5: QM3}, we have 
\eq
\|Q_M^{a,\pm}(t)\|_{\s^2_x(\R^9)\to \s^2_x(\R^9)}\lesssim_{c_1,c,\alpha}\frac{1}{\langle M+t\rangle^2}.
\eeq
This completes the proof.

\end{proof}
Now let us prove \eqref{main:eq}. 
\begin{proof}[Proof of Proposition \ref{main:prop}] Given $a=(jk)(l)\in L$, we let
\eq
\A^{\pm}_{a,c,c_1,j}(t,M):=\langle x_j-x_k\rangle^{-5/2} F(|x_l-x_j|\geq M) \tilde{P}_{lj,l}^\pm P_c(H_{a})P_c(H_a)\bar{F}_\tau(H_a,c)e^{\pm itH_a}\langle x_l-x_j\rangle^{-2} 
\eeq
and  
\begin{align}
&\bar{\A}^{\pm}_{a,c,c_1,k}(t,M):=\langle x_j-x_k\rangle^{-5/2}\nonumber\\
&\times F(|x_j-x_k|\leq c_1M) F(|x_l-x_j|\geq M) \tilde{P}_{lj,l}^\pm P_c(H_{a})P_c(H_a)\bar{F}_\tau(H_a,c)e^{\pm itH_a}\langle x_l-x_k\rangle^{-2}. 
\end{align}
Due to the definition of $\tau(v_T)$ for $v_T\in \R^3$, we have 
\eq
e^{itH_a} \bar{F}_\tau(H_a,c)F(|v_a|\leq \frac{1}{\langle M\rangle^{\beta}})=e^{itH_a} \bar{F}_\tau(H_a,c)F(|v_a|\leq \frac{1}{\langle M\rangle^{\beta}})\bar{F}_{\tau_{j}}(v_j(P_j), c_2)
\eeq
for some $c_2>0$. Hence, by using Lemma \ref{Lem: main: 2}, it suffices to show that 
$$
\| \A^{\pm}_{a,c,j}(t,M)\|_{\s^2_x(\R^9)\to \s^2_x(\R^9)}\lesssim_{\epsilon,c,c_1}\frac{1}{\langle t+M\rangle^{\frac{3}{2}(1-\epsilon)}}
$$
and
$$
\| \bar{\A}^{\pm}_{a,c,k}(t,M)\|_{\s^2_x(\R^9)\to \s^2_x(\R^9)}\lesssim_{\epsilon,c,c_1}\frac{1}{\langle t+M\rangle^{\frac{3}{2}(1-\epsilon)}}
$$
for all $\epsilon \in (0,1)$. Let us look at $\A^{+}_{a,c,j}(t,M)$ first. Decompose $\A^{+}_{a,c,j}(t,M)$ into four pieces:
\eq
\A^{+}_{a,c,j}(t,M)=\sum\limits_{m=1}^4\A^{+}_{a,c,m}(t,M),
\eeq
where
\begin{subequations}
\begin{align}
\A^{+}_{a,c,1}(t,M):=&\langle x_j-x_k\rangle^{-5/2}F(|x_l-x_j|\geq M)\tilde{P}_{lj,l}^+P_c(H_a)\bar{F}_\tau(H_a,c)e^{itH_a}\nonumber\\
&\times F_\tau (H_0,\frac{c}{10})\langle x_l-x_j\rangle^{-2},
\end{align}
\begin{align}
\A^{+}_{a,c,2}(t,M):=&\langle x_j-x_k\rangle^{-5/2}F(|x_l-x_j|\geq M)\tilde{P}_{lj,l}^+P_c(H_a)\bar{F}_\tau(H_a,c)e^{itH_a}\nonumber\\
&\times \bar{F}_\tau (H_0,\frac{c}{10})\bar{F}_{\tau_j}(v_j(P_j),c_1)\langle x_l-x_j\rangle^{-2},
\end{align}
\begin{align}
\A^{+}_{a,c,3}(t,M):=&\langle x_j-x_k\rangle^{-5/2}F(|x_l-x_j|\geq M)\tilde{P}_{lj,l}^+P_c(H_a)\bar{F}_\tau(H_a,c)e^{itH_a}\nonumber\\
&\times F_{\tau^a}(H^a,c_1)\bar{F}_\tau (H_0,\frac{c}{10})F_{\tau_j}(v_j(P_j),c_1)\langle x_l-x_j\rangle^{-2},
\end{align}
and 
\begin{align}
\A^{+}_{a,c,4}(t,M):=&\langle x_j-x_k\rangle^{-5/2}F(|x_l-x_j|\geq M)\tilde{P}_{lj,l}^+P_c(H_a)\bar{F}_\tau(H_a,c)e^{itH_a}\nonumber\\
&\times \bar{F}_{\tau^a}(H^a,c_1)\bar{F}_\tau (H_0,\frac{c}{10})F_{\tau_j}(v_j(P_j),c_1)\langle x_l-x_j\rangle^{-2},
\end{align}
with $c_1=c_1(c)>0$ satisfying: ($H^d_0:=\omega_{j'}(P_j)+\omega_{k'}(P_k), d=(j'k')(l')\in L) $
\eq
F_{\tau^d}(H^d_0,10c_1)\bar{F}_\tau(H_0, \frac{c}{10})F_{\tau_{j'}}(v_{j'}(P_{j'}),c_1)=0,\quad\forall d=(j'k')(l')\in L,
\eeq
and 
\eq
F_{\tau_{l'}}(v_{l'}(P_{l'}),10^{10}c_1)\bar{F}_\tau(H_0,\frac{c}{10})F_{\tau_{j'}}(v_{j'}(P_{j'}),c_1)=0,\quad \forall d=(j'k')(l')\in L.
\eeq
\end{subequations}
For $\A^{+}_{a,c,1}(t,M)$, by employing Assumption \ref{asp: subHH}:
\eq
\| \langle x_j-x_k\rangle^{2}\bar{F}_\tau(H_a,c)F_\tau(H_0,\frac{c}{10})\|_{\s^2_x(\R^9)\to \s^2_x(\R^9)}\lesssim_{c}1
\eeq
and Assumption \ref{asp: subH}:
\eq
\| \langle x_j-x_k\rangle^{-2}P_c(H_a)\bar{F}_\tau (H_a,c)e^{itH_a}\langle x_j-x_k\rangle^{-2}\|_{\s^2_x(\R^9)\to \s^2_x(\R^9)}\lesssim_{c}\frac{1}{\langle t\rangle^{3/2}},
\eeq
we have 
\begin{align}
    \| \A^{+}_{a,c,1}(t,M)&\|_{\s^2_x(\R^9)\to \s^2_x(\R^9)}\leq  \| F(|x_l-x_j|\geq M)\tilde{P}_{lj,l}^+\|_{\s^2_x(\R^9)\to \s^2_x(\R^9)}\nonumber\\
    &\times \|\langle x_j-x_k\rangle^{-5/2}P_c(H_a)\sqrt{\bar{F}_\tau (H_a,c)}e^{itH_a}\langle x_j-x_k\rangle^{-2}\|_{\s^2_x(\R^9)\to \s^2_x(\R^9)}\nonumber\\
    &\times \| \langle x_j-x_k\rangle^2 \sqrt{\bar{F}_\tau (H_a,c)}F_\tau (H_0,\frac{c}{10})\|_{\s^2_x(\R^9)\to \s^2_x(\R^9)}\|\langle x_j-x_l\rangle^{-2}\|_{\s^2_x(\R^9)\to \s^2_x(\R^9)}\nonumber\\
    \lesssim_{c}& \frac{1}{\langle t\rangle^{3/2}}.\label{Aac1}
\end{align}
For $\A^{+}_{a,c,2}(t,M)$, by employing Lemma \ref{Lem: main: 2}:
\begin{align}
\| \langle x_j-x_k\rangle^{-5/2}&P_c(H_a)\bar{F}_\tau(H_a,c)e^{itH_a}\bar{F}_{\tau_j}(v_j(P_j),c_1)\bar{F}_\tau(H_0,\frac{c}{10})\langle x_j-x_l\rangle^{-2}\|_{\s^2_x(\R^9)\to \s^2_x(\R^9)}\nonumber\\
\lesssim_{c}&\frac{1}{\langle t\rangle^{3/2}},
\end{align}
we have 
\begin{align}
    &\| \A^{+}_{a,c,2}(t,M)\|_{\s^2_x(\R^9)\to \s^2_x(\R^9)}\leq  \| F(|x_l-x_k|\geq M)\tilde{P}_{lj,l}^+\|_{\s^2_x(\R^9)\to \s^2_x(\R^9)}\nonumber\\
    &\times \| \langle x_j-x_k\rangle^{-5/2}P_c(H_a)\bar{F}_\tau(H_a,c)e^{itH_a}\bar{F}_{\tau_j}(v_j(P_j),c_1)\bar{F}_\tau (H_0,\frac{c}{10})\langle x_j-x_l\rangle^{-2}\|_{\s^2_x(\R^9)\to \s^2_x(\R^9)}\nonumber\\
    & \lesssim_{c}\frac{1}{\langle t\rangle^{3/2}}.\label{Aac2}
\end{align}
For $\A^{+}_{a,c,3}(t,M)$, since 
\eq
F_{\tau^a}(H^a,c_1)\bar{F}_\tau (H_0, \frac{c}{10})F_{\tau_j}(v_j(P_j),c_1)=F_{\tau^a}(H^a,c_1)\bar{F}_{\tau^a}(H^a_0,10c_1)\bar{F}_\tau (H_0, \frac{c}{10})F_{\tau_j}(v_j(P_j),c_1),
\eeq
by employing Assumption \ref{asp: subHH}:
\eq
\| \langle x_j-x_k\rangle^2 F_{\tau^a}(H^a, c_1)\bar{F}_{\tau^a}(H_0^a, 10 c_1)\|_{\s^2_x(\R^9)\to \s^2_x(\R^9)}\lesssim_{ c_1}1,
\eeq
and by using Assumption \ref{asp: subHH}:
\eq
\| \langle x_j-x_k\rangle^{-2}P_c(H_a)e^{itH_a}\bar{F}_\tau(H_a,c)\langle x_j-x_k\rangle^{-2}\|_{\s^2_x(\R^9)\to \s^2_x(\R^9)}\lesssim_{c}\frac{1}{\langle t\rangle^{3/2}},
\eeq
we have 
\begin{align}
    \| \A^{+}_{a,c,3}(t,M)&\|_{\s^2_x(\R^9)\to \s^2_x(\R^9)}\leq  \|F(|x_l-x_j|\geq M)\tilde{P}_{lj,l}^+\|_{\s^2_x(\R^9)\to \s^2_x(\R^9)}\nonumber\\
    &\times \|\langle x_j-x_k\rangle^{-5/2}P_c(H_a)e^{itH_a}\bar{F}_\tau(H_a,c)\langle x_j-x_k\rangle^{-2}\|_{\s^2_x(\R^9)\to \s^2_x(\R^9)}\nonumber\\
    &\times \| \langle x_j-x_k\rangle^2 F_{\tau^a}(H^a,c_1)\bar{F}_{\tau^a}(H_0^a,10c_1)\|_{\s^2_x(\R^9)\to \s^2_x(\R^9)}\nonumber\\
    & \times \| \bar{F}_\tau (H_0,\frac{c}{10})F_{\tau_j}(v_j(P_j),c_1)\langle x_l-x_j\rangle^{-2}\|_{\s^2_x(\R^9)\to \s^2_x(\R^9)}\nonumber\\
    &\lesssim_{c,c_1}\frac{1}{\langle t\rangle^{3/2}}.\label{Aac3}
\end{align}
For $\A^{+}_{a,c,4}(t,M)$, by employing Lemma \ref{Lem: main: 3} (take $\alpha =1-\epsilon$), we have 
\begin{align}
  &\|\langle x_j-x_k\rangle^{-5/2}P_c(H_a) e^{i(t^\alpha+M^\alpha) H_a} \bar{F}_\tau(H_a,c)\bar{F}_{\tau^a}(H^a,c_1)(e^{i(t-(t^\alpha+M^\alpha)))H_a}\nonumber\\
    &-e^{i(t-(t^\alpha+M^\alpha)))H_0})\times \bar{F}_\tau (H_0,\frac{c}{10})F_{\tau_j}(v_j(P_j),c_1)\|_{\s^2_x(\R^9)\to \s^2_x(\R^9)}\nonumber\\
    \lesssim_{ c,c_1}&\frac{1}{\langle M^\alpha +t^\alpha\rangle^{3/2} }\|\langle \eta\rangle^4 V_{jk}(\eta)\|_{\s^\infty_\eta(\R^3)}\nonumber\\
    \lesssim_{c,c_1}& \frac{1}{\langle M +t\rangle^{\frac{3}{2}(1-\epsilon)} }\|\langle \eta\rangle^4 V_{jk}(\eta)\|_{\s^\infty_\eta(\R^3)}
\end{align}
and by using Lemma \ref{Lem: main: 4}:
\begin{align}
\| &F(|x_l-x_j|\geq M) P_{lj,l}^+P_c(H_a)e^{i(t^\alpha+M^\alpha)H_a}\bar{F}_\tau(H_a,c)\bar{F}_{\tau^a}(H^a,c_2)\nonumber\\
&\times e^{i(t-t^\alpha-M^\alpha)H_0}\bar{F}_\tau(H_0,\frac{c}{10})F_{\tau_j}(v_j(P_j),c_1)\langle x_l-x_j\rangle^{-2} \|_{\s^2_x(\R^9)\to \s^2_x(\R^9)}\nonumber\\
&\lesssim_{ \alpha,c,c_1}\frac{1}{\langle t+M\rangle^2}, 
\end{align}
by taking $\alpha =1-\epsilon$, we have 
\begin{align}
\| \A^+_{a,c,4}&(t,M)\|_{\s^2_x(\R^9)\to \s^2_x(\R^9)}\leq 
\| \langle x_j-x_k\rangle^{-5/2}P_c(H_a) e^{i(t^\alpha+M^\alpha) H_a} \bar{F}_\tau(H_a,c)\bar{F}_{\tau^a}(H^a,c_1)\nonumber\\
    &\times(e^{i(t-(t^\alpha+M^\alpha)))H_a}-e^{i(t-(t^\alpha+M^\alpha)))H_0})\times \bar{F}_\tau (H_0,\frac{c}{10})F_{\tau_j}(v_j(P_j),c_1)\|_{\s^2_x(\R^9)\to \s^2_x(\R^9)}\nonumber\\
    &+\| F(|x_l-x_j|\geq M) P^+_{lj,l}P_c(H_a)e^{i(t^\alpha+M^\alpha)H_a}\bar{F}_\tau(H_a,c)\bar{F}_{\tau^a}(H^a,c_2)\nonumber\\
&\times e^{i(t-t^\alpha-M^\alpha)H_0}\bar{F}_\tau(H_0,\frac{c}{10})F_{\tau_j}(v_j(P_j),c_1)\langle x_l-x_j\rangle^{-2} \|_{\s^2_x(\R^9)\to \s^2_x(\R^9)}\nonumber\\
\lesssim_{\epsilon,c,c_1}&\frac{1}{(M+t)^{\frac{3}{2}(1-\epsilon)}}\|\langle \eta\rangle^4 V_{jk}(\eta)\|_{\s^\infty_\eta(\R^3)}+\frac{1}{\langle M+t\rangle^{2} }.\label{Aac4}
\end{align}
When $\sigma \geq 4$, by using \eqref{Aac1}, \eqref{Aac2}, \eqref{Aac3} and \eqref{Aac4}, we get the desired estimate for $\A^+_{a,c,c_1,j}$. Similarly, we can get the desired estimates for $\A^-_{a,c,c_1,j}$ and $\bar{\A}^\pm_{a,c,c_1,k}$. We finish the proof. 
\end{proof}

\subsection{Proof of Proposition \ref{Prop: end}}\label{subsection: prop4.1}
In this part, we prove Proposition \ref{Prop: end}. Let us explain how to prove it first. Fix $a=(jk)(l)\in L$ and choose $\psi\in \s^2_x(\R^9)$ with $\|\psi\|_{\s^2_x(\R^9)}=1$. Recall that 
\begin{align}
\psi_{a,M,c}(x,t)=&\langle x_j-x_k\rangle^{-5/2}F(|x_j-x_k|\leq c_1M)F(|x_j-x_l|>M) \bar{F}_{\tau}(H_a,c)e^{-itH}\psi
\end{align}
\begin{align}
   \psi_{a,M,c}^+(x,t)=&\langle x_j-x_k\rangle^{-5/2}F(|x_j-x_k|\leq c_1M)\bar{F}(|x_j-x_l|>M) \tilde{P}_{lj,l}^+\bar{F}_{\tau}(H_a,c)e^{-itH_a}\Omega_{a}^*\psi
\end{align}
and
\begin{align}
   \psi_{a,M,c}^-(x,t)=&\langle x_j-x_k\rangle^{-5/2}F(|x_j-x_k|\leq c_1M)\bar{F}(|x_j-x_l|>M) \tilde{P}_{lj,l}^-\bar{F}_{\tau}(H_a,c)e^{-itH_a}\psi.
\end{align}
Let 
\begin{align}
\tilde{\psi}_{a,M,c}(x,t)=&\langle x_j-x_k\rangle^{-5/2}F(|x_j-x_k|\leq c_1M)F(|x_j-x_l|>M) P_c^M(H_a)\bar{F}_{\tau}(H_a,c)e^{-itH}\psi
\end{align}
\begin{align}
   &\tilde{\psi}_{a,M,c}^+(x,t)\nonumber\\
   =&\langle x_j-x_k\rangle^{-5/2}F(|x_j-x_k|\leq c_1M)F(|x_j-x_l|>M) \tilde{P}_{lj,l}^+P_c^M(H_a)\bar{F}_{\tau}(H_a,c)e^{-itH_a}\Omega_{a}^*\psi
\end{align}
and
\begin{align}
   &\tilde{\psi}_{a,M,c}^+(x,t)\nonumber\\
   =&\langle x_j-x_k\rangle^{-5/2}F(|x_j-x_k|\leq c_1M)F(|x_j-x_l|>M) \tilde{P}_{lj,l}^-P_c^M(H_a)\bar{F}_{\tau}(H_a,c)e^{-itH_a}\psi.
\end{align}
By using \eqref{asp: subH: eq1} in Assumption \ref{asp: subH} and the unitarity of $e^{-itH}$ on $\s^2_x(\R^9)$, we have that for some constant $C=C(c,\sup\limits_{p\in \R^3}|\tau(p)|)>0$, 
\begin{align}
&\| \psi_{a,M,c}(x,t)-\tilde{\psi}_{a,M,c}(x,t)\|_{\s^2_x(\R^9)}\nonumber\\
\leq& \| \bar{F}(|x_j-x_l|>M)F(|x_j-x_k|\leq c_1M)\|_{\s^2_x(\R^9)\to \s^2_x(\R^9)}\nonumber\\
&\times\| \langle x_j-x_k\rangle^{-5/2}\bar{P}_c^M(H_a)\bar{F}_\tau(H_a,c)\|_{\s^2_x(\R^9)\to \s^2_x(\R^9)} \| e^{-itH}\psi\|_{\s^2_x(\R^9)}\nonumber\\
\leq & \frac{C}{M^{3\beta/2}},\label{psitotpsi1}
\end{align}
\begin{align}
&\| \psi_{a,M,c}^+(x,t)-\tilde{\psi}_{a,M,c}^+(x,t)\|_{\s^2_x(\R^9)}\nonumber\\
\leq& \| \bar{F}(|x_j-x_l|>M)F(|x_j-x_k|\leq c_1M)\tilde{P}_{lj,l}^+\|_{\s^2_x(\R^9)\to \s^2_x(\R^9)} \nonumber\\
&\times \| \langle x_j-x_k\rangle^{-5/2}\bar{P}_c^M(H_a)\bar{F}_\tau(H_a,c)\|_{\s^2_x(\R^9)\to \s^2_x(\R^9)}\| e^{-itH_a}\Omega_{\alpha,l}^{jk,*}\psi\|_{\s^2_x(\R^9)}\nonumber\\
\leq & \frac{C}{M^{3\beta/2}},\label{psitotpsi2}
\end{align}
and 
\begin{align}
&\| \psi_{a,M,c}^-(x,t)-\tilde{\psi}_{a,M,c}^-(x,t)\|_{\s^2_x(\R^9)}\nonumber\\
\leq& \| \bar{F}(|x_j-x_l|>M)F(|x_j-x_k|\leq c_1M)\tilde{P}_{lj,l}^-\|_{\s^2_x(\R^9)\to \s^2_x(\R^9)} \nonumber\\
&\times \| \langle x_j-x_k\rangle^{-5/2}\bar{P}_c^M(H_a)\bar{F}_\tau(H_a,c)\|_{\s^2_x(\R^9)\to \s^2_x(\R^9)} \| e^{-itH_a}\psi\|_{\s^2_x(\R^9)}\nonumber\\
\leq & \frac{C}{M^{3\beta/2}}.\label{psitotpsi3}
\end{align}
When $\beta>1/3$, we have 
\begin{multline}
\| \psi_{a,M,c}(x,t)-\tilde{\psi}_{a,M,c}(x,t)\|_{\s^2_x(\R^9)}+\| \psi_{a,M,c}^+(x,t)-\tilde{\psi}_{a,M,c}^+(x,t)\|_{\s^2_x(\R^9)}\\
+\| \psi_{a,M,c}^-(x,t)-\tilde{\psi}_{a,M,c}^-(x,t)\|_{\s^2_x(\R^9)}\leq \frac{C}{M^{1/2+\delta}}\label{psiaMcr}
\end{multline}
for some $\delta= \frac{3\beta}{2}-\frac{1}{2}>0$. Hence, it suffices to show that for some $\delta>0$,
\eq\label{eq psiM}
\| \tilde{\psi}_{a,M,c}(x,t)-\tilde{\psi}_{a,M,c}^+(x,t)-\tilde{\psi}_{a,M,c}^-(x,t)\|_{\s^2_x(\R^9)}\leq \frac{C}{M^{\frac{1}{2}+\delta}}\|\psi\|_{\s^2_x(\R^9)}
\eeq
for some constant $C(c,\delta)>0$. Now let us explain how to prove \eqref{eq psiM}. Break $\tilde{\psi}_{a,M,c}(x,t)$ into two pieces:
\begin{align}
    \tilde{\psi}_{a,M,c}(x,t)=& \langle x_j-x_k\rangle^{-5/2}F(|x_j-x_k|\leq c_1M)F(|x_j-x_l|>M) \tilde{P}_{lj,l}^+P_c^M(H_a)\bar{F}_{\tau}(H_a,c)e^{-itH}\psi\nonumber\\
    &+\langle x_j-x_k\rangle^{-5/2}F(|x_j-x_k|\leq c_1M)F(|x_j-x_l|>M) \tilde{P}_{lj,l}^-P_c^M(H_a)\bar{F}_{\tau}(H_a,c)e^{-itH}\psi\nonumber\\
    =:&  \tilde{\psi}_{a,M,c,+}(x,t)+\tilde{\psi}_{a,M,c,-}(x,t).
\end{align}
For $\tilde{\psi}_{a,M,c,+}(x,t)$, we use $\tilde{\psi}_{a,M,c}^+(x,t)$ to approximate it. To be precise, using Duhamel's formula with respect to $e^{-itH_a}$, we obtain
\begin{align}
    &\tilde{\psi}_{a,M,c,+}(x,t)-\tilde{\psi}_{a,M,c}^+(x,t)\nonumber\\
    =
    &i\int_t^\infty ds \mathscr{A}_{a,c,c_1,j}^+(s-t,M)\langle x_l-x_j\rangle^2 V_{jl}(x_j-x_l)e^{-isH}\psi\nonumber\\
    &+i\int_t^\infty ds \mathscr{A}_{a,c,c_1,k}^+(s-t,M)\langle x_l-x_k\rangle^2V_{kl}(x_k-x_l)e^{-isH}\psi\nonumber\\
    =:&\tilde{\psi}_{a,M,c,+,r1}(x,t)+\tilde{\psi}_{a,M,c,+,r2}(x,t).
\end{align}
By using Proposition \ref{main:prop} and the unitarity of $e^{-itH}$ on $\s^2_x(\R^9)$, we have that for all $\epsilon\in (0,\frac{1}{2})$ and for some $C=C(\epsilon, c)>0$
\begin{align}\label{psitotpsir1}
    \| \tilde{\psi}_{a,M,c,+,r1}(x,t)\|_{\s^2_x(\R^9)}\leq & \int_t^\infty ds \| \mathscr{A}_{a,c,c_1,j}^+(s-t,M)\|_{\s^2_x(\R^9)\to \s^2_x(\R^9)} \|\langle \eta\rangle^2 V_{jl}(\eta)\|_{\s^2_\eta(\R^3)}\|\psi\|_{\s^2_x(\R^9)}\nonumber\\
   \leq &\int_t^\infty ds \frac{C}{\langle M+s-t\rangle^{\frac{3}{2}(1-\epsilon)}}\|\langle \eta\rangle^2 V_{jl}(\eta)\|_{\s^2_\eta(\R^3)}\nonumber\\
   \leq & \frac{C}{ M^{\frac{1}{2}-\frac{3\epsilon}{2}}}\|\langle \eta\rangle^2 V_{jl}(\eta)\|_{\s^2_\eta(\R^3)},
\end{align}
and 
\begin{align}\label{psitotpsir2}
    \| \tilde{\psi}_{a,M,c,+,r2}(x,t)\|_{\s^2_x(\R^9)}\leq & \int_t^\infty ds \|\mathscr{A}_{a,c,c_1,k}^+(s-t,M)\|_{\s^2_x(\R^9)\to \s^2_x(\R^9)} \|\langle \eta\rangle^2 V_{kl}(\eta)\|_{\s^2_\eta(\R^3)}\|\psi\|_{\s^2_x(\R^9)}\nonumber\\
   \leq &\int_t^\infty ds \frac{C}{\langle M+s-t\rangle^{\frac{3}{2}(1-\epsilon)}}\|\langle \eta\rangle^2 V_{kl}(\eta)\|_{\s^2_\eta(\R^3)}\nonumber\\
   \leq & \frac{C}{M^{\frac{1}{2}-\frac{3\epsilon}{2}}}\|\langle \eta\rangle^2 V_{kl}(\eta)\|_{\s^2_\eta(\R^3)}.
\end{align}
Therefore, due to \eqref{psitotpsi1}-\eqref{psitotpsi3} and \eqref{psitotpsir1}-\eqref{psitotpsir2}, by using Assumption \ref{asp: subH}, we conclude that for all $\epsilon \in (0,1/2)$ and some constant $C=C(\epsilon,c )$,
\eq
\|\psi_{a,M,c}(x,t) \|_{\s^2_x(\R^9)}\leq \frac{C}{M^{\frac{1}{2}-\epsilon} }.\label{psiaMc}
\eeq
Now we are going to estimate $\tilde{\psi}_{a,M,c,+,r1}(x,t)$ and $\tilde{\psi}_{a,M,c,+,r2}(x,t)$ by using \eqref{psiaMc}. For $\tilde{\psi}_{a,M,c,+,r1}(x,t)$, Break it into two pieces:
\begin{align}
    \tilde{\psi}_{a,M,c,+,r1}(x,t)=& i\int_t^\infty ds \mathscr{A}_{a,c,c_1,j}^+(s-t,M)\langle x_l-x_j\rangle^2 \chi(|x_j-x_l|\leq \langle t-s\rangle^{1/2})V_{jl}(x_j-x_l)e^{-isH}\psi\nonumber\\
    &+i\int_t^\infty ds \mathscr{A}_{a,c,c_1,j}^+(s-t,M)\langle x_l-x_j\rangle^2 \chi(|x_j-x_l|> \langle t-s\rangle^{1/2})V_{jl}(x_j-x_l)e^{-isH}\psi\nonumber\\
    =:&\tilde{\psi}_{a,M,c,+,r11}(x,t)+\tilde{\psi}_{a,M,c,+,r12}(x,t).
\end{align}
By using Assumption \ref{asp: subHH} and the unitarity of $e^{-isH}$ on $\s^2_x(\R^9)$, we have 
\begin{align}
   & \| \tilde{\psi}_{a,M,c,+,r11}(x,t)\|_{\s^2_x(\R^9)}\nonumber\\
   \leq & \int_t^\infty ds \| \A_{a,c,c_1,j}^+(s-t,M)\langle x_l-x_j\rangle^{-2}\chi(|x_j-x_l|\leq \langle t-s\rangle^{1/2})\|_{\s^2_x(\R^9)\to \s^2_x(\R^9)}\nonumber\\
    & \times \| \langle \eta\rangle^4 V_{jl}(\eta)\|_{\s^\infty_\eta(\R^3)}\| e^{-isH}\psi\|_{\s^2_x(\R^9)}\nonumber\\
    \leq & \int_t^\infty ds \frac{C}{\langle M+s-t\rangle^2}\| \langle \eta\rangle^4 V_{jl}(\eta)\|_{\s^\infty_\eta(\R^3)}\| \psi\|_{\s^2_x(\R^9)}\nonumber\\
    \leq & \frac{C}{M}\| \langle \eta\rangle^4 V_{jl}(\eta)\|_{\s^\infty_\eta(\R^3)}\label{end: eqr11}
\end{align}
for some constant $C=C(c)>0$. For $\tilde{\psi}_{a,M,c,+,r12}(x,t)$, by using \eqref{psiaMc}, Assumption \ref{asp: subHH} and Lemma \ref{main:prop}, we have that for all $\epsilon \in (0, \frac{1}{5})$ and some constant $C=C(\epsilon, c)>0$,
\begin{align}
    \|\tilde{\psi}_{a,M,c,+,r12}(x,t) \|_{\s^2_x(\R^9)}\leq & \int_t^\infty ds \| \A^+_{a,c,c_1,j}(s-t,M)\|_{\s^2_x(\R^9)\to \s^2_x(\R^9)}\| \langle \eta\rangle^5 V_{jl}(\eta)\|_{\s^\infty_\eta(\R^3)}\nonumber\\
    & \| \psi_{a,M_1,c }(x,t)\vert_{M_1=\frac{\langle t-s\rangle^{1/2}}{10}}\|_{\s^2_x(\R^9)}\nonumber\\
    \leq & \int_t^\infty ds \frac{C}{\langle M+s-t\rangle^{\frac{3}{2}(1-\epsilon)}}\| \langle \eta\rangle^5 V_{jl}(\eta)\|_{\s^\infty_\eta(\R^3)}\times \frac{1}{\langle t-s\rangle^{\frac{1}{2}-\epsilon}}\nonumber\\
    \leq & \frac{C}{M^{1-\frac{5}{2}\epsilon}}\| \langle \eta\rangle^5 V_{jl}(\eta)\|_{\s^\infty_\eta(\R^3)}.\label{end: eqr12}
\end{align}
Therefore, due to \eqref{end: eqr11} and \eqref{end: eqr12}, we have 
\eq
\| \tilde{\psi}_{a,M,c,+,r1}(x,t)\|_{\s^2_x(\R^9)}\leq \frac{C}{M^{\frac{1}{2}+\delta}}\| \langle \eta\rangle^5 V_{jl}(\eta)\|_{\s^\infty_\eta(\R^3)}\label{end: psi+r1}
\eeq
for all $\epsilon\in (0, \frac{1}{2})  $ and some constant $C=C(c,\epsilon)>0$. Similarly, we have 
\eq
\| \tilde{\psi}_{a,M,c,+,r2}(x,t)\|_{\s^2_x(\R^9)}\leq \frac{C}{M^{\frac{1}{2}+\epsilon}}\| \langle \eta\rangle^5 V_{kl}(\eta)\|_{\s^\infty_\eta(\R^3)}\label{end: psi+r2}
\eeq
for all $\epsilon\in (0, \frac{1}{2})  $ and some constant $C=C(c,\epsilon)>0$. Due to \eqref{end: psi+r1} and \eqref{end: psi+r2}, we have that for all $\epsilon\in (0, \frac{1}{2})$ and some constant $C=C(c,\epsilon)>0$
\eq
\| \tilde{\psi}_{a,M,c,+}(x,t)-\tilde{\psi}_{a,M,c}^+(x,t)\|_{\s^2_x(\R^9)}\leq \frac{C}{M^{\frac{1}{2}+\epsilon}}.\label{end: psir1}
\eeq
Similarly, we have 
for all $\epsilon\in (0, \frac{1}{2})$ and some constant $C=C(c,\epsilon)>0$
\eq
\| \tilde{\psi}_{a,M,c,-}(x,t)-\tilde{\psi}_{a,M,c}^-(x,t)\|_{\s^2_x(\R^9)}\leq \frac{C}{M^{\frac{1}{2}+\epsilon}}.\label{end: psir2}
\eeq
Based on \eqref{end: psir1}, \eqref{end: psir2} and \eqref{psiaMcr}, we conclude that 
\eq
\| \tilde{\psi}_{a,M,c,+}(x,t)-\tilde{\psi}_{a,M,c}^+(x,t)\|_{\s^2_x(\R^9)}\leq \frac{C}{M^{\frac{1}{2}+\epsilon}}\label{0end: eq}
\eeq
for all $\epsilon \in (0, \frac{1}{4})$ and some $C=C(\epsilon,c)>0$. We finish the proof.  

\subsection{Proof of Lemma \ref{Lem: free}}
\begin{proof}[Proof of Lemma \ref{Lem: free}] Following the discussion in subsection \ref{sec: outline}, the conclusion of Lemma \ref{Lem: free} by using Propositions \ref{Prop: end} and \ref{main:prop}.
    
\end{proof}
\subsection{Proof of Proposition \ref{Prop: Pmu}}
\begin{proof}[Proof of Proposition \ref{Prop: Pmu}] Following the discussion in subsection \ref{sec:outline}, we conclude that $P_\mu P_{sc}=0$ by using Lemmas \ref{Lem: free} and \ref{mainLem: Pba}.
    
\end{proof}

\noindent{\textbf{Acknowledgements}: }A. S. was partially supported by Simons Foundation Grant number 851844 and NSF grants DMS-2205931. X. W. was partially supported by DMS-2204795, Humboldt Fellowship, NSF CAREER DMS-2044626/DMS-2303146.

\bibliographystyle{amsplain}
\bibliography{bib}

\begin{bibdiv}
\begin{biblist}

\bib{BreeS2019}{article}{
      author={Breeling, M.},
      author={Soffer, A.},
       title={Three-body dispersive scattering},
        date={2019},
     journal={arXiv preprint arXiv:1901.09438},
}

\bib{CS1988}{article}{
      author={Constantin, P.},
      author={Saut, J.-C.},
       title={Local smoothing properties of dispersive equations},
        date={1988},
        ISSN={0894-0347,1088-6834},
     journal={J. Amer. Math. Soc.},
      volume={1},
      number={2},
       pages={413\ndash 439},
         url={https://doi.org/10.2307/1990923},
      review={\MR{928265}},
}

\bib{cycon2009schrodinger}{book}{
      author={Cycon, H.~L.},
      author={Froese, R.~G.},
      author={Kirsch, W.},
      author={Simon, B.},
       title={Schr{\"o}dinger operators: With application to quantum mechanics and global geometry},
   publisher={Springer},
        date={2009},
}

\bib{Der1990}{article}{
      author={Derezi\'{n}ski, J.},
       title={Criteria for the {K}ato smoothness with respect to a dispersive {$N$}-body {S}chr\"{o}dinger operator},
        date={1990},
        ISSN={0022-2488,1089-7658},
     journal={J. Math. Phys.},
      volume={31},
      number={4},
       pages={842\ndash 850},
         url={https://doi.org/10.1063/1.528818},
      review={\MR{1044891}},
}

\bib{D1993}{article}{
      author={Derezi\'{n}ski, J.},
       title={Asymptotic completeness of long-range {$N$}-body quantum systems},
        date={1993},
        ISSN={0003-486X},
     journal={Ann. of Math. (2)},
      volume={138},
      number={2},
       pages={427\ndash 476},
         url={https://doi-org.proxy.libraries.rutgers.edu/10.2307/2946615},
      review={\MR{1240577}},
}

\bib{GJY2004}{article}{
      author={Galtbayar, A.},
      author={Jensen, A.},
      author={Yajima, K.},
       title={Local time-decay of solutions to {S}chr\"{o}dinger equations with time-periodic potentials},
        date={2004},
        ISSN={0022-4715,1572-9613},
     journal={J. Statist. Phys.},
      volume={116},
      number={1-4},
       pages={231\ndash 282},
         url={https://doi.org/10.1023/B:JOSS.0000037203.79298.ec},
      review={\MR{2083143}},
}

\bib{gerard1991mourre}{inproceedings}{
      author={G{\'e}rard, C.},
       title={The mourre estimate for regular dispersive systems},
        date={1991},
   booktitle={Annales de l'ihp physique th{\'e}orique},
      volume={54},
       pages={59\ndash 88},
}

\bib{CL2002}{book}{
      author={G\'{e}rard, C.},
      author={\L~aba, I.},
       title={Multiparticle quantum scattering in constant magnetic fields},
      series={Mathematical Surveys and Monographs},
   publisher={American Mathematical Society, Providence, RI},
        date={2002},
      volume={90},
        ISBN={0-8218-2919-X},
         url={https://doi.org/10.1198/10857110260141265},
      review={\MR{1871447}},
}

\bib{HS1}{article}{
      author={Hunziker, W.},
      author={Sigal, I.~M.},
       title={The quantum {$N$}-body problem},
        date={2000},
        ISSN={0022-2488},
     journal={J. Math. Phys.},
      volume={41},
      number={6},
       pages={3448\ndash 3510},
         url={https://doi-org.proxy.libraries.rutgers.edu/10.1063/1.533319},
      review={\MR{1768629}},
}

\bib{HSS1999}{article}{
      author={Hunziker, W.},
      author={Sigal, I.~M.},
      author={Soffer, A.},
       title={Minimal escape velocities},
        date={1999},
        ISSN={0360-5302,1532-4133},
     journal={Comm. Partial Differential Equations},
      volume={24},
      number={11-12},
       pages={2279\ndash 2295},
         url={https://doi.org/10.1080/03605309908821502},
      review={\MR{1720738}},
}

\bib{JK1979}{article}{
      author={Jensen, A.},
      author={Kato, T.},
       title={Spectral properties of {S}chr\"{o}dinger operators and time-decay of the wave functions},
        date={1979},
        ISSN={0012-7094,1547-7398},
     journal={Duke Math. J.},
      volume={46},
      number={3},
       pages={583\ndash 611},
         url={http://projecteuclid.org/euclid.dmj/1077313577},
      review={\MR{544248}},
}

\bib{LS2021}{article}{
      author={Liu, B.},
      author={Soffer, A.},
       title={The {L}arge {T}imes {A}symptotics of {NLS} type equations},
        date={2021},
     journal={submitted},
}

\bib{Liu-S}{article}{
      author={Liu, B.},
      author={Soffer, A.},
       title={The large time asymptotic solutions of nonlinear schrödinger type equations},
        date={2023},
     journal={Applied Numerical Mathematics},
}

\bib{M1979}{article}{
      author={Mourre, E.},
       title={Link between the geometrical and the spectral transformation approaches in scattering theory},
        date={1979},
        ISSN={0010-3616},
     journal={Comm. Math. Phys.},
      volume={68},
      number={1},
       pages={91\ndash 94},
         url={http://projecteuclid.org.proxy.libraries.rutgers.edu/euclid.cmp/1103905269},
      review={\MR{539739}},
}

\bib{RS4}{book}{
      author={Reed, M.},
      author={Simon, B.},
       title={Methods of modern mathematical physics. {I}. {F}unctional analysis},
   publisher={Academic Press, New York-London},
        date={1972},
      review={\MR{493419}},
}

\bib{RS1978}{book}{
      author={Reed, M.},
      author={Simon, B.},
       title={Methods of modern mathematical physics. {IV}. {A}nalysis of operators},
   publisher={Academic Press [Harcourt Brace Jovanovich, Publishers], New York-London},
        date={1978},
        ISBN={0-12-585004-2},
      review={\MR{493421}},
}

\bib{S1990}{article}{
      author={Sigal, I.~M.},
       title={On long-range scattering},
        date={1990},
        ISSN={0012-7094},
     journal={Duke Math. J.},
      volume={60},
      number={2},
       pages={473\ndash 496},
         url={https://doi-org.proxy.libraries.rutgers.edu/10.1215/S0012-7094-90-06019-3},
      review={\MR{1047762}},
}

\bib{SS1}{article}{
      author={Sigal, I.~M.},
      author={Soffer, A.},
       title={The {$N$}-particle scattering problem: asymptotic completeness for short-range systems},
        date={1987},
        ISSN={0003-486X},
     journal={Ann. of Math. (2)},
      volume={126},
      number={1},
       pages={35\ndash 108},
         url={https://doi-org.proxy.libraries.rutgers.edu/10.2307/1971345},
      review={\MR{898052}},
}

\bib{SS1990}{article}{
      author={Sigal, I.~M.},
      author={Soffer, A.},
       title={Long-range many-body scattering. {A}symptotic clustering for {C}oulomb-type potentials},
        date={1990},
        ISSN={0020-9910},
     journal={Invent. Math.},
      volume={99},
      number={1},
       pages={115\ndash 143},
         url={https://doi-org.proxy.libraries.rutgers.edu/10.1007/BF01234414},
      review={\MR{1029392}},
}

\bib{SS1993}{article}{
      author={Sigal, I.~M.},
      author={Soffer, A.},
       title={Asymptotic completeness for {$N\leq 4$} particle systems with the {C}oulomb-type interactions},
        date={1993},
        ISSN={0012-7094},
     journal={Duke Math. J.},
      volume={71},
      number={1},
       pages={243\ndash 298},
         url={https://doi-org.proxy.libraries.rutgers.edu/10.1215/S0012-7094-93-07110-4},
      review={\MR{1230292}},
}

\bib{SS1994}{article}{
      author={Sigal, I.~M.},
      author={Soffer, A.},
       title={Asymptotic completeness of {$N$}-particle long-range scattering},
        date={1994},
        ISSN={0894-0347},
     journal={J. Amer. Math. Soc.},
      volume={7},
      number={2},
       pages={307\ndash 334},
         url={https://doi-org.proxy.libraries.rutgers.edu/10.2307/2152761},
      review={\MR{1233895}},
}

\bib{SW1998}{article}{
      author={Soffer, A.},
      author={Weinstein, M.~I.},
       title={Time dependent resonance theory},
        date={1998},
        ISSN={1016-443X,1420-8970},
     journal={Geom. Funct. Anal.},
      volume={8},
      number={6},
       pages={1086\ndash 1128},
         url={https://doi.org/10.1007/s000390050124},
      review={\MR{1664792}},
}

\bib{Sof-W5}{article}{
      author={Soffer, A.},
      author={Wu, X.},
       title={{L}ocal {D}ecay {E}stimates},
        date={2022},
     journal={arXiv preprint arXiv:2211.00500},
}

\bib{SW20222}{article}{
      author={Soffer, A.},
      author={Wu, X.},
       title={{O}n {T}he large {T}ime {A}symptotics of {K}lein-{G}ordon type equations with {G}eneral {D}ata-{I}},
        date={2022},
     journal={arXiv preprint arXiv:2204.11261},
}

\bib{SW2}{article}{
      author={Soffer, A.},
      author={Wu, X.},
       title={On the large time asymptotics of klein-gordon type equations with general data-i},
        date={2022},
     journal={arXiv preprint arXiv:2204.11261},
}

\bib{SW20221}{article}{
      author={Soffer, A.},
      author={Wu, X.},
       title={On {T}he large {T}ime {A}symptotics of {S}chr\"{o}dinger type equations with general data},
        date={2022},
     journal={arXiv preprint arXiv:2203.00724},
}

\bib{SW1}{article}{
      author={Soffer, A.},
      author={Wu, X.},
       title={On the large time asymptotics of schr\"{o}dinger type equations with general data},
        date={2022},
     journal={arXiv preprint arXiv:2203.00724},
}

\bib{TH1993}{article}{
      author={Tamura, H.},
       title={Asymptotic completeness for four-body {S}chr\"{o}dinger operators with short-range interactions},
        date={1993},
        ISSN={0034-5318},
     journal={Publ. Res. Inst. Math. Sci.},
      volume={29},
      number={1},
       pages={1\ndash 21},
         url={https://doi-org.proxy.libraries.rutgers.edu/10.2977/prims/1195167541},
      review={\MR{1208027}},
}

\bib{zielinski1997asymptotic}{article}{
      author={Zielinski, L.},
       title={Asymptotic completeness for multiparticle dispersive charge transfer models},
        date={1997},
     journal={journal of functional analysis},
      volume={150},
      number={2},
       pages={453\ndash 470},
}

\end{biblist}
\end{bibdiv}

\vspace{2mm}

\noindent
\textsc{Soffer: Department of Mathematics, Rutgers University, Piscataway, NJ 08854, U.S.A.}\\
{\em email: }\textsf{\bf soffer@math.rutgers.edu}

\medskip\noindent
\textsc{Wu: Department of Mathematics, Texas A \& M University, College Station, TX 77843, U.S.A.}\\
{\em email: }\textsf{\bf xw292@tamu.edu}


\end{document}